\documentclass[11pt,twoside]{article}
\usepackage{fancyhdr}
\usepackage[colorlinks,citecolor=blue,urlcolor=blue,linkcolor=blue,bookmarks=false]{hyperref}
\usepackage{amsfonts,epsfig,graphicx}
\usepackage{afterpage}
\usepackage{amsmath,amssymb,amsthm} 
\usepackage{fullpage}
\usepackage[T1]{fontenc} 
\usepackage{epsf} 
\usepackage{graphics} 
\usepackage{amsfonts,amsmath}
\usepackage[sort,numbers]{natbib} 
\usepackage{psfrag,xspace}
\usepackage{color,etoolbox}

\theoremstyle{plain}

\newtheorem{theo}{Theorem}[section]

\newtheorem{lem}{Lemma}[section]
\newtheorem{prop}{Proposition}[section]
\newtheorem{cor}{Corollary}[section]

\theoremstyle{definition} 

\newtheorem{nota}{Notation}[section]
\newtheorem{de}{Definition}[section]
\newtheorem{exa}{Example}[section]
\newtheorem{as}{Assumption}[section]
\newtheorem{alg}{Algorithm}[section]

\newcommand{\btheo}{\begin{theo}}
\newcommand{\bde}{\begin{de}}
\newcommand{\ble}{\begin{lem}}
\newcommand{\bpr}{\begin{prop}}
\newcommand{\bno}{\begin{nota}}
\newcommand{\bex}{\begin{exa}}
\newcommand{\bcor}{\begin{cor}}
\newcommand{\spro}{\begin{proof}}
\newcommand{\bas}{\begin{as}}
\newcommand{\balg}{\begin{alg}}

\newcommand{\etheo}{\end{theo}}
\newcommand{\ede}{\end{de}}
\newcommand{\ele}{\end{lem}}
\newcommand{\epr}{\end{prop}}
\newcommand{\eno}{\end{nota}}
\newcommand{\eex}{\end{exa}}
\newcommand{\ecor}{\end{cor}}
\newcommand{\fpro}{\end{proof}}
\newcommand{\eas}{\end{as}}
\newcommand{\ealg}{\end{alg}}

\theoremstyle{plain}

\newtheorem{theos}{Theorem}
\newtheorem{props}{Proposition}
\newtheorem{lems}{Lemma}
\newtheorem{cors}{Corollary}

\theoremstyle{definition}
\newtheorem{exas}{Example}
\newtheorem{algs}{Algorithm}
\newtheorem{asss}{Assumption}
\newtheorem{defns}{Definition}
\newcommand{\btheos}{\begin{theos}}
\newcommand{\etheos}{\end{theos}}
\newcommand{\bprops}{\begin{props}}
\newcommand{\eprops}{\end{props}}
\newcommand{\bdes}{\begin{defns}}
\newcommand{\edes}{\end{defns}}
\newcommand{\blems}{\begin{lems}}
\newcommand{\elems}{\end{lems}}
\newcommand{\bcors}{\begin{cors}}
\newcommand{\ecors}{\end{cors}}
\newcommand{\bexs}{\begin{exas}}
\newcommand{\eexs}{\end{exas}}
\newcommand{\balgs}{\begin{algs}}
\newcommand{\ealgs}{\end{algs}}
\newcommand{\bass}{\begin{asss}}
\newcommand{\eass}{\end{asss}}

\newlength{\widebarargwidth}
\newlength{\widebarargheight}
\newlength{\widebarargdepth}

\newcommand{\matsnorm}[2]{|\!|\!| #1 | \! | \!|_{{#2}}}

\newcommand{\opnorm}[1]{\ensuremath{\matsnorm{#1}{\text{\footnotesize{\mbox{op}}}}}}
\newenvironment{carlist}
 {\begin{list}{$\bullet$}
 {\setlength{\topsep}{0in} \setlength{\partopsep}{0in}
  \setlength{\parsep}{0in} \setlength{\itemsep}{\parskip}
  \setlength{\leftmargin}{0.07in} \setlength{\rightmargin}{0.08in}
  \setlength{\listparindent}{0in} \setlength{\labelwidth}{0.08in}
  \setlength{\labelsep}{0.1in} \setlength{\itemindent}{0in}}}
 {\end{list}}

\newcommand{\bcar}{\begin{carlist}}
\newcommand{\ecar}{\end{carlist}}

\newcommand{\Prob}{\ensuremath{\mathbb{P}}}

\makeatletter
\long\def\@makecaption#1#2{
        \vskip 0.8ex
        \setbox\@tempboxa\hbox{\small {\bf #1:} #2}
        \parindent 1.5em  
        \dimen0=\hsize
        \advance\dimen0 by -3em
        \ifdim \wd\@tempboxa >\dimen0
                \hbox to \hsize{
                        \parindent 0em
                        \hfil 
                        \parbox{\dimen0}{\def\baselinestretch{0.96}\small
                                {\bf #1.} #2
                                } 
                        \hfil}
        \else \hbox to \hsize{\hfil \box\@tempboxa \hfil}
        \fi
        }
\makeatother

\long\def\comment#1{}

\makeatletter
\def\@cite#1#2{[\if@tempswa #2 \fi #1]}
\makeatother

\makeatletter
\long\def\barenote#1{
    \insert\footins{\footnotesize
    \interlinepenalty\interfootnotelinepenalty 
    \splittopskip\footnotesep
    \splitmaxdepth \dp\strutbox \floatingpenalty \@MM
    \hsize\columnwidth \@parboxrestore
    {\rule{\z@}{\footnotesep}\ignorespaces
      #1\strut}}}
\makeatother

\newcommand{\bit}{\begin{itemize}}
\newcommand{\eit}{\end{itemize}}
\newcommand{\ben}{\begin{enumerate}}
\newcommand{\een}{\end{enumerate}}

\newcommand{\bear}{\begin{eqnarray}}
\newcommand{\eear}{\end{eqnarray}}

\newcommand{\order}{{\mathcal{O}}}

\newcommand{\inprod}[2]{\ensuremath{\langle #1 , \, #2 \rangle}}

\newcommand{\obs}{\ensuremath{y}}

\newcommand{\Exs}{\ensuremath{{\mathbb{E}}}}

\newcommand{\beq}{\begin{quotation}}
\newcommand{\enq}{\end{quotation}}

\newcommand{\estart}{\begin{equation}}
\newcommand{\eend}{\end{equation}}

\newcommand{\widgraph}[2]{\includegraphics[keepaspectratio,width=#1]{#2}}

\newcommand{\defn}{\ensuremath{:  =}}

\newcommand{\bec}{\begin{center}}
\newcommand{\enc}{\end{center}}

\newcommand{\beit}{\begin{itemize}}
\newcommand{\enit}{\end{itemize}}

\newcommand{\been}{\begin{enumerate}}
\newcommand{\enen}{\end{enumerate}}

\newcommand{\comsl}{\begin{slide}}
\newcommand{\comspor}{\begin{slide*}}

\newcommand{\comsld}[2]{\begin{slide}[#1,#2]}
\newcommand{\comspord}[2]{\begin{slide*}[#1,#2]}

\newcommand{\mendsl}{\end{slide}}
\newcommand{\mendspo}{\end{slide*}}

\newcommand{\real}{\ensuremath{{\mathbb{R}}}}

\DeclareMathOperator{\sign}{sign}

\newcommand{\Ycar}{\ensuremath{V}} 
\newcommand{\Dmat}{\ensuremath{D}}
\newcommand{\Event}{\ensuremath{\mathcal{E}}}

\newcommand{\mismain}{V}

\newcommand{\HACKPAR}{\ensuremath{a}}
\newcommand{\BRAZIL}{\ensuremath{\xi}}
\newcommand{\Wop}{\ensuremath{W}}
\newcommand{\What}{\ensuremath{\widehat{\Wop}}}

\newcommand{\Fil}{\ensuremath{\mathcal{F}}}
\newcommand{\QFUNTHREE}[3]{\ensuremath{Q_{#3}(#1|#2)}}

\newcommand{\vecone}{1}

\newcommand{\memsamp}[1]{\ensuremath{M_{n}(#1)}}

\newcommand{\mgradsamp}[1]{\ensuremath{G_{n}(#1)}}

\newcommand{\red}[1]{\textcolor{red}{#1}}

\newtoggle{figs}
\newtoggle{cmts}

\toggletrue{figs}
\togglefalse{cmts}

\iftoggle{cmts}{
\newcommand{\mjwcomment}[1]{{\bf{{\red{{MJW --- #1}}}}}}
\newcommand{\mjwcommenta}[1]{{\bf{{\color{yellow}{{MJW --- #1}}}}}}
\newcommand{\sbcomment}[1]{{\bf{{\red{{SB --- #1}}}}}}
\newcommand{\sbcommenta}[1]{{\bf{{\color{yellow}{{SB --- #1}}}}}}}
{\newcommand{\mjwcomment}[1]{}
\newcommand{\mjwcommenta}[1]{}
\newcommand{\sbcomment}[1]{}
\newcommand{\sbcommenta}[1]{}}

\newcommand{\fdens}[1]{\ensuremath{f_{#1}}}
\newcommand{\gdens}[1]{\ensuremath{g_{#1}}}
\newcommand{\ParSet}{\ensuremath{\Omega}}
\newcommand{\ParSpace}{\ParSet}

\newcommand{\pmsstar}{\ensuremath{\pms^*}}
\newcommand{\pmstar}{\ensuremath{\pmsstar}}
\newcommand{\stepsize}{\alpha} 

\newcommand{\MOGSNR}{\ensuremath{\eta}}
\newcommand{\MORSNR}{\ensuremath{\eta}}
\newcommand{\MYSNR}{\ensuremath{\mbox{SNR}}}
\newcommand{\tstar}{\ensuremath{t^*}}

\newcommand{\mprob}{\ensuremath{\mathbb{P}}}

\newcommand{\EMupdate}{\ensuremath{M}}
\newcommand{\PopEMupdate}{\ensuremath{\EMupdate}}

\newcommand{\specq}{\ensuremath{q}}
\newcommand{\HACKQ}{\ensuremath{\specq}}
\newcommand{\Balltwor}{\ensuremath{\mathbb{B}_2(r;\theta^*)}}
\newcommand{\BalltworONE}[1]{\ensuremath{\mathbb{B}_2(#1;\theta^*)}}
\newcommand{\BalltwoTWO}[2]{\ensuremath{\mathbb{B}_2(#1;#2)}}

\newcommand{\Mem}{\ensuremath{M}}
\newcommand{\MemSamp}{\ensuremath{\Mem_\numobs}} 

\newcommand{\Tcrit}{\ensuremath{T}}

\newcommand{\thetasamp}{\ensuremath{\pms}}

\newcommand{\memsubsamp}[1]{\ensuremath{M_{\numobs/T}(#1)}}

\newcommand{\rateemunif}{\varepsilon_M^{\mathrm{unif}}}
\newcommand{\rategradunif}{\varepsilon_G^{\mathrm{unif}}}

\newcommand{\deltaconstmor}{32}
\newcommand{\snrconstmor}{\eta}

\newcommand{\PopEM}{\ensuremath{M}}
\newcommand{\HACKRATEEM}{\ensuremath{\rateem \Big(\frac{\numobs}{\Tfinal},
    \frac{\delta}{\Tfinal} \Big)}}
\newcommand{\pmstil}{\ensuremath{\tilde{\pms}}}
\newcommand{\OrGrad}{\ensuremath{T}}
\newcommand{\PopGradEM}{\ensuremath{G}}
\newcommand{\condens}{k}
\newcommand{\pmsone}{\ensuremath{\pms}}
\newcommand{\pmstwo}{\ensuremath{\pms'}}
\newcommand{\HiddenSpace}{\ensuremath{\mathcal{Z}}}
\newcommand{\ObsSpace}{\ensuremath{\mathcal{Y}}}
\newcommand{\Hidden}{\ensuremath{Z}}
\newcommand{\Obs}{\ensuremath{Y}}
\newcommand{\hidden}{\ensuremath{z}}
\newcommand{\Qfun}{\ensuremath{Q}}
\newcommand{\Qhat}{\ensuremath{\Qfun_\numobs}}
\newcommand{\QFUN}[2]{\ensuremath{\Qfun(#1\vert#2)}}
\newcommand{\QFUNHAT}[2]{\ensuremath{\Qhat(#1\vert#2)}}

\newcommand{\Grem}{\ensuremath{G}}
\newcommand{\ncon}{\ensuremath{\xi}}
\newcommand{\GremSamp}{\ensuremath{\Grem_\numobs}}
\newcommand{\NORMAL}{\ensuremath{\mathcal{N}}}
\newcommand{\usedim}{\ensuremath{d}}
\newcommand{\missing}{\ensuremath{\rho}}
\newcommand{\xtil}{\tilde{x}}
\newcommand{\HACKRATEGRAD}{\ensuremath{\rategrad\Big(\frac{\numobs}{\Tfinal},
    \frac{\delta}{\Tfinal} \Big)}}
\newcommand{\Ytil}{\ensuremath{\widetilde{Y}}}
\newcommand{\DelPrime}{\ensuremath{\widetilde{\Delta}}}
\newcommand{\thetastar}{\ensuremath{\theta^*}}
\newcommand{\sigmaSGD}{\ensuremath{\sigma_{\mathrm{G}}}}

\newcommand{\Psit}{\Psi}
\newcommand{\newpsi}{\Psi}
\newcommand{\mm}{\mu_\theta}
\newcommand{\SIVA}{\ensuremath{\delta}}

\newcommand{\regnoise}{\ensuremath{v}}

\newcommand{\pmsit}[1]{\ensuremath{\pms^{#1}}}

\newcommand{\Tfinal}{\ensuremath{T}}

\newcommand{\thetatil}{\ensuremath{\widetilde{\theta}}}

\newcommand{\MYEXP}[1]{\ensuremath{e^{#1}}}

\newcommand{\sivacon}{\ensuremath{c}}

\newcommand{\thetahat}{\ensuremath{\widehat{\theta}}}

\newcommand{\DelTil}{\ensuremath{\widetilde{\Delta}}} 
\newcommand{\SpecSet}{\ensuremath{\mathbb{A}}}

\newcommand{\Sphere}[1]{\ensuremath{\mathbb{S}^{#1}}}

\newcommand{\rade}[1]{\ensuremath{\varepsilon_{#1}}}

\newcommand{\LASSI}{\ensuremath{c}}

\newcommand{\TAUCON}{\ensuremath{c_\tau}}

\newcommand{\plaincon}{\ensuremath{c}}

\newcommand{\qvec}{\ensuremath{v}}
\newcommand{\qvechat}{\ensuremath{\widehat{\qvec}}}
\newcommand{\Sig}{\ensuremath{\Sigma}}
\newcommand{\SigHat}{\ensuremath{\widehat{\Sig}}}

\newcommand{\CovMat}{\ensuremath{\Sigma}}

 \newcommand{\myinter}{\omega} 
\newcommand{\ccon}{\ensuremath{\xi_1}}
\newcommand{\ctwo}{\ensuremath{\xi_2}}
\newcommand{\mytemp}{\pi}

\newcommand{\fsavg}{\ensuremath{\frac{1}{n} \sum_{i=1}^n}}
\newcommand{\fsavgtwo}[2]{\ensuremath{\frac{#1}{#2 n} \sum_{i=1}^n}}
\newcommand{\numobs}{\ensuremath{n}}

\newcommand{\beqs}{\begin{equation*}}
\newcommand{\eeqs}{\end{equation*}}

\newcommand{\contractparam}{\ensuremath{\kappa}}
\newcommand{\smoothparam}{\ensuremath{\mu}}
\newcommand{\strongparam}{\ensuremath{\lambda}}

\newcommand{\mism}{\ensuremath{U}}

\DeclareSymbolFont{extraup}{U}{zavm}{m}{n}
\DeclareMathSymbol{\vardiamond}{\mathalpha}{extraup}{87}

\newcommand{\xobs}{\ensuremath{x_{\mathrm{obs}}}}
\newcommand{\xmis}{\ensuremath{x_{\mathrm{mis}}}}
\newcommand{\bobs}{\ensuremath{\theta_{\mathrm{obs}}}}
\newcommand{\bmis}{\ensuremath{\theta_{\mathrm{mis}}}}
\newcommand{\zobs}{\ensuremath{z_{\mathrm{obs}}}}

\newcommand{\xobsi}{\ensuremath{x_{\mathrm{obs},i}}}

\newcommand{\Xobs}{\ensuremath{X_{\mathrm{obs}}}}

\newcommand{\sobs}{\ensuremath{| \xobs |}}

\newcommand{\eone}{\ensuremath{\mathcal{E}_1}}
\newcommand{\etwo}{\ensuremath{\mathcal{E}_2}}
\newcommand{\ethree}{\ensuremath{\mathcal{E}_3}}
\newcommand{\efour}{\ensuremath{\mathcal{E}_4}}
\newcommand{\efive}{\ensuremath{\mathcal{E}_5}}
\newcommand{\esix}{\ensuremath{\mathcal{E}_6}}

\newcommand{\rateem}{\ensuremath{\varepsilon_M}}
\newcommand{\rategrad}{\ensuremath{\varepsilon_G}}

\newcommand{\pms}{\ensuremath{\theta}}
\newcommand{\pmss}{\ensuremath{\Omega}}

\newcommand{\fp}{\ensuremath{{\theta^{\ast}}}}

\newcommand{\gopt}{\ensuremath{\theta^*}}
\newcommand{\gvec}{\ensuremath{\theta}}

\newcommand{\rvec}{\ensuremath{\theta}}
\newcommand{\mopt}{\ensuremath{\theta^*}}
\newcommand{\mvec}{\ensuremath{\theta}}

\newcommand{\gsonenospace}{GS}
\newcommand{\gsone}{GS$\;$}
\newcommand{\gstwonospace}{FOS}
\newcommand{\gstwo}{FOS$\;$}

\newcommand{\gsonetext}{{{Gradient Stability}}}
\newcommand{\gstwotext}{{{First-order Stability}}}

\newcommand{\gsoneparam}{\gamma}
\newcommand{\gstwoparam}{\gamma}

\newcommand{\qfun}[2]{\ensuremath{Q (#1|#2)}}

\newcommand{\qfunone}[1]{\ensuremath{Q(\cdot|#1)}}

\newcommand{\mem}[1]{\ensuremath{M(#1)}}
\newcommand{\mgrad}[1]{\ensuremath{G(#1)}}
\newcommand{\memzero}{\ensuremath{M}}
\newcommand{\mgradzero}{\ensuremath{G}}

\newcommand{\mga}[1]{\ensuremath{T(#1)}}
\newcommand{\mgazero}{\ensuremath{T}}

\newcommand{\norm}[1]{\|#1 \|_2}

\begin{document}

\begin{center} {\LARGE{\bf{Statistical guarantees for the EM algorithm: \\ From population to sample-based analysis
 }}} \\

\vspace*{.3in}

{\large{
\begin{tabular}{ccccc}
Sivaraman Balakrishnan$^\dagger$ && Martin J. Wainwright$^{\dagger,\ast}$ && Bin Yu$^{\dagger,\ast}$ \\
\end{tabular}

\vspace*{.1in}

\begin{tabular}{ccc}
Department of Statistics$^{\dagger}$ 
& & Department of Electrical Engineering and Computer Sciences$^{\ast}$
\end{tabular}
\begin{tabular}{c}
University of California, Berkeley \\
Berkeley, CA  94720
\end{tabular}

\vspace*{.2in}

\begin{tabular}{c}
{\texttt{$\{$sbalakri,wainwrig,binyu$\}$@berkeley.edu}}
\end{tabular}
}}

\vspace*{.2in}

\today
\vspace*{.2in}
\begin{abstract}
We develop a general framework for proving rigorous guarantees on the
performance of the EM algorithm and a variant known as gradient
EM. 
Our analysis is divided into two parts: a treatment of these
algorithms at the population level (in the limit of infinite data),
followed by results that apply to updates based on a finite set of
samples.  First, we characterize the domain of attraction of any
global maximizer of the population likelihood.  This characterization
is based on a novel view of the EM updates as a perturbed form of
likelihood ascent, or in parallel, of the gradient EM updates as a
perturbed form of standard gradient ascent.  Leveraging this
characterization, we then provide non-asymptotic guarantees on the EM
and gradient EM algorithms when applied to a finite set of samples.
We develop consequences of our general theory for three canonical
examples of incomplete-data problems: mixture of Gaussians, mixture of
regressions, and linear regression with covariates missing completely
at random.  In each case, our theory guarantees that with a suitable
initialization, a relatively small number of EM (or gradient EM) steps
will yield (with high probability) an estimate that is within
statistical error of the MLE.  We provide simulations to confirm this
theoretically predicted behavior.
\end{abstract}
\end{center}

\section{Introduction}
\label{SecIntroduction}
Data problems with missing values, corruptions, and latent variables
are common in practice. From a computational standpoint, computing the
maximum likelihood estimate (MLE) in such incomplete data problems can
be quite complex.  To a certain extent, these concerns have been
assuaged by the development of the expectation-maximization (EM)
algorithm, along with growth in computational resources.  The EM
algorithm is widely applied to incomplete data problems, and there is
now a very rich literature on its behavior
(e.g.,~\cite{dlr,hardem,mengecm,ecme,aem1,aem2,neal99,sem1,sem2,sem3,sem4}).
However, a major issue is that in most models, although the MLE is
known to have good statistical properties, the EM algorithm is only
guaranteed to return a local optimum.  The goal of this paper is to
address this potential gap between statistical and computational
guarantees in application of the EM algorithm.

The EM algorithm has a lengthy and rich history.  Various algorithms
of the EM-type were analyzed in early work
(e.g.,\cite{hartley,healy,rubin74,beale75,orchard1972,sundberg,baum1970}),
before \citet{dlr} introduced the EM algorithm in its modern general
form.  Among other results, they established its well-known
monotonicity properties.  The subsequent work of \citet{wu1983}
established some of the most general convergence results known for the
EM algorithm; see also the more recent papers~\cite{tseng,hero}.
Together with other results, \citet{wu1983} showed that if the likelihood is
unimodal and certain regularity conditions hold, then the EM algorithm
converges to the unique global optimum.  However, in most interesting
cases of the EM algorithm, the likelihood function is multi-modal, in
which case the behavior of the EM algorithm remains a little more
mysterious.  Indeed, despite its popularity and widespread practical
effectiveness, the EM algorithm is often considered a ``sensible
heuristic'' with little or no theoretical backing.

One interesting observation with the EM algorithm is given a
``suitable'' initialization, it often converges to a statistically
useful estimate.  For instance, in application to a mixture of
regressions problem (see Section~\ref{SecExaMOR} for more details),
Chaganty and Liang~\cite{spectral_mor} empirically demonstrate good performance
for a two-stage estimator, in which the method of methods is used as
an initialization, and then the EM algorithm is applied to refine
this initial estimator.  Although encouraging, this type of behavior
is not well understood in a quantitative sense, especially how EM
fixed points reached by this type of two-stage estimator 
are related to the global maximizers of the population
likelihood.  The goal of this paper is to address this question, and
to develop some general tools for characterizing fixed points of the
suitably initialized sample-based EM algorithm, and their relation to 
maximum likelihood estimates.

Some two-stage estimators have recently been analyzed in work on 
alternating minimization algorithms (see e.g. \cite{kmo10,jns13,hardem_mor,altmin_pr}) 
which show that at least in certain special cases optimization methods 
can be locally effective despite non-convexity. Most directly related 
to our work is the paper of \citet{hardem_mor} which considers a 
special (degenerate) noiseless case of the EM algorithm for the 
mixtures of regressions problem. Results for the noisy mixtures of 
regressions problem follow from our general treatment of the 
EM algorithm (see Section~\ref{SecExaMOR}).

In some settings, performing an exact M-step is computationally
burdensome, in which case a natural alternative is some form of
generalized EM updates.  In such an algorithm, instead of performing
an exact maximization, we simply choose a parameter value that does
not decrease the likelihood.  In addition to the standard EM updates,
we also analyze a particular case of such an algorithm, known as
gradient EM, based on a taking a single gradient step per iteration.

Our main results concern the population EM and gradient EM algorithms
and their finite-sample counterparts. 
Our first set of results (Theorems~\ref{ThmEM} and~\ref{ThmGradEM})
give conditions under which the population algorithms are contractive
to the MLE, in a ball around the MLE.
These results are completely deterministic. 
This population-level analysis is based on viewing these algorithms
as perturbed versions of certain ``oracle'' algorithms which 
are known to be contractive around the MLE. 
Our second set of results (Theorem~\ref{ThmSampleEM},
Theorem~\ref{ThmSampleGradEM} and Theorem~\ref{ThmStochasticEM})
concern the sample-based EM and gradient EM algorithms which 
approximate the 
population-based algorithms 
using a subset of samples at each step. We give conditions under
which these sample operators 
converge to an $\varepsilon$-ball around the population MLE.
These results involve probabilistic bounds on the deviations
between the iterates of the population and sample-based algorithms.

The remainder of this paper is organized as follows.
Section~\ref{SecBackground} provides an introduction to the EM and
gradient EM algorithms, as well as a description of the three examples
treated in detail in this paper---namely, Gaussian mixture models
(Section~\ref{SecExaGMM}), mixture of regressions
(Section~\ref{SecExaMOR}), and regression with missing covariates
(Section~\ref{SecExaMissing}). Section~\ref{SecMain} is devoted to our
general convergence results on both the EM and gradient EM algorithms. 
In Section~\ref{SecConsequences}, we revisit the three model classes
previously introduced, and illustrate the use of our general theory by
deriving some concrete corollaries. In concrete examples our theory
gives a characterization of the quality of initialization needed and 
the rate of convergence of the EM and gradient EM algorithms. 
We complement these theoretical
results with simulations that confirm various aspects of the
theoretical predictions.  In order to promote readability, we defer
the more technical aspects of proofs to the appendices.


\section{Background and model examples}
\label{SecBackground}

We begin with basic background on the EM algorithm and its variants,
along with a number of specific models that we revisit later in the
paper.


\subsection{EM algorithm and its relatives}
Let $\Obs$ and $\Hidden$ be random variables taking values in the
sample spaces $\mathcal{\ObsSpace}$ and $\mathcal{\HiddenSpace}$,
respectively.  Suppose that the pair $(\Obs,\Hidden)$ has a joint
density function $\fdens{\pmsstar}$ that belongs to some parameterized
family $\{\fdens{\theta} \, \mid \, \pms \in \ParSet \}$, for a non-empty compact
convex set $\ParSet$.  Rather than
observing the complete data $(\Obs, \Hidden)$, we observe
only component $\Obs$.  Thus, the component $\Hidden$ corresponds to
the missing or latent structure in the data.  

Our goal is to obtain an estimate of the unknown parameter $\pmstar$
via maximum likelihood---namely, to compute some $\thetahat \in
\ParSet$ maximizing the function $\pms \mapsto \gdens{\pms}(\obs)$,
where
\begin{align}
\gdens{\pms}(\obs) & = \int_{\HiddenSpace} \fdens{\pms}(\obs,
\hidden) d \hidden
\end{align}
is the density function of the observed variable $Y$. Throughout this
paper, we assume that $\pmstar$ is a maximizer of the population
likelihood, but not that $\pmstar$ is a \emph{unique} maximizer.
Uniqueness is often violated in mixture models for which parameters
are typically only identifiable up to permutation.  In the examples
that we consider, this non-identifiability will be resolved by
appropriate initialization conditions.

In many settings, it can be difficult or computationally expensive to
evaluate the log likelihood of the observed data, but relatively easy
to compute the log likelihood $\log f_\theta(\obs, \hidden)$ of both
the latent and observed variables.  The EM algorithm is well-suited to
such settings.  For each $\pmsone \in \ParSpace$, let
$\condens_{\pmsone}(\hidden \mid \obs)$ denote the conditional density
of $\hidden$ given $\obs$.  A straightforward application of Jensen's
inequality then shows that the log likelihood at $\pmstwo
\in \ParSpace$ can be lower bounded as
\begin{align}
\label{EqnJensen}
\log g_{\pmstwo}(\obs) & \geq \underbrace{\int_{\HiddenSpace}
  \condens_{\pmsone}(\hidden \mid \obs) \log f_{\pmstwo}(\obs,\hidden)
  d \hidden}_{\QFUN{\pmstwo}{\pmsone}} - \int_{\HiddenSpace}
\condens_{\pmsone}(\hidden \mid \obs) \log \condens_{\pmsone}(\hidden
\mid \obs) d\hidden,
\end{align}
with equality holding when $\pmsone = \pmstwo$.  Thus, we have a
family of lower bounds on the log likelihood, and the EM algorithm
successively maximizes this lower bound ($M$-step), and then
reevaluates the lower bound at the new parameter value ($E$-step).

 
\paragraph{Standard EM updates:}
With this notation, it is easy to specify the EM
iterations.  The update $\pms^t \rightarrow \pms^{t+1}$ consists of
the following two steps.
\begin{itemize}
\item E-step: Evaluate the expectation in equation~\eqref{EqnJensen}
  to compute $\QFUN{\cdot}{\pms^t}$.
\item M-step: Compute the maximizer $\pms^{t+1} = \arg \max
  \limits_{\pms' \in \ParSpace} \QFUN{\pms'}{\pms^t}$.
\end{itemize}
\noindent For future use, it is convenient to introduce the mapping
$\memzero: \ParSpace \rightarrow \ParSpace$ given by
\begin{align}
\label{EqnEMMAP}
\mem{\theta} & \defn \arg \max \limits_{\pms' \in \ParSet}
\QFUN{\pms'}{\pms}.
\end{align}
With this choice, the $M$-step corresponds to the update
$\pms^{t+1} = \mem{\pms^t}$. 


\paragraph{Generalized EM updates:}
In a \emph{generalized EM} algorithm, the requirements of the $M$-step
are relaxed: instead of finding the exact optimum, the algorithm is
required only to find a value $\pms^{t+1} \in \ParSet$ such that
\begin{align}
\label{EqnGeneralizedEM}
\QFUN{\pms^{t+1}}{\pms^t} & \geq \QFUN{\pms^t}{\pms^t}.
\end{align}
Depending on how $\pms^{t+1}$ is chosen, this requirement actually
defines a family of algorithms.


\paragraph{Gradient EM updates:}
A closely related variant of the generalized EM updates is what we
refer to as the \emph{gradient EM updates}, applicable in the case
when the function $\QFUN{\cdot}{\pms^t}$ is differentiable at each
iteration $t$.  Given a step size $\stepsize > 0$, these updates take
the form
\begin{equation}
\pms^{t+1} = \pms^t + \stepsize \nabla \QFUN{\pms^t}{\pms^t},
\end{equation}
where the gradient is taken in the first argument of $\Qfun$.  For
ease of notation, we define the mapping $\mgradzero\colon \pmss
\rightarrow \pmss$ by
\begin{align}
\label{EqnGradEMOperator}
\mgrad{\pms} = \pms + \stepsize \nabla \QFUN{\pms}{\pms}.
\end{align}
An iteration of gradient EM can now be written compactly as
$\theta^{t+1} = \mgrad{\theta^t}$.

There is a natural extension that includes a constraint arising from
the parameter space $\ParSpace$, in which the update is projected back
onto the constraint set
\footnote{To avoid pathologies additionally assume 
that the constraint set is closed.}.  
For simplicity, we focus on 
unconstrained problems in this paper, but all of our
results extend in a straightforward way to constrained examples by
incorporating the additional Euclidean projection.  For appropriate
choices of the step size parameter $\stepsize$, the gradient EM
updates guarantee the ascent condition~\eqref{EqnGeneralizedEM}, so
that it is a particular case of a generalized EM algorithm.

\paragraph{Population versus sample updates:}

Let us now make an important distinction, namely, that between the
population and sample-based versions of the EM updates.  Up to this
point, we have suppressed dependence on the number of observed samples
$\numobs$.  The population form of the (gradient) EM updates are an
``oracle version'', in which we effectively observe an infinite number
of samples, and consequently, the function $Q(\cdot|\theta)$ takes the form
\begin{align}
\label{EqnPopulationQ}
\QFUN{\theta'}{\theta} = \int_{\ObsSpace} \Big(\int_{\HiddenSpace}
\condens_{\pmsone}(\hidden \mid \obs) \log f_{\pmstwo}(\obs,\hidden) d
\hidden \Big) g_{\pmstar}(\obs) d \obs.
\end{align}
From here onwards, we use the notation $\PopEM$ and $\PopGradEM$ for
the EM and gradient EM operators, respectively, both defined at the
population level.

In the classical statistical settings, we observe only $\numobs$ i.i.d. samples
$\{\obs_i\}_{i=1}^\numobs$ of the $\Obs$ component. 
Under the i.i.d. assumption, we define the function
\begin{align}
\label{EqnSampleQ}
\QFUNHAT{\theta'}{\theta} & = \frac{1}{\numobs} \sum_{i=1}^\numobs
\Big ( \int_{\HiddenSpace} \condens_{\pmsone}(\hidden \mid \obs_i)
\log f_{\pmstwo}(\obs_i,\hidden) d \hidden \Big),
\end{align}
so that the expectation over $\Obs$ in equation~\eqref{EqnPopulationQ}
is replaced by the empirical expectation defined by the samples.  The
function $\Qhat$ defines an analog of the population EM
operator~\eqref{EqnEMMAP}, namely
\begin{align}
\label{EqnMemSamp}
\MemSamp(\theta) & = \arg \max_{\theta' \in \ParSpace}
\QFUNHAT{\theta'}{\theta}.
\end{align}
In an analogous fashion, we define the sample-based analog of the
gradient EM operator~\eqref{EqnGradEMOperator}, namely
\begin{align}
\label{EqnGremSamp}
\GremSamp(\pms) & \defn \pms + \stepsize \nabla
\QFUNHAT{\theta}{\theta},
\end{align}
where $\stepsize > 0$ is an appropriately chosen step size parameter.


\subsection{Illustrative examples}
\label{SecIllustrative}

The EM algorithm is popular and a variety of examples can be found in
the literature.  In this section, we review three specific models
analyzed in this paper, and derive the form of the population and
sample-based updates, both for the usual EM algorithm and the gradient
EM algorithm.

\subsubsection{Gaussian mixture models}
\label{SecExaGMM}

An isotropic, balanced two-component Gaussian mixture model can be
specified by a density of the form
\begin{align}
\label{EqnMixtureDensity}
f_\theta(y) & = \frac{1}{2} \phi(y; \thetastar, \sigma^2 I_{\usedim})
+ \frac{1}{2} \phi(y; -\thetastar, \sigma^2 I_{\usedim}),
\end{align}
where $\phi(\cdot \, ; \mu, \Sigma)$ denotes the density of a
$\NORMAL(\mu, \Sigma)$ random vector in $\real^d$.  Here we have
assumed that the components are equally weighted; with the variance
$\sigma^2$ known, the goal is to estimate the unknown mean vector
$\thetastar$.  In this example, the hidden variable $\Hidden \in
\{0,1\}$ is an indicator variable for the underlying mixture
component---that is
\begin{align*}
(\Obs \mid \Hidden = 0) \sim \NORMAL(-\thetastar, \sigma^2 I_\usedim), \quad
  \mbox{and} \quad (\Obs \mid \Hidden = 1) \sim \NORMAL(\thetastar,
  \sigma^2 I_\usedim).
\end{align*}

Suppose that we are given $\numobs$ i.i.d. samples
$\{y_i\}_{i=1}^\numobs$ drawn from the mixture
density~\eqref{EqnMixtureDensity}.  The complete data $\{(y_i,
z_i)\}_{i=1}^\numobs$ corresponds to the original samples along with
the component indicator variables $z_i \in \{0,1\}$.  The sample-based
function $\Qhat$ takes the form
\begin{align}
\label{EqnQfunhatGMM}
\QFUNHAT{\theta'}{\theta} = - \fsavgtwo{1}{2} \big[ w_{\theta}(y_i)
  \norm{y_i - \theta'}^2 + (1 - w_{\theta}(y_i)) \norm{y_i + \theta'}^2 \big], 
  \end{align}
where $w_{\theta}(y) \defn \MYEXP{- \frac{\norm{\theta - y}^2}{2
    \sigma^2}} \Big[ \MYEXP{ - \frac{\norm{\theta - y}^2}{2 \sigma^2} }
    + \MYEXP{ - \frac{\norm{\theta + y}^2}{2 \sigma^2}} \Big]^{-1}$.

\paragraph{EM updates:}
This example is especially simple in that the EM operator $\MemSamp:
\real^\usedim \rightarrow \real$ has a closed form solution, given by
\begin{subequations}
\begin{align}
\label{EqnMemSampMOG}
\MemSamp(\theta) & \defn \arg \max_{\theta' \in \real^\usedim}
\QFUNHAT{\theta'}{\theta} \; = \;  \frac{2}{\numobs}
  \sum_{i=1}^\numobs w_{\theta}(y_i) y_i - \frac{1}{\numobs} \sum_{i=1}^\numobs y_i.
\end{align}
The population EM operator $\Mem: \real^\usedim \rightarrow
\real^\usedim$ is defined analogously
\begin{align}
\label{EqnMemMOG}
\mem{\theta} & = 2 \Exs \big[ w_{\theta}(Y) Y \big],
\end{align}
\end{subequations}
where the empirical expectation has been replaced by expectation under
the mixture distribution~\eqref{EqnMixtureDensity}.

\paragraph{Gradient EM updates:} On the other hand, the sample-based
and population gradient EM operators with step size $\stepsize > 0$
are given by
\begin{align}
\GremSamp(\theta) = \theta + \stepsize \: \Big \{ \frac{1}{\numobs}
\sum_{i=1}^\numobs (2 w_{\theta}(y_i) - 1) y_i - \theta  \Big \}, ~~
\mbox{and} ~~ \Grem(\theta) = \theta + \stepsize \big[\: 2\Exs \big[
 w_{\theta}(Y) Y\big] - \theta \big].
\end{align}

\noindent We return to analyze the EM updates for the Gaussian mixture
model in Section~\ref{SecGMMAnalysis}.


\subsubsection{Mixture of regressions}
\label{SecExaMOR}

We now consider the mixture of regressions model, as has been analyzed
in some recent work~\cite{convex_mor,hardem_mor,spectral_mor}.  In the
standard linear regression model, we observe i.i.d. samples of the
pair $(Y, X) \in \real \times \real^{\usedim}$ linked via the equation
\begin{align}
\label{EqnLinearRegression}
y_i & = \inprod{x_i}{\thetastar} + \regnoise_i,
\end{align}
where $\regnoise_i \sim \NORMAL(0, \sigma^2)$ is the observation noise
assumed to be independent of $x_i$,
$x_i \sim \NORMAL(0,I)$ are the design vectors and $\thetastar \in
\real^\usedim$ is the unknown regression vector to be estimated. In
the mixture of regressions problem, there are two underlying choices
of regression vector---say $\thetastar$ and $-\thetastar$---and we
observe a pair $(y_i, x_i)$ drawn from the
model~\eqref{EqnLinearRegression} with probability $\frac{1}{2}$, and
otherwise generated according to the alternative regression model $y_i
= \inprod{x_i}{-\thetastar} + \regnoise_i$.  Here the hidden variables
$\{z_i\}_{i=1}^\numobs$ correspond to labels of the underlying
regression model: say $z_i = 1$ when the data is generated according
to the model~\eqref{EqnLinearRegression}, and $z_i = 0$ otherwise. In
this symmetric form, the mixture of regressions model is closely
related to models for phase retrieval, albeit over $\real^\usedim$, as
considered in a line of recent work
(e.g.,~\cite{candes_pr,altmin_pr,balan_pr}).

\paragraph{EM updates:}
Define the weight function
\begin{subequations}
\begin{align}
\label{EqnMORweight}
w_{\theta}(x,y) = \frac{ \exp \big( \frac{- (y -
    \inprod{x}{\theta})^2 }{2\sigma^2} \big)} { \exp \big( \frac{-(y
    - \inprod{x}{\theta})^2}{2\sigma^2} \big) + \exp \big( \frac{-(y
    + \inprod{x}{\theta})^2 }{2\sigma^2} \big)}.
\end{align}
In terms of this notation, the sample EM update is based on maximizing
the function
\begin{align}
\label{EqnMORQfun}
\QFUN{\theta'}{\theta} & = - \frac{1}{2 \numobs} \sum_{i=1}^\numobs
\Big( w_{\theta} (x_i,y_i) (y_i- \inprod{x_i}{\theta'})^2 + (1 -
w_{\theta} (x_i,y_i)) (y_i + \inprod{x_i}{\theta'})^2 \Big).
\end{align}
\end{subequations}
Again, there is a closed form solution to this maximization problem:
more precisely, the sample EM operator $\MemSamp: \real^\usedim
\rightarrow \real^\usedim$ takes the form
\begin{subequations}
\begin{align}
\label{EqnMemSampMOR}
\MemSamp(\theta) & = \Big( \sum_{i=1}^n x_i x_i^T \Big)^{-1} \Big(
\sum_{i=1}^\numobs (2 w_{\rvec}(x_i,y_i) - 1)y_i x_i\Big).
\end{align}
Similarly, by an easy calculation, we find that the
population EM operator \mbox{$\Mem: \real^\usedim \rightarrow
  \real^\usedim$} has the form
\begin{align}
\label{EqnMemMOR}
\Mem(\theta) & = 2 \Exs \big[ w_{\theta}(X,Y) Y X \big],
\end{align}
\end{subequations}
where the expectation is taken over the joint distribution of the pair
$(Y, X) \in \mathbb{R} \times \mathbb{R}^d$.

\paragraph{Gradient EM updates:}
On the other hand, the gradient EM operators are given by
\begin{subequations}
\begin{align}
\label{EqnGradientMOR}
\GremSamp(\theta) & = \theta + \stepsize \Big \{ \frac{1}{\numobs}
\sum_{i=1}^\numobs \Big[ (2w_{\theta} (x_i,y_i) - 1) y_i x_i - x_i x_i^T \theta \Big] \Big \}, \quad \mbox{and} \quad \\
\Grem(\theta) & = \theta + \stepsize \, 2\Exs \Big[ w_{\theta} (X, Y)
  \, YX - \theta \Big],
\end{align}
\end{subequations}
where $\stepsize > 0$ is a step size parameter. \\

\noindent We return to analyze the EM updates for the mixture of
regressions model in Section~\ref{SecMORAnalysis}.


\subsubsection{Linear regression with missing covariates}
\label{SecExaMissing}

Our first two examples involved mixture models in which the class
membership variable was hidden.  Another canonical use of the EM
algorithm is in cases with corrupted or missing data.  In this
section, we consider a particular instantiation of such a problem,
namely that of linear regression with the covariates missing
completely at random.

As introduced in Section~\ref{SecExaMOR}, in standard linear
regression, we observe response-covariate pairs \mbox{$(y_i, x_i) \in
  \real \times \real^\usedim$} generated according to the linear
model~\eqref{EqnLinearRegression}. In the missing data extension of
this problem, instead of observing the covariate vector $x_i \in
\real^\usedim$ directly, we observe the corrupted version $\xtil_i \in
 \real^\usedim$ with components
\begin{align}
\xtil_{ij} & = \begin{cases} x_{ij} & \mbox{with probability
    $1-\missing$} \\ \ast & \mbox{with probability $\missing$,}
\end{cases}
\end{align}
where $\missing \in [0,1)$ is the probability of missingness.

In this example, the E-step involves imputing the mean and covariance
of the jointly Gaussian distribution of covariate-response pairs.  For
a given sample $(x,y)$, let $\xobs$ denote the observed portion of
$x$, and let $\bobs$ denote the corresponding sub-vector of $\theta$.
Define the missing portions $\xmis$ and $\bmis$ in an analogous
fashion.  With this notation, the EM algorithm imputes the conditional
mean and conditional covariance using the current parameter estimate
$\theta$.  Using properties of joint Gaussians, the conditional mean
of $X$ given $(\xobs, y)$ is found to be
\begin{subequations}
\begin{align}
\label{EqnMuMissing}
\mu_{\mvec}(\xobs,y) \defn \begin{bmatrix} \Exs (\xmis | \xobs, y,
  \mvec) \\ \xobs \end{bmatrix} = \begin{bmatrix} \mism_{\mvec} \zobs
  \\ \xobs \end{bmatrix},
\end{align}
where 
\begin{align}
\label{EqnMuMissingAux}
\mism_{\mvec} & = \frac{1}{\norm{\bmis}^2 + \sigma^2}
\begin{bmatrix}
- \bmis \: \bobs^T & \bmis
\end{bmatrix} 
\quad \mbox{and} \quad \zobs \defn \begin{bmatrix} \xobs \\ y
\end{bmatrix}
\in \real^{\sobs + 1}.
\end{align}
Similarly, the conditional second moment matrix takes the form
\begin{align}
\label{EqnCovMissing} 
\CovMat_{\mvec}(\xobs,y) & \defn \Exs \Big[ XX^T \mid \xobs, y, \mvec
  \Big] =
\begin{bmatrix} I & \mism_{\mvec} \zobs \xobs^T \\ \xobs
    \zobs^T \mism_{\mvec}^T & \xobs \xobs^T \end{bmatrix}.
\end{align}
\end{subequations}
In writing all these expressions, we have
assumed that the coordinates are permuted so that the missing values
are in the first block.

We now have the necessary notation in place to describe the EM and
gradient EM updates.  For a given parameter $\theta$, the EM update is
based on maximizing
\begin{align}
\label{EqnMissingQsamp}
\QFUNHAT{\theta'}{\theta} & \defn - \fsavgtwo{1}{2}
\inprod{\theta'}{\CovMat_{\mvec}(\xobsi,y_i) \theta'} + \fsavg y_i
\inprod{\mu_{\mvec}(\xobsi,y_i)}{\theta'}.
\end{align}
Again, this optimization problem has an explicit solution, so that the
sample-based EM operator is given by
\begin{subequations}
\begin{align}
\label{EqnMemSampMissing}
M_\numobs(\theta) \defn \Big[ \sum_{i=1}^\numobs
  \CovMat_{\theta}(\xobsi,y_i) \Big]^{-1} \Big[ \sum_{i=1}^n y_i
  \mu_{\mvec}(\xobsi,y_i)\Big],
\end{align}
accompanied by its population counterpart
\begin{align}
\label{EqnMemMissing}
M(\theta) & \defn \big \{ \Exs \big[
  \CovMat_{\theta}(X_{\mathrm{obs}},Y) \big] \big \}^{-1} \; \Exs
\big[ Y \mu_{\mvec}(X_{\mathrm{obs}},Y) \big].
\end{align}
\end{subequations}
On the other hand, the gradient EM algorithm with step size
$\stepsize$ takes the form
\begin{subequations}
\begin{align}
\label{EqnMissingGremSamp}
G_\numobs(\theta) & = \theta + \stepsize \Big\{ \frac{1}{\numobs}
\sum_{i=1}^\numobs \big[ y_i \mu_{\theta}(\xobsi, y_i) -
\CovMat_{\mvec}(\xobsi, y_i) \theta \big] \Big\},
\end{align}
along with the population counterpart
\begin{align}
\label{EqnMissingGrem}
G(\theta) & = \theta + \stepsize \Exs \Big[Y \mu_{\theta}(\Xobs, Y) -
  \CovMat_{\mvec}(\Xobs, Y) \theta \Big],
\end{align}
\end{subequations}
We return to analyze the gradient EM updates for this model in
Section~\ref{SecMissingAnalysis}.

\section{General convergence results}
\label{SecMain}

We now turn to analysis of the EM algorithm and gradient EM
algorithms. In both cases, we let $\fp$ denote a maximizer of the
population likelihood. In this section, we give general sufficient
conditions under which the population algorithms converge to $\fp$ and
under which the sample-based algorithms converge to an
$\varepsilon$-ball around $\fp$.  Our analysis of each algorithm is
organized as follows:

   Our first result in Sections~\ref{SecEMAnalysis} and
  \ref{SecGREMAnalysis} concern the population EM and gradient EM
  operators respectively.  Theorems~\ref{ThmEM} and~\ref{ThmGradEM}
  give conditions under which the population operators are contractive
  on a ball containing the fixed point $\fp$, say $\Balltwor = \{
  \pms \in \ParSpace \, \mid \, \|\pms - \fp\|_2 \leq r\}$ for some
  \mbox{radius $r$.}  This population-level analysis is developed by
  viewing the population operators as perturbed versions of
  \emph{oracle} operators which are known to be contractive around
  $\fp$. Our conditions which relate the population EM and gradient EM
  operators to the oracle operators are then verified in concrete
  examples in Section~\ref{SecConsequences}. The analysis here is
  entirely deterministic.

   Our second result in Sections~\ref{SecEMAnalysis} and \ref{SecGREMAnalysis} 
concern the sample-based EM and gradient EM operators.
These sample-based operators approximate the population-based update using a
subset of samples at each step. 
  Theorem~\ref{ThmSampleEM}
  and Theorem~\ref{ThmSampleGradEM} for sample-based EM and gradient
  EM, respectively give conditions under which the sample-based operator
 is guaranteed to converge to an
  $\varepsilon$-ball around the fixed point $\fp$.  These results involve 
  probabilistic bounds on the deviations between the population-based
  and sample-based operators.  In addition, for gradient EM, we also
  analyze a stochastic update that uses a single sample per update in
  the flavor of stochastic approximation algorithms (see
  Theorem~\ref{ThmStochasticEM} in Section~\ref{SecStochEMAnalysis}).
%


\subsection{Analysis of EM algorithm}
\label{SecEMAnalysis}

Let us begin with analysis of the standard EM updates, starting with
the population version before turning to a sample-based version.

\subsubsection{Guarantees for population-level EM}
Recall that we always assume that the vector $\fp$ maximizes the
population likelihood.  It is a classical fact~\cite{em_textbook} that
it must then satisfy the condition
\begin{align}
\label{EqnSelfConsistency}
\fp & = \arg \max_{\pms \in \pmss} \QFUN{\pms}{\fp},
\end{align}
a property known as \emph{self-consistency}.  For this reason, the
function $\specq(\cdot) \defn \QFUN{\cdot}{\fp}$ plays an important
role in our analysis. 

We assume throughout this section that the function $q$ is
\emph{$\strongparam$-strongly concave}, meaning that
\begin{align}
\label{EqnStronglyConcave}
\HACKQ(\pms_1) - \HACKQ(\pms_2) - \inprod{\nabla
  \HACKQ(\pms_2)}{\pms_1-\pms_2} & \leq -\frac{\strongparam}{2} \norm{
  \pms_1- \pms_2}^2,
\end{align}
for all pairs $(\pms_1, \pms_2)$ in a neighborhood of $\thetastar$.
As we will illustrate, this condition holds in most concrete
instantiations of EM, including the three model classes introduced in
the previous section.

For any fixed $\theta$, in order to relate the population EM updates
to the fixed point $\fp$, we require control on the two gradient
mappings $\nabla \HACKQ(\cdot) = \nabla \qfunone{\fp}$ and $\nabla
\qfunone{\theta}$.  These mappings are central in characterizing the
fixed point $\fp$ and the update $\PopEMupdate(\theta)$ respectively.
Indeed, by virtue of the self-consistency
property~\eqref{EqnSelfConsistency} and the convexity of $\ParSet$,
the fixed point satisfies the first-order optimality condition
\begin{align}
\label{EqnOracleBasic}
\inprod{\nabla \qfun{\fp}{\fp}}{\pms' - \fp} & \leq 0 \qquad \mbox{for all
  $\pms' \in \ParSpace$.}
\end{align}
Similarly, for any $\pms \in \ParSpace$, since $\PopEMupdate(\theta)$
maximizes the function $\pms' \mapsto Q(\pms'|\pms)$, we have
\begin{align}
\label{EqnUpdateBasic}
\inprod{\nabla \qfun{\PopEMupdate(\pms)}{\pms}}{\pms' -
  \PopEMupdate(\theta)} & \leq 0 \qquad \mbox{for all $\pms'
  \in \ParSpace$.}
\end{align}
Equations~\eqref{EqnOracleBasic} and \eqref{EqnUpdateBasic} are sets
of inequalities that \emph{characterize} the points $\mem{\pms}$ and
$\fp$.  Thus, at an intuitive level, in order to establish that
$\mem{\pms}$ and $\fp$ are close, it suffices to verify that these two
characterizations are close in a suitable sense.  We also note that
inequalities similar to the condition~\eqref{EqnUpdateBasic} are often
used as a starting point in the classical analysis of M-estimators
(e.g., see van de Geer~\cite{vandeGeer}). In the analysis of EM, we
obtain additional leverage from the self-consistency
condition~\eqref{EqnSelfConsistency} that characterizes $\thetastar$.

With this intuition in mind, we introduce the following regularity
condition in order to relate conditions~\eqref{EqnUpdateBasic}
and~\eqref{EqnSelfConsistency}: The condition involves a Euclidean
ball of radius $r$ around the fixed point $\fp$, given by
\begin{align}
\Balltwor & \defn \big \{ \theta \in \ParSpace \, \mid \|\theta -
\fp\|_2 \leq r \big \}.
\end{align}
\bdes[\gstwotext~(\gstwonospace)] The functions \mbox{$\{
  \qfunone{\pms}, \pms \in \Omega \}$} satisfy
condition~\gstwo($\gstwoparam$) over $\Balltwor$ if
\begin{align}
\label{EqnGSTWO}
\norm{\nabla \qfun{\mem{\pms}}{\fp} - \nabla \QFUN{\mem{\pms}}{\pms} }
& \leq \gstwoparam \norm{\pms - \fp} \qquad \mbox{for all $\pms \in
  \Balltwor$.}
\end{align}
\edes

\vspace*{.02in}

\noindent 
To provide some high-level intuition, observe the
condition~\eqref{EqnGSTWO} is always satisfied at the fixed point
$\theta^*$, in particular with parameter $\gstwoparam = 0$.
Intuitively then, by allowing for a strictly positive parameter
$\gstwoparam$, one might expect that this condition would hold in a
local neighborhood $\Balltwor$ of the fixed point $\thetastar$, as
long as the functions $Q(\cdot|\theta)$ and the map $M$ are
sufficiently regular.

As a concrete example, recall the Gaussian mixture model first
introduced in Section~\ref{SecExaGMM}.  For this model, the
condition~\eqref{EqnGSTWO} is equivalent to
\begin{align*}
\Exs \Big[ 2 \big( w_\theta(Y) - w_{\thetastar}(Y) \big) Y \Big] &
\leq \gstwoparam \, \|\theta - \thetastar\|_2,
\end{align*}
where $w_\theta$ was previously defined following
equation~\eqref{EqnQfunhatGMM}.  Given that the function $\theta
\mapsto w_\theta(y)$ is smooth in $\theta$, provided 
that $\gamma$ is not too small, it is reasonable to expect
that this condition will hold in a neighborhood of $\thetastar$, and
we confirm this intuition in Corollary~\ref{CorMOGPop} to follow. \\

\vspace*{.05in}

\noindent Under the conditions we have introduced, the following
result guarantees that the population EM operator is locally
contractive:
\btheos
\label{ThmEM}
For some radius $r > 0$ and pair $(\gstwoparam,\strongparam)$ such that
 $0 \leq \gstwoparam < \strongparam$,
suppose that the function \mbox{$Q(\cdot|\theta^*)$} is
$\strongparam$-strongly concave~\eqref{EqnStronglyConcave}, and that the
FOS$(\gstwoparam)$ condition~\eqref{EqnGSTWO} holds on the ball
$\Balltwor$.  Then the population EM operator $\PopEM$ is contractive
over $\Balltwor$, in particular with
\begin{align*}
\norm{ \mem{\pms} - \fp} \leq \frac{\gstwoparam}{\strongparam}
\norm{\pms - \fp} \qquad \mbox{for all $\pms \in \Balltwor$.}
\end{align*}
 \etheos
\noindent As an immediate consequence, under the conditions of the
theorem, for any initial point \mbox{$\pms^0 \in \Balltwor$,} the
population EM sequence $\{\pms^t\}_{t=0}^\infty$ exhibits linear
convergence---viz.
\begin{align}
\norm{\pms^t - \fp} & \leq \Big(\frac{\gstwoparam}{\strongparam}
\Big)^t \; \norm{\pms^0 - \fp} \qquad \mbox{for all $t = 1, 2,
  \ldots$.}
\end{align}
\begin{proof}
Since both $\mem{\pms}$ and $\fp$ are in $\pmss$, we may apply
condition~\eqref{EqnOracleBasic} with $\pms' = \mem{\pms}$ and
condition~\eqref{EqnUpdateBasic} with $\pms' = \fp$.  Doing so, adding
the resulting inequalities and then performing some algebra yields the
condition
\begin{align}
\label{EqnBasicAdd}
\inprod{\nabla  \qfun{\mem{\pms}}{\fp} - \nabla \qfun{\fp}{\fp}}{\fp -
  \mem{\pms} } & \leq \inprod{\nabla  \qfun{\mem{\pms}}{\fp}- \nabla
  \QFUN{\mem{\pms}}{ \pms} }{\fp - \mem{\pms} }.
\end{align}
Now the $\strongparam$-strong concavity
condition~\eqref{EqnStronglyConcave} implies that the left-hand side
is lower bounded as%
\begin{subequations}
\begin{align}
\label{EqnHanaSleepOne}
\inprod{\nabla  \qfun{\mem{\pms}}{\fp} - \nabla \qfun{\fp}{\fp}}{\fp -
  \mem{\pms} } & \geq \strongparam \|\fp - \mem{\pms}\|_2^2.
\end{align}
On the other hand, the FOS($\gstwoparam$) condition together
with the Cauchy-Schwarz inequality implies that
the right-hand side is upper bounded as
\begin{align}
\label{EqnHanaSleepTwo}
\inprod{\nabla  \qfun{\mem{\pms}}{\fp} - \nabla \QFUN{\mem{\pms}}{
    \pms}}{\fp -\mem{\pms}} \leq \gstwoparam \norm{\fp - \mem{\pms}}
\norm{\pms - \fp},
\end{align}
\end{subequations}
Combining inequalities~\eqref{EqnHanaSleepOne}
and~\eqref{EqnHanaSleepTwo} with the original
bound~\eqref{EqnBasicAdd} yields 
\begin{align*}
\strongparam \|\fp - \mem{\pms}\|_2^2 & \leq \gstwoparam \norm{\fp -
  \mem{\pms}} \norm{\pms - \fp},
\end{align*}
and canceling terms completes the proof.
\end{proof}


\subsubsection{Guarantees for sample-based EM}

We now turn to theoretical results on sample-based versions of the EM
algorithm.  More specifically, we consider two forms of the EM
algorithm, the first being the standard form in which the operator
$\MemSamp: \ParSpace \rightarrow \ParSpace$, as previously
defined~\eqref{EqnMemSamp}, is applied repeatedly, thereby generating
the sequence $\pmsit{t+1} = \MemSamp(\pmsit{t})$.  We also analyze a
sample-splitting
\footnote{From a practical point, a potential advantage of sample
  splitting is that each iteration may be cheaper, since it is based
  on a smaller sample size.  In contrast, a disadvantage is that it
  can be difficult to correctly specify the number of iterations in
  advance.}
version of the EM algorithm, in which given a total of $\numobs$
samples and $\Tfinal$ iterations, we divide the full data set into
$\Tfinal$ subsets of size $\lfloor \numobs/\Tfinal \rfloor$, and then
perform the updates $\pmsit{t+1} = \Mem_{\numobs/\Tfinal}(\pmsit{t})$,
using a fresh subset of samples at each iteration.

For a given sample size $\numobs$ and tolerance parameter $\delta \in
(0,1)$, we let $\rateem(\numobs, \delta)$ be the smallest scalar such
that, for any fixed $\pms \in \Balltwor$, we have
\begin{align}
\label{EqnRateEM}
\norm{\memsamp{\pms}- \mem{\pms}} & \leq \rateem(\numobs,\delta)
\end{align}
with probability at least $1-\delta$.  This tolerance
parameter~\eqref{EqnRateEM} enters our analysis of the
sample-splitting form of EM.  On the other hand, in order to analyze
the standard sample-based form of EM, we require a stronger condition,
namely one in which the bound~\eqref{EqnRateEM} holds uniformly over
the ball $\Balltwor$.  Accordingly, we let $\rateemunif(\numobs,
\delta)$ be the smallest scalar for which
\begin{align}
\label{EqnRateEMUniform}
\sup_{\pms \in \Balltwor} \norm{\memsamp{\pms}- \mem{\pms}} & \leq
\rateemunif(\numobs,\delta)
\end{align}
with probability at least $1 - \delta$.  With these definitions, we
have the following guarantees:
\btheos
\label{ThmSampleEM}
Suppose that the population EM operator $\PopEM: \ParSpace
\rightarrow \ParSpace$ is contractive with parameter $\kappa \in
(0,1)$ on the ball $\Balltwor$, and the initial vector $\thetasamp^0$
belongs to $\Balltwor$.
\begin{enumerate}
\item[(a)] If the sample size $\numobs$ is large enough to ensure that
\begin{subequations}
\begin{align}
\label{EqnSamplingBoundUniform}
\rateemunif(\numobs,\delta) & \leq (1 - \kappa) r,
\end{align}
then the EM iterates $\{\theta^t\}_{t=0}^\infty$ satisfy the bound
\begin{align}
\label{EqnSampleEMRateUniform}
\norm{\thetasamp^t - \fp} & \leq \kappa^t \norm{ \pms^0 - \fp} +
\frac{1}{1 - \kappa} \; \rateemunif(\numobs,\delta)
\end{align}
with probability at least $1 - \delta$.
\end{subequations}
\item[(b)] For a given iteration number $\Tfinal$, suppose the sample
  size $\numobs$ is large enough to ensure that%
\begin{subequations}
\begin{align}
\label{EqnSamplingBound}
\HACKRATEEM & \leq (1 - \kappa) r.
\end{align}
Then the sample-splitting EM iterates $\{\pmsit{t}\}_{t=0}^\Tfinal$
based on $\frac{\numobs}{\Tfinal}$ samples per round satisfy the bound
\begin{align}
\label{EqnSampleEMRate}
\norm{\thetasamp^t - \fp} & \leq \kappa^t \norm{ \pms^0 - \fp} +
\frac{1}{1 - \kappa} \; \rateem \Big(\frac{\numobs}{\Tfinal},
\frac{\delta}{\Tfinal} \Big).
\end{align}
\end{subequations}
\end{enumerate}
\etheos
\noindent Figure~\ref{FigStatPrec} provides an illustration of the
behavior predicted by Theorem~\ref{ThmSampleEM}: both algorithms are
expected to show geometric convergence to the target parameter
$\thetastar$, up to some tolerance.  
For the bound~\eqref{EqnSampleEMRateUniform}
note that the first term is decreasing in $t$, whereas the second
term is independent of $t$.
Thus, for a fixed sample size $\numobs$,
the bounds in Theorem~\ref{ThmSampleEM} suggests a reasonable choice
of the number of iterations.  In particular, focusing on the standard
EM algorithm, consider any positive integer\footnote{ As will
  be clarified in the sequel, such a choice of $\Tcrit$ exists in
  various concrete models considered here.}  $\Tcrit$ such that
\begin{align}
\label{EqnTcritChoiceUniform}
\Tcrit \geq \log_{1/\contractparam} \frac{ (1 - \contractparam) \, \norm{\theta^0 - \fp}
}{\rateemunif ( \numobs, \delta)}.
\end{align}
This choice ensures that the first term in the
bound~\eqref{EqnSampleEMRateUniform} is dominated by the second term,
and hence that
\begin{align}
\norm{\pms^\Tcrit - \fp} & \leq \frac{2}{1 - \kappa} \;
\rateemunif(\numobs, \delta),
\end{align}
with probability at least $1 - \delta$. For the sample-splitting update
in~\eqref{EqnSampleEMRate} the first term 
is decreasing in $t$, whereas the second
term is increasing in $t$. In this case, a similar conclusion holds when $\Tcrit$
is chosen to be the smallest positive integer such that
\begin{align}
\label{EqnTcritChoice}
\Tcrit \geq \log_{1/\contractparam} \frac{ (1 - \contractparam) \, \norm{\theta^0 - \fp}
}{\rateem \big( \frac{\numobs}{\Tcrit}, \frac{\delta}{\Tcrit} \big)}.
\end{align}

In order to obtain readily interpretable bounds for
specific models, it only remains to establish the
$\kappa$-contractivity of the population operator, and to compute
either the function $\rateem$ or the function
$\rateemunif$. \\ \iftoggle{figs}{
\begin{figure}[h]
\begin{center}
\widgraph{0.6\textwidth}{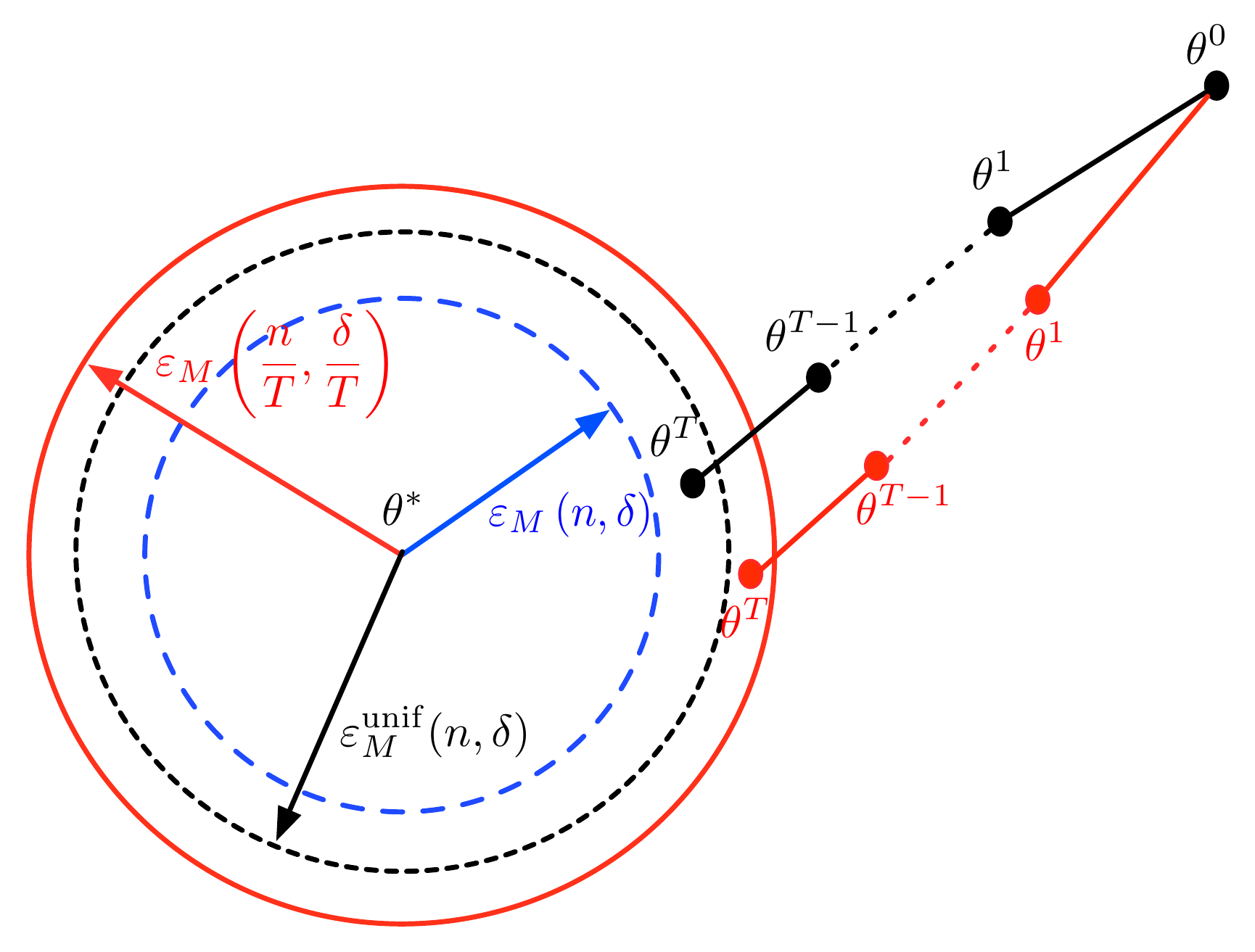}
\end{center}
\caption{An illustration of Theorem~\ref{ThmSampleEM}. The first part
  of the theorem describes the geometric convergence of iterates of
  the EM algorithm to the ball of radius
  $\order(\rateemunif(\numobs,\delta))$ (in black). The second part
  describes the geometric convergence of the sample-splitting EM
  algorithm to the ball of radius $\order(\rateem(\numobs/\Tcrit,
  \delta/\Tcrit))$ (in red). In typical examples the ball
  to which sample-splitting EM converges is only a
  logarithmic factor larger than the ball
  $\order(\rateem(\numobs,\delta))$ (in blue).}
\label{FigStatPrec}
\end{figure}}

\noindent Let us now turn to the proof of the theorem.
\begin{proof}
We give a detailed proof of the claim~\eqref{EqnSampleEMRate}, from
which it will be clear that the claim~\eqref{EqnSampleEMRateUniform}
follows by a nearly identical argument.  For any iteration $s \in
\{1, 2, \ldots, \Tfinal \}$, we have
\begin{align}
\label{EqnPizza}
\norm{\memsubsamp{\pms^s}- \mem{\pms^s}} & \leq \rateem \Big(
\frac{\numobs}{\Tfinal}, \frac{\delta}{\Tfinal} \Big)
\end{align}
with probability at least $1 - \frac{\delta}{\Tfinal}$.  Consequently,
by a union bound over all $\Tfinal$ indices, the
bound~\eqref{EqnPizza} holds uniformly with probability at least $1 -
\delta$.  We perform the remainder of our analysis under this event.

It suffices to show that
\begin{align}
\label{EqnInductionClaim}
\|\pms^{s+1} - \fp\|_2 & \leq \kappa \|\pms^s - \fp\|_2 + \HACKRATEEM
\qquad \mbox{for each iteration $s \in \{1, 2, \ldots, \Tfinal-1\}$.}
\end{align}
Indeed, when this bound holds, we may iterate it to show that
\begin{align*}
\|\pms^t - \fp\|_2 & \leq \kappa \|\pms^{t-1} - \fp\|_2 + \HACKRATEEM
\\ & \leq \kappa \Big \{ \kappa \|\pms^{t-2} - \fp\|_2 + \HACKRATEEM
\Big \} + \HACKRATEEM \\
& \leq \kappa^t \|\pms^0 - \fp\|_2 + \Big \{ \sum_{s = 0}^{t-1}
\kappa^s \Big \} \HACKRATEEM \\
& \leq \kappa^t \|\pms^0 - \fp\|_2 + \frac{1}{1 - \kappa} \,
\HACKRATEEM,
\end{align*}
where the final step follows by summing the geometric series.

It remains to prove the claim~\eqref{EqnInductionClaim}, and we do so
via induction on the iteration number.  Beginning with $s = 1$, we
have
\begin{align*}
\|\pms^{1} - \fp\|_2 \; = \; \| \memsubsamp{\pms^0} - \fp \|_2 &
\stackrel{\mathrm{(i)}}{\leq} \| \mem{\pms^0} - \fp\|_2 + \|
\memsubsamp{\pms^0} - \mem{\pms^0}\|_2 \\
& \stackrel{\mathrm{(ii)}}{\leq} \kappa \|\pms^0 - \fp \|_2 + \HACKRATEEM,
\end{align*}
where step (i) follows by triangle inequality, whereas step (ii)
follows from the bound~\eqref{EqnPizza}, and the contractivity of the
population operator applied to $\pms^0 \in \Balltwor$.  By our
initialization condition and the bound~\eqref{EqnSamplingBound}, note
that we are guaranteed that $\|\pms^1 - \fp\|_2 \leq r$.  

In the induction from $s \mapsto s+1$, suppose that $\|\pms^s -
\fp\|_2 \leq r$, and the bound~\eqref{EqnInductionClaim} holds at
iteration $s$.  The same argument then implies that the
bound~\eqref{EqnInductionClaim} also holds for iteration $s+1$, and
that $\|\pms^{s+1} - \fp\|_2 \leq r$, thus completing the proof.
\end{proof}

\subsection{Analysis of gradient EM algorithm}
\label{SecGREMAnalysis}

We now turn to analysis of the gradient EM algorithm.  As before, we
separate our analysis into two parts, the first
(Theorem~\ref{ThmGradEM}) addressing the behavior of the
population-level operator, and the second
(Theorems~\ref{ThmSampleGradEM} and~\ref{ThmStochasticEM}) providing
guarantees for sample-based updates.

\subsubsection{Guarantees for population-level gradient EM}

Recall that the gradient EM algorithm generates a sequence of iterates
$\{\theta^t\}_{t=0}^\infty$ via the recursion $\pms^{t+1} =
\mgrad{\pms^t}$, where
\begin{align}
\label{EqnGradEMTwo}
\mgrad{\pms} & \defn \pms + \stepsize \nabla \QFUN{\pms}{\pms}.
\end{align}
Here $\stepsize > 0$ is a step size parameter to be chosen.  For
analyzing gradient EM, we also require an additional condition on
the function $\specq(\pms) = \QFUN{\pms}{\fp}$, previously defined in
Section~\ref{SecEMAnalysis}.  In addition to the $\strongparam$-strong
concavity assumption~\eqref{EqnStronglyConcave}, we also assume that
$\specq$ is $\smoothparam$-smooth, 
meaning
that
\begin{align}
\label{EqnSmoothness}
\HACKQ(\pms_1) - \HACKQ(\pms_2) - \inprod{\nabla
  \HACKQ(\pms_2)}{\pms_1-\pms_2} & \geq -\frac{\smoothparam}{2} \norm{
  \pms_1- \pms_2}^2,
\end{align}
for all pairs $(\pms_1, \pms_2).$

In order to gain intuition into the gradient EM algorithm, it is
instructive to compare its iterates with those of standard gradient
ascent on the function $\specq$.  Gradient ascent on $\specq$
performs the updates $\pmstil^{t+1} = \mgazero(\pmstil^t)$, 
where
\begin{align}
\label{EqnStandardGrad}
\mga{\pms} & \defn \pms + \stepsize \nabla q(\theta).
\end{align}
Under the stated strong concavity and smoothness assumptions, it is a
standard result from optimization
theory~\cite{bubeck_co,Bertsekas_nonlin,Nesterov} that the gradient
operator $\OrGrad: \ParSpace \rightarrow \ParSpace$ with step size
choice $\stepsize = \frac{2}{\smoothparam + \strongparam}$ is
contractive, in particular with
\begin{align}
\label{EqnOrGradContract}
\norm{\mga{\pms} - \fp} & \leq \Big( \frac{\smoothparam -
  \strongparam}{\smoothparam + \strongparam} \Big) \norm{\pms -
  \fp} \qquad \mbox{for all $\pms \in \Balltwor$.}
\end{align}
Intuitively, then, if the function $\qfunone{\pms}$ is ``close
enough'' to the function $q(\cdot) = \qfunone{\fp}$, then the gradient
EM operator might be expected to satisfy a similar contractivity
condition.  The closeness requirement is formalized in the following
condition:
\bdes[\gsonetext~(\gsonenospace)] The functions \mbox{$\{
  \qfunone{\pms}, \pms \in \Omega \}$} satisfy
condition~\gsone($\gsoneparam$) over $\Balltwor$ if
\begin{align}
\label{EqnGSONE}
\norm{ \nabla \qfun{\pms}{\fp} - \nabla \QFUN{\pms}{\pms}} & \leq
\gsoneparam \norm{\pms - \fp} \qquad \mbox{for all $\pms \in \Balltwor$.}
\end{align}
\edes
\noindent See Figure~\ref{FigGS} for an illustration of this
condition.  We give concrete examples of this condition and its
verification in Section~\ref{SecConsequences}.  As with the FOS
condition observe that the GS condition is always satisfied at the fixed
point $\theta^*$, i.e. for $r =0$ with $\gamma=0$. Allowing for
strictly positive $\gamma$, if the functions
$Q(\cdot|\theta)$ are sufficiently regular we expect the condition to
hold in a region around $\theta^*$.  Observe that this condition
involves the gradient of the functions $Q(\cdot|\theta)$ and
$Q(\cdot|\fp)$ at $\theta$, as opposed to $\mem{\theta}$ in the case
of the FOS condition.  For this reason, it can be easier to verify for
specific models. \\
\iftoggle{figs}{
\begin{figure}
\begin{center}
\widgraph{0.7\textwidth}{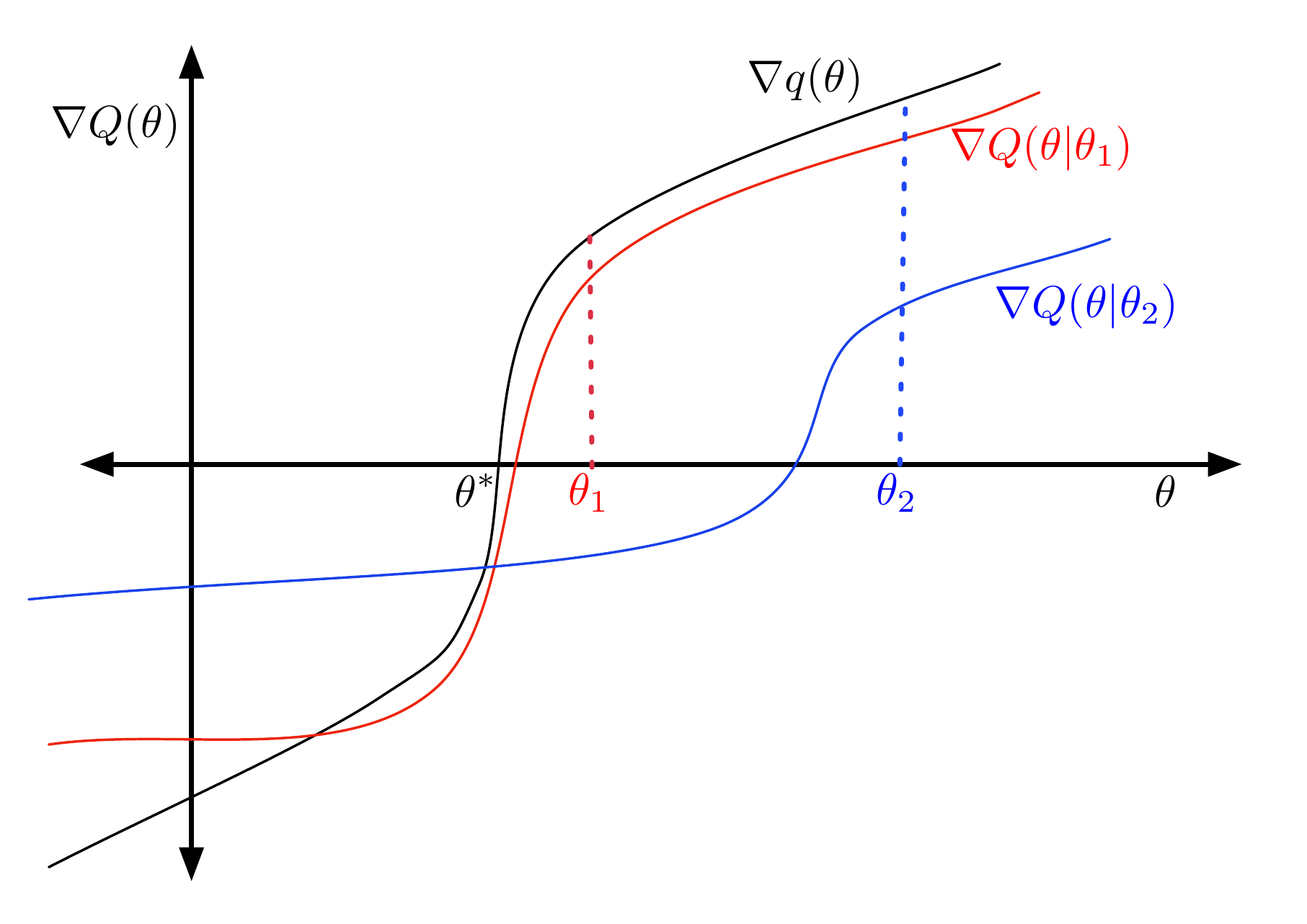}
\end{center}
\caption{Illustration of the gradient stability
  condition~\eqref{EqnGSONE}: for a point $\theta_1$ close to the
  population optimum $\thetastar$, the gradients $\nabla
  Q(\theta_1|\theta_1)$ and $\nabla q(\theta_1)$ must be close,
  whereas for a point $\theta_2$ distant from $\fp$, the gradients
  $\nabla Q(\theta_2| \theta_2)$ and $\nabla q(\theta_2)$ can be quite
  different.}
\label{FigGS}
\end{figure}}

\vspace*{.02in}

\noindent Under this condition, the following result guarantees local
contractivity of the gradient EM operator~\eqref{EqnGradEMTwo}:
\btheos
\label{ThmGradEM}
For some radius $r > 0$,
and a triplet $(\gsoneparam,\strongparam,\smoothparam)$
such that $0 \leq \gsoneparam < \strongparam \leq \smoothparam$,
suppose that the function \mbox{$\specq(\theta) = \qfun{\theta}{\fp}$}
is $\strongparam$-strongly concave~\eqref{EqnStronglyConcave},
$\smoothparam$-smooth~\eqref{EqnSmoothness}, and that the
GS$(\gsoneparam)$ condition~\eqref{EqnGSONE} holds on the ball
$\Balltwor$.  Then the population gradient EM operator $\PopGradEM$
with step size $\stepsize = \frac{2}{\smoothparam + \strongparam}$ is
contractive over $\Balltwor$, in particular with
\begin{align}
\norm{\mgrad{\pms} - \fp} & \leq \Big( 1 - \frac{2 \strongparam -
  2\gsoneparam}{\smoothparam + \strongparam} \Big) \; \norm{\pms -
  \fp} \qquad \mbox{for all $\pms \in \Balltwor$.}
\end{align}
\etheos
\noindent As an immediate consequence, under the conditions of the
theorem, for any initial point \mbox{$\pms^0 \in \Balltwor$,} the
population gradient EM sequence $\{\pms^t\}_{t=0}^\infty$ exhibits
linear convergence---viz.
\begin{align}
\norm{\pms^t - \fp} & \leq \Big( 1 - \frac{2 \strongparam -
  2\gsoneparam}{\smoothparam + \strongparam} \Big)^t \; \norm{\pms^0 -
  \fp} \qquad \mbox{for all $t = 1, 2, \ldots$.}
\end{align}
\begin{proof}
By definition of the gradient EM update~\eqref{EqnGradEMTwo}, we have
\begin{align*}
\norm{\mgrad{\pms} - \fp} & = \norm{\pms + \stepsize \nabla
  \QFUN{\pms}{\pms} - \fp} \\
& \stackrel{\mathrm{(i)}}{\leq} \underbrace{\norm{ \pms + \stepsize \nabla
    \qfun{\pms}{\fp} - \fp}}_{\|\OrGrad(\pms) - \fp \|_2} + \stepsize
\norm{ \nabla \QFUN{\pms}{\pms} - \nabla \qfun{\pms}{\fp}} \\
& \stackrel{\mathrm{(ii)}}{\leq} \Big( \frac{\smoothparam -
    \strongparam}{\smoothparam + \strongparam} \Big) \norm{\pms - \fp}
  + \stepsize \gsoneparam \norm{\pms - \fp}.
\end{align*}
where step (i) follows from the triangle inequality, and step (ii)
uses the contractivity of $\OrGrad$ from
equation~\eqref{EqnOrGradContract}, and condition~\gsonenospace.
Substituting $\stepsize = \frac{2}{\smoothparam + \strongparam}$ and
performing some algebra yields the claim.
\end{proof}


\subsubsection{Guarantees for sample-based gradient EM}

In this section, in parallel with our earlier analysis of sample-based
version of the EM algorithm, we analyze two sample-based variants of
the gradient EM algorithm, the first when the update operator
$\GremSamp$ is computed using all $\numobs$ samples and 
applied repeatedly, and the second based on sample-splitting.

We begin by introducing quantities that measure the deviations of the
sample operator $\GremSamp$ from the population version $\Grem$.  For
a given sample size $\numobs$ and tolerance parameter $\delta \in
(0,1)$, we let $\rategrad(\numobs, \delta)$ be the smallest scalar
such that, for any fixed vector $\pms \in \Balltwor$,
\begin{align}
\label{EqnRateGrad}
\norm{\GremSamp(\pms)- \Grem(\pms)} & \leq \rategrad(\numobs,\delta)
\end{align}
with probability at least $1-\delta$.  The uniform analogue of this
deviation is defined similarly: the quantity $\rategradunif(\numobs,
\delta)$ is the smallest scalar for which
\begin{align}
\label{EqnRateGradUniform}
\sup_{\pms \in \Balltwor} \norm{\GremSamp(\pms)- \Grem(\pms)} & \leq
\rategradunif(\numobs,\delta)
\end{align}
with probability at least $1 - \delta$.
\btheos
\label{ThmSampleGradEM}
Suppose that the population gradient EM operator $\Grem: \ParSpace
\rightarrow \ParSpace$ is contractive with parameter $\kappa \in
(0,1)$ on the ball $\Balltwor$, and the initial vector $\thetasamp^0$
belongs to $\Balltwor$.
\begin{enumerate}
\item[(a)] If the sample size $\numobs$ is large enough to ensure that
\begin{subequations}
\begin{align}
\label{EqnSamplingBoundGremUniform}
\rategradunif(\numobs,\delta) & \leq (1 - \kappa) r,
\end{align}
then the gradient EM iterates $\{\theta^t\}_{t=0}^\infty$ satisfy the
bound
\begin{align}
\label{EqnSampleGREMRateUniform}
\norm{\thetasamp^t - \fp} & \leq \kappa^t \norm{ \pms^0 - \fp} +
\frac{1}{1 - \kappa} \; \rategradunif (\numobs, \delta)
\end{align}
\end{subequations}
with probability at least $1 - \delta$.
\item[(b)] If the sample size $\numobs$ is large enough to ensure that
\begin{subequations}
\begin{align}
\label{EqnSamplingBoundGrem}
\HACKRATEGRAD & \leq (1 - \kappa) r,
\end{align}
then the sample-splitting gradient EM iterates
$\{\theta^t\}_{t=0}^\Tfinal$ based on $\frac{\numobs}{\Tfinal}$
samples per round satisfy the bound
\begin{align}
\label{EqnSampleGREMRate}
\norm{\thetasamp^t - \fp} & \leq \kappa^t \norm{ \pms^0 - \fp} +
\frac{1}{1 - \kappa} \; \rategrad \Big(\frac{\numobs}{\Tfinal},
\frac{\delta}{\Tfinal} \Big)
\end{align}
\end{subequations}
with probability at least $1 - \delta$.
\end{enumerate}
\etheos
Note that the guarantees~\eqref{EqnSampleGREMRateUniform} and
\eqref{EqnSampleGREMRate} are identical to the earlier
bounds~\eqref{EqnSampleEMRateUniform} and \eqref{EqnSampleEMRate} from
Theorem~\ref{ThmSampleEM}, modulo the replacements of $(\rateem,
\rateemunif)$ by $(\rategrad, \rategradunif)$.  We omit the proofs,
since they follow from essentially the same argument as
Theorem~\ref{ThmSampleEM}.  Thus, in order to obtain interpretable
bounds for gradient EM applied to specific models, it only remains to
establish the $\kappa$-contractivity of the population operator, and
to compute the functions $\rategrad$ or $\rategradunif$.

\subsubsection{Stochastic version of gradient EM}
\label{SecStochEMAnalysis}
In this section, we analyze a sample-based variant of gradient EM that
is inspired by stochastic approximation.  It can be viewed as an
extreme form of sample-splitting, in which we use only a single sample
per iteration, but compensate for the noisiness using a decaying step
size. Throughout this section we assume 
that (a lower bound  
on) the radius of convergence $r$ of the population operator is known
to the algorithm\footnote{This assumption can be restrictive in practice.
We believe the requirement can be eliminated by a more judicious
choice of the step-size parameter in the first few iterations.}.

In particular, given a sequence of positive step sizes
$\{\stepsize^t\}_{t=0}^\infty$, we analyze the recursion
\begin{align}
\label{EqnStochasticEM}
\theta^{t+1} & = \Pi \Big( \theta^t + \stepsize^t \nabla
\QFUNTHREE{\theta^t}{\theta^t}{1} \Big),
\end{align}
where the gradient $\nabla \QFUNTHREE{\theta^t}{\theta^t}{1}$ is
computed using a single fresh sample at each iteration.  Here $\Pi$
denotes the projection onto the Euclidean ball 
$\BalltwoTWO{\frac{r}{2}}{\theta^0}$
of radius $\frac{r}{2}$ centered at the initial iterate
$\theta^0$.  Thus, given any initial vector $\theta^0$ in the ball of
radius $r/2$ centered at $\thetastar$, we are guaranteed that all
iterates remain within an $r$-ball of $\thetastar$.  The following
result is stated in terms of the constant $\BRAZIL \defn \frac{2
  \smoothparam \strongparam}{\strongparam + \smoothparam} -
\gsoneparam > 0$, and the uniform variance $\sigmaSGD^2 \defn \sup
\limits_{\theta \in \Balltwor} \Exs \| \nabla
\QFUNTHREE{\theta}{\theta}{1}\|_2^2$.
\btheos
\label{ThmStochasticEM}
For a triplet $(\gsoneparam,\strongparam,\smoothparam)$
such that $0 \leq \gsoneparam < \strongparam \leq \smoothparam$,
suppose that the population function $\specq$ is
$\strongparam$-strongly concave~\eqref{EqnStronglyConcave},
$\smoothparam$-smooth~\eqref{EqnSmoothness}, and satisfies the
GS$(\gsoneparam)$ condition~\eqref{EqnGSONE} over the ball $\Balltwor$.
Then given an initialization $\theta^0 \in \BalltworONE{\frac{r}{2}}$, the
stochastic EM gradient updates~\eqref{EqnStochasticEM} with step size
$\stepsize^t \defn \frac{3}{2 \BRAZIL \, (t + 2) }$ satisfy the bound
\begin{align}
\label{EqnStochGradEM}
\Exs[\|\theta^{t} - \thetastar\|_2^2] & \leq \frac{9
  \sigmaSGD^2}{\BRAZIL^2} \, \frac{1}{(t+2)} + \Big(\frac{2}{t+2}
\Big)^{3/2} \: \|\theta^0 - \thetastar\|_2^2 \qquad \mbox{for
  iterations $t = 1, 2, \ldots$.}
\end{align}
\etheos 
\noindent While the stated claim~\eqref{EqnStochGradEM} provides
bounds in expectation, it is also possible to obtain high-probability
results.\footnote{Although we do not consider this extension here,
  stronger exponential concentration results follow from controlling
  the moment generating function of the random variable $\sup_{\theta
    \in \Balltwor} \| \nabla\QFUNTHREE{\theta}{\theta}{1}\|_2^2$.
  For instance, see \citet{nemirovski} for such results in
  the context of stochastic optimization.}

\begin{proof}
In order to prove this theorem we first
establish a recursion on the expected
mean-squared error. As with Theorem~\ref{ThmGradEM} this
result is established by relating the population
gradient EM operator to the gradient ascent 
operator on the function $q(\cdot)$.
This key recursion
along with some algebra will yield
the theorem.
\blems
\label{LKR}
Given the stochastic EM gradient iterates with step sizes
$\{\stepsize^t \}_{t=0}^\infty$, the error \mbox{$\Delta^{t+1} \defn
  \theta^{t+1} - \thetastar$} at iteration $t+1$ satisfies the recursion
\begin{align}
\label{EqnBrazilChile}
\Exs [ \|\Delta^{t+1}\|_2^2] & \leq \Big \{ 1 - \stepsize^t \BRAZIL
\Big \} \Exs [\|\Delta^t\|_2^2] + (\stepsize^t)^2 \sigmaSGD^2,
\end{align}
where $\sigmaSGD^2 = \sup
\limits_{\theta \in \Balltwor} \Exs[\| \nabla
  \QFUNTHREE{\theta}{\theta}{1}\|_2^2]$.
\elems
\noindent We prove this lemma in Appendix~\ref{AppStochEM}. 

\vspace*{.1in}

Using this result, we can now complete the proof of the
bound~\eqref{EqnStochGradEM}.  With the step size choice $\stepsize^t
\defn \frac{\HACKPAR}{\BRAZIL \, (t + 2) }$ where $\HACKPAR =
\frac{3}{2}$, unwrapping the recursion~\eqref{EqnBrazilChile} yields
\begin{align}
\label{EqnColombia}
\Exs[\|\Delta^{t+1}\|_2^2] & \leq \frac{\HACKPAR^2
  \sigmaSGD^2}{\BRAZIL^2} \sum_{\tau =2}^{t+1} \Big \{ \frac{1}{\tau^2}
\prod_{\ell =\tau+1}^{t+2} \Big(1 - \frac{\HACKPAR}{\ell} \Big) \Big
\} + \frac{\HACKPAR^2
  \sigmaSGD^2}{\BRAZIL^2 (t+2)^2} + \prod_{\ell =2}^{t+2} \Big(1 - \frac{\HACKPAR}{\ell} \Big) \;
\Exs[\|\Delta^0\|_2^2].
\end{align}
In order to bound these terms we use the following fact: 
For any $\HACKPAR \in (1,2)$, we have 
$$\prod_{\ell = \tau +1}^{t+2}
\Big(1 - \frac{\HACKPAR}{\ell} \Big) \leq \Big ( \frac{\tau+1}{t + 3}
\Big)^\HACKPAR.$$ 
\noindent See Noorshams and Wainwright~\cite{NooWai13a} for a proof.
Using this fact in Equation~\eqref{EqnColombia} yields
\begin{align*}
\Exs[\|\Delta^{t+1}\|_2^2] & \leq \frac{\HACKPAR^2 \sigmaSGD^2}{\BRAZIL^2
  \, (t+3)^\HACKPAR} \sum_{\tau =2}^{t+2} \frac{(\tau
  +1)^\HACKPAR}{\tau^2} + \Big(\frac{2}{t+3} \Big)^{\HACKPAR}
\Exs[\|\Delta^0\|_2^2] \\
& \leq \frac{2 \HACKPAR^2 \sigmaSGD^2}{\BRAZIL^2 \, (t+3)^\HACKPAR} \;
\sum_{\tau =2}^{t+2} \frac{1}{\tau^{2 - \HACKPAR}} +
\Big(\frac{2}{t+3} \Big)^{\HACKPAR} \: \Exs[\|\Delta^0\|_2^2].
\end{align*}
Finally, applying the integral upper bound $\sum \limits_{\tau
  =2}^{t+2} \frac{1}{\tau^{2 - \HACKPAR}} \leq \int_1^{t+2}
\frac{1}{x^{2-\HACKPAR}} dx \; \leq 2 (t+3)^{\HACKPAR - 1}$ yields the
claim~\eqref{EqnStochGradEM}.
\end{proof}

\noindent In order to obtain guarantees for stochastic gradient EM
applied to specific models, it only remains to prove the concavity and
smoothness properties of the population function $q$, and to bound the
uniform variance $\sigmaSGD$. \\


\paragraph{A summary:}
For the convenience of the reader, let us now summarize the theorems
given in this section, including the assumptions on which they rely
and the results that they provide. \\

\vspace*{.05in}

\begin{centering}
{\begin{small}
\begin{tabular}{|| l | l | c ||}
\hline \multicolumn{1}{||c}{Condition} & \multicolumn{1}{|c}{Result}&
\multicolumn{1}{|c||}{Thm.} \\ \hline 
Strong concavity of $q$ and FOS
& Pop. contractivity of EM {\color{red}(R1)} & Thm.~\ref{ThmEM}
\\
Bound on $\rateemunif$ and {\color{red}(R1)} & Fin.-sample bound for
EM & Thm.~\ref{ThmSampleEM} \\
Bound on $\rateem$ and {\color{red}(R1)} & Fin.-sample bound for
sample splitting EM & Thm.~\ref{ThmSampleEM} \\
Strong concavity, smoothness of $q$ and GS & Pop.  contractivity of
grad. EM {\color{red} (R2)} & Thm.~\ref{ThmGradEM} \\
Bound on $\rategradunif$ and {\color{red} (R2)} & Fin.-sample bound
for grad. EM & Thm.~\ref{ThmSampleGradEM} \\
Bound on $\rategrad$ and {\color{red} (R2)} & Fin.-sample bound for
sample splitting grad. EM & Thm.~\ref{ThmSampleGradEM} \\
Bound on $\sigmaSGD$ and {\color{red} (R2)} & Fin.-sample bound for
stochastic gradient EM & Thm.~\ref{ThmStochasticEM} \\ \hline
\end{tabular}
  \end{small}
}
\end{centering}


\section{Consequences for specific models}
\label{SecConsequences}

In the previous section, we provided a number of general theorems on
the behavior of the EM algorithm as well as the gradient EM algorithm,
at both the population and sample levels.  In this section, we develop
some concrete consequences of this general theory for the three
specific model classes previously introduced in
Section~\ref{SecIllustrative}.


\subsection{Gaussian mixture models}
\label{SecGMMAnalysis}

We begin by analyzing the EM updates for the Gaussian mixture model
previously introduced in Section~\ref{SecExaGMM}.  Our first result
(Corollary~\ref{CorMOGPop}) establishes contractivity for the
population operator~\eqref{EqnMemMOG}, whereas our second result
(Corollary~\ref{CorMOGUniform}) provides bounds for the sample-based
EM updates.

Recall that our mixture model consists of two equally weighted
components, with distributions $\NORMAL(\thetastar, \sigma^2 I)$ and
$\NORMAL(-\thetastar, \sigma^2 I)$ respectively.  The difficulty of
estimating this mixture model can be characterized by the
signal-to-noise ratio $\frac{\|\thetastar\|_2}{\sigma}$, and our
analysis requires a lower bound of the form
\begin{align}
\label{EqnMOGSNR}
\frac{\|\thetastar\|_2}{\sigma} & > \MOGSNR,
\end{align}
for a sufficiently large constant $\MOGSNR > 0$.  Past work by
\citet{redner84} provides evidence for the necessity of this
assumption: for Gaussian mixtures with low signal-to-noise ratio, they
show that the ML solution has large variance and furthermore verify
empirically that the convergence of the EM algorithm can be quite
slow.  Other researchers~\cite{xu_gmm1,xu_gmm2} also provide
theoretical justification for the slow convergence of EM on poorly
separated Gaussian mixtures. \\

\noindent With the signal-to-noise ratio lower bound $\MOGSNR$ defined
above we have the following guarantee:
\bcors[Population contractivity for Gaussian mixtures]
\label{CorMOGPop}
Consider a Gaussian mixture model for which the SNR
condition~\eqref{EqnMOGSNR} holds for a sufficiently large $\MOGSNR$.
Then there is a universal constant $c > 0$ such that
the population EM operator~\eqref{EqnMemMOG} is $\kappa$-contractive
over the ball $\Balltwor$ with
\begin{align}
r = \frac{\|\thetastar\|_2}{4}, \quad \mbox{and} \quad \kappa(\MOGSNR)
\leq \MYEXP{-c \MOGSNR^2}.
\end{align}
\ecors

\noindent This corollary guarantees that when the SNR is sufficiently
large, then the MLE $\thetastar$ has a basin of attraction that is at
least a constant fraction of the signal strength.  Moreover, the
convergence rate of the population updates is geometric, with the
contraction factor $\kappa$ decreasing exponentially in the
signal-to-noise ratio.  The proof of Corollary~\ref{CorMOGPop}
involves establishing that for a sufficiently large SNR, the strong
concavity and \gstwo($\gstwoparam$) conditions hold for a Gaussian
mixture model, so that Theorem~\ref{ThmEM} can be applied.  Although
the proof structure is conceptually straightforward, the details are
quite technical, so that we defer it to
Appendix~\ref{AppCorMOGPop}. \\

Based on the population-level contractivity guaranteed by
Corollary~\ref{CorMOGPop}, we can also establish guarantees for the
standard EM sequence $\theta^{t+1} = \MemSamp(\theta^t)$, where the
sample-based operator $\MemSamp$ was previously defined in
equation~\eqref{EqnMemSampMOG}.  This guarantee involves the function
\mbox{$\varphi(\sigma; \|\thetastar\|_2) \defn \norm{\thetastar}
  \sqrt{\norm{\thetastar}^2 + \sigma^2}$,} as well as positive universal
constants $(c, \plaincon_1, \plaincon_2)$.
\bcors[Sample-based EM guarantees for Gaussian mixtures]
\label{CorMOGUniform}
In addition to the conditions of Corollary~\ref{CorMOGPop}, suppose
that the sample size is lower bounded as $\numobs \geq \sivacon_1
\usedim \log(1/\delta)$.  Then given any initialization $\theta^0 \in
\BalltworONE{\frac{\|\thetastar\|_2}{4}}$, there is a contraction
coefficient $\kappa(\MOGSNR) \leq e^{-c \MOGSNR^2}$ such that the
standard EM iterates $\{\theta^t\}_{t=0}^\infty$ satisfy the bound
\begin{align}
\label{EqnMOGUniform}
\norm{\theta^t - \thetastar} & \leq \kappa^t \norm{\theta^0 -
  \thetastar} + \frac{\sivacon_2}{1 - \kappa} \varphi(\sigma;
\|\thetastar\|_2) \; \sqrt{ \frac{\usedim}{\numobs} \; \log(1/\delta)}
\end{align}
with probability at least $1 - \delta$.
\ecors 

\noindent See Appendix~\ref{AppMOGUniform} for the proof of this
result. In Appendix~\ref{AppCorMOGSam}, we also give guarantees for EM
with sample-splitting which achieves better dependence on
$\norm{\theta^*}$ and $\sigma$ with an easier proof at the cost of
additional logarithmic factors in sample complexity.

A related result of \citet{dasgupta} shows that when the SNR is
sufficiently high a modified EM algorithm, with an intermediate
pruning step, reaches a near-optimal solution in two iterations.  On
one hand, the SNR condition in our corollary is significantly weaker,
requiring only that it is larger than a fixed constant independent of
dimension (as opposed to scaling with $d$), but their theory is
developed for more general $k$-mixtures. \\
\afterpage{
\iftoggle{figs}{
\begin{figure}
\begin{center}
\begin{tabular}{ccc}
\widgraph{0.45 \textwidth}{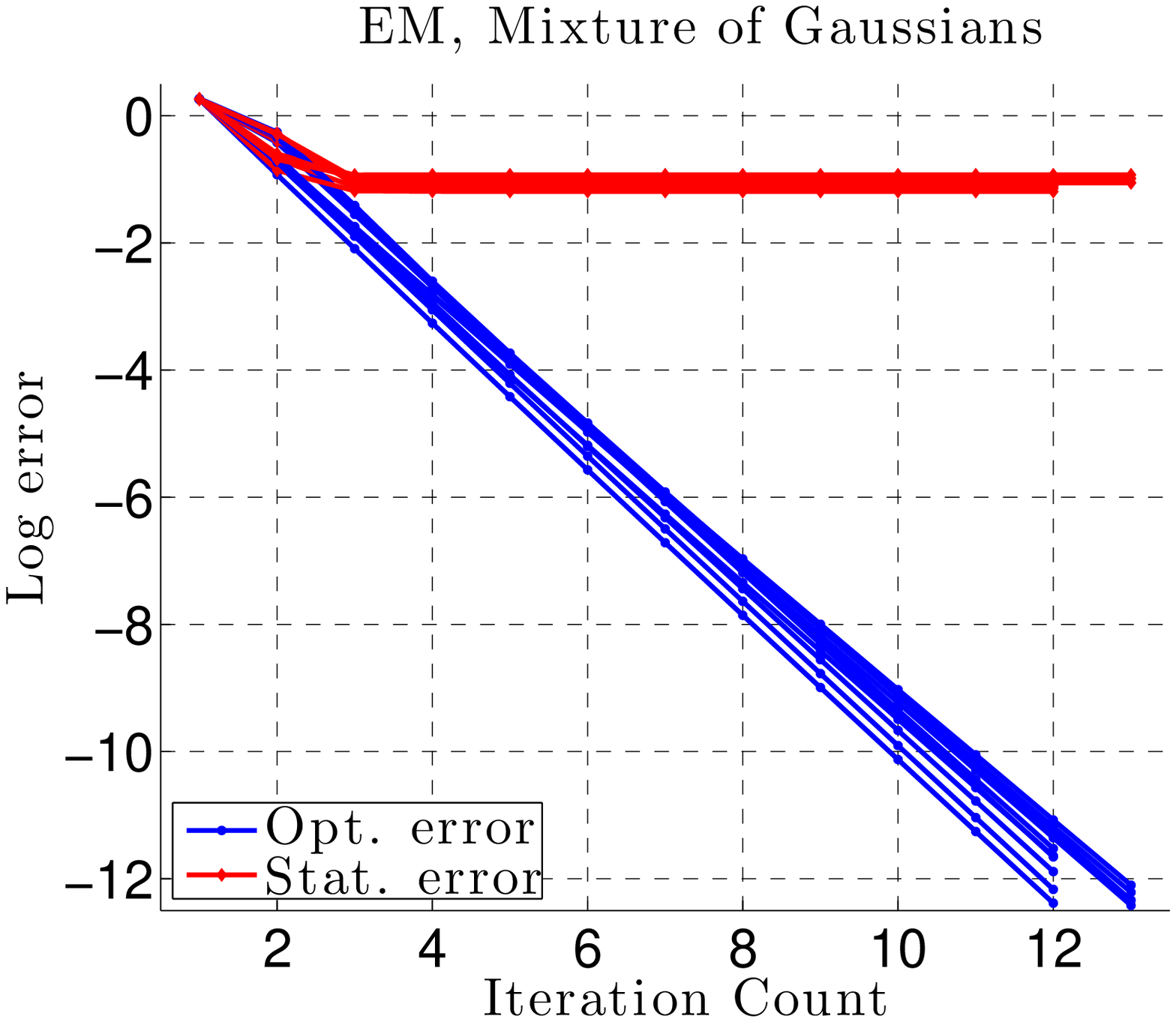} & &
 \widgraph{0.45 \textwidth}{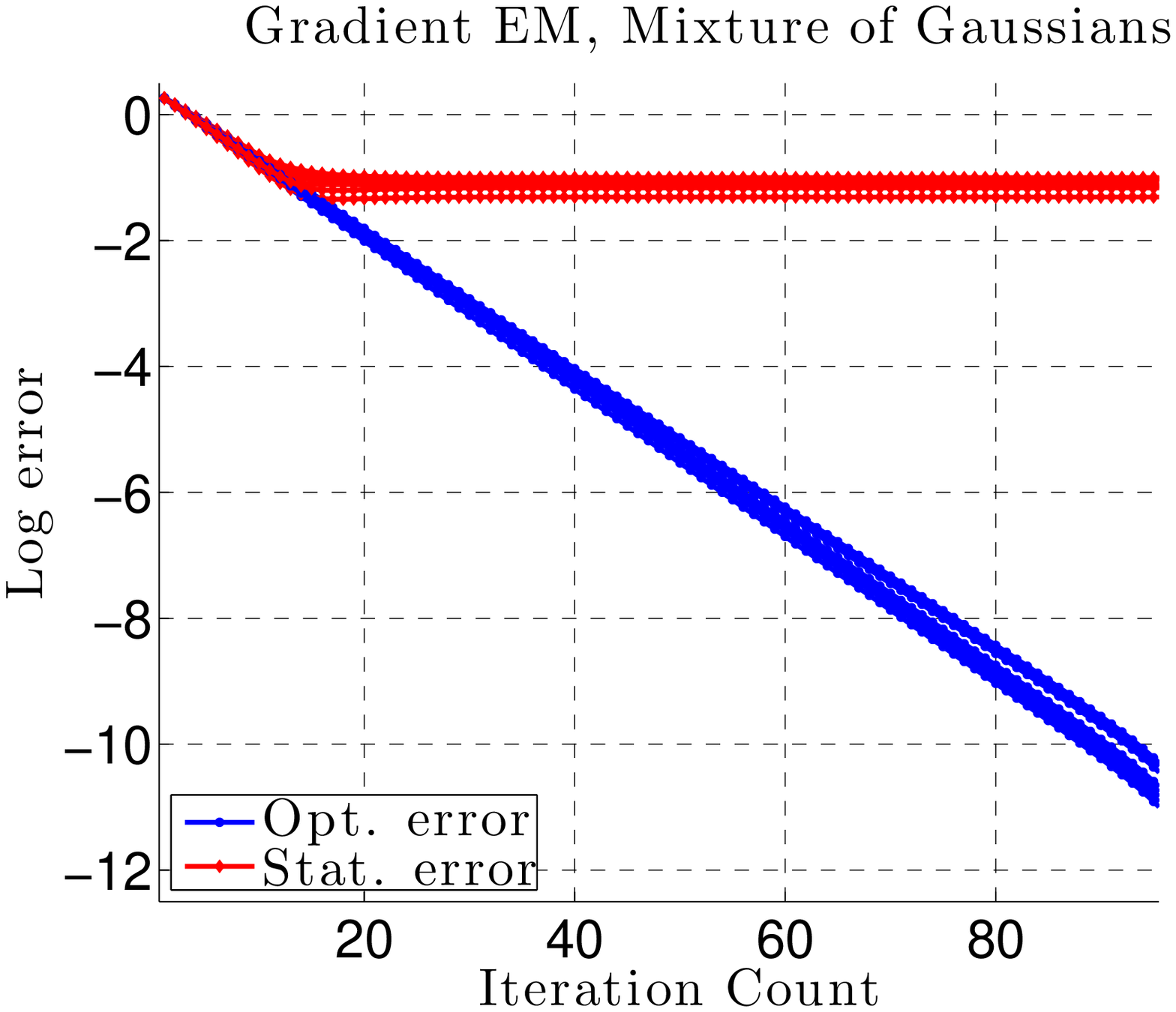} \\
(a) & & (b)
\end{tabular}
\end{center}
\caption{Plots of the iteration count versus log optimization error
  $\log(\norm{\theta^t - \thetahat})$ and log statistical error
  $\log(\norm{\theta^t - \fp})$.  (a) Results for the EM
  algorithm\protect\footnotemark.  (b) Results for the gradient EM
  algorithm.  Each plot shows $10$ different problem instances with
  dimension $\usedim = 10$, sample size $\numobs = 1000$ and
  signal-to-noise ratio $\frac{\norm{\thetastar}}{\sigma} = 2$.  The
  optimization error decays geometrically up to numerical precision,
  whereas the statistical error decays geometrically before leveling
  off.  }
\label{FigMOGIteration}
\end{figure}
\footnotetext{In this and subsequent figures we show simulations for
  the standard (i.e. not sample-splitting) versions of the EM and
  gradient EM algorithms.}}}

The bound~\eqref{EqnMOGUniform} provides a rough guide of how many
iterations are required: consider the smallest positive integer such
that
\begin{subequations}
\begin{align}
\label{EqnMOGSampleTcrit}
\Tcrit & \geq \log_{1/\kappa} \Big(\frac{\|\theta^0 - \thetastar\|_2(1 - \kappa)}{\varphi(\sigma;
    \|\thetastar\|_2)} \sqrt{ \frac{\numobs}{\usedim} \;
    \frac{1}{\log(1/\delta)}} \Big).
\end{align}
With this choice, we are guaranteed that the iterate $\theta^\Tcrit$
satisfies the bound
\begin{align}
\label{EqnMOGSampleBound}
\|\theta^T - \thetastar\|_2 & \leq \frac{(1 +
  \plaincon_2)\varphi(\sigma; \|\thetastar\|_2)}{1- \kappa} \;
\sqrt{\frac{\usedim}{\numobs} \; \log(1/\delta)}
\end{align}
\end{subequations}
with probability at least $1 - \delta$.  Treating $\sigma$ and
$\norm{\theta^*}$ as fixed there is no point in performing additional
iterations, since by standard minimax results, any estimator of
$\thetastar$ based on $\numobs$ samples must have $\ell_2$-error of
the order $\sqrt{\frac{\usedim}{\numobs}}$.  Of course, the iteration
choice~\eqref{EqnMOGSampleTcrit} is not computable based only on data,
since it depends on unknown quantities such as $\thetastar$ and the
contraction coefficient $\kappa$.  However, as a rough guideline, it
suggests that the iteration complexity should grow logarithmically in
the ratio $\numobs/\usedim$.

Corollary~\ref{CorMOGUniform} makes a number of qualitative
predictions that can be tested.  To begin, it predicts that the
statistical error $\|\theta^t - \thetastar\|_2$ should decrease
geometrically, and then level off at a plateau.
Figure~\ref{FigMOGSNR} shows the results of simulations designed to
test this prediction: for dimension $\usedim = 10$ and sample size
$\numobs = 1000$, we performed $10$ trials with the standard EM
updates applied to Gaussian mixture models with SNR
$\frac{\|\thetastar\|_2}{\sigma} = 2 $.  In panel (a), the red curves
plot the log statistical error versus the iteration number, whereas
the blue curves show the log optimization error versus iteration.  As
can be seen by the red curves, the statistical error decreases
geometrically before leveling off at a plateau.  On the other hand,
the optimization error decreases geometrically to numerical tolerance.
Panel (b) shows that the gradient EM updates have a qualitatively
similar behavior for this model, although the overall convergence rate
appears to be slower.

\iftoggle{figs}{
\begin{figure}
\begin{center}
\widgraph{0.5\textwidth}{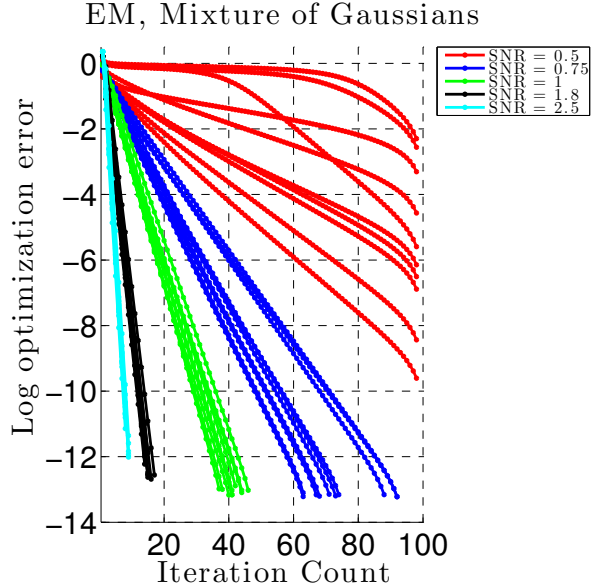}
\end{center}
\caption{Plot of the iteration count versus the (log) optimization
  error $\log(\norm{\theta^t - \thetahat})$ for different values of
  the SNR $\frac{\norm{\thetastar}}{\sigma}$.  For each SNR, we
  performed $10$ independent trials of a Gaussian mixture model with
  dimension $d = 10$ and sample size $\numobs = 1000$.  Larger values
  of SNR lead to faster convergence rates, consistent with
  Corollary~\ref{CorMOGUniform}.}
\label{FigMOGSNR}
\end{figure}}

In conjunction with Corollary~\ref{CorMOGPop},
Corollary~\ref{CorMOGUniform} also predicts that the convergence rate
should increase as the signal-to-noise ratio
$\frac{\|\thetastar\|_2}{\sigma}$ is increased.
Figure~\ref{FigMOGSNR} shows the results of simulations designed to
test this prediction: again, for mixture models with dimension
$\usedim = 10$ and sample size $\numobs = 1000$, we applied the
standard EM updates to Gaussian mixture models with varying SNR
$\frac{\|\thetastar\|_2}{\sigma}$.  For each choice of SNR, we
performed $10$ trials, and plotted the log optimization error $\log
\|\theta^t - \thetahat\|_2$ versus the iteration number.  As expected,
the convergence rate is geometric (linear on this logarithmic scale),
and the rate of convergence increases as the SNR grows\footnote{To be
  clear, Corollary~\ref{CorMOGUniform} predicts geometric convergence
  of the statistical error $\|\theta^t - \thetastar\|_2$, whereas
  these plots show the optimization error $\|\theta^t -
  \thetahat\|_2$.  However, the analysis underlying
  Corollary~\ref{CorMOGUniform} can also be used to show geometric
  convergence of the optimization error.}.
%


\subsection{Mixtures of regressions}
\label{SecMORAnalysis}

In this section, we analyze the EM and gradient EM algorithms for the
mixture of regressions (MOR) model, previously introduced in
Section~\ref{SecExaMOR}.  As in our analysis of the Gaussian mixture
model, our theory applies when the signal-to-noise ratio is
sufficiently large, as enforced by a condition of the form
\begin{align}
\label{EqnMORSNR}
\frac{\|\thetastar\|_2}{\sigma} & > \MORSNR
\end{align}
Under a suitable lower bound on this quantity, our first result
guarantees that the population level operators~\eqref{EqnMemMOR} 
and~\eqref{EqnGradientMOR}
are 
locally contractive.
\bcors[Population contractivity for MOR]
\label{CorMORPop}
Consider any mixture of regressions model satisfying the SNR
condition~\eqref{EqnMORSNR} for a sufficiently large constant
$\MORSNR$. Then the population EM operator $M$ from
 equation~\eqref{EqnMemMOR} and the population gradient EM
 operator $G$ from equation~\eqref{EqnGradientMOR} are 
 $\kappa$-contractive over the ball
  $\Balltwor$ with
\begin{align}
\label{EqnMORRadius}
r = \frac{\|\thetastar\|_2}{\deltaconstmor}, \quad
\mbox{and} \quad \kappa \leq \frac{1}{2}.
\end{align}
\ecors
\noindent 
As shown in the proof, the contraction coefficient $\kappa$ is again a
decreasing function of the SNR parameter $\MORSNR$.  However, its
functional form is not as explicit as in the Gaussian mixture case.
The proof of Corollary~\ref{CorMORPop} involves verifying that the
function $q$ for the MOR model satisfies the required concavity,
smoothness,
\gsonenospace($\gsoneparam$) and 
\gstwonospace($\gstwoparam$) conditions.  It is quite technically involved,
so that we defer it to Appendix~\ref{AppCorMORPop}.\\

Let us now provide guarantees for a sample-splitting version of the EM
updates.  Recall that sample-based EM operator was previously defined
in equation~\eqref{EqnMemSampMOR}.  For a given sample size $\numobs$
and iteration number $\Tcrit$, suppose that we split\footnote{To
  simplify exposition, assume that $\numobs/\Tcrit$ is an integer.}
our full data set into $\Tcrit$ subsets, each of size
$\numobs/\Tcrit$.  We then generate the sequence $\theta^{t+1} =
\Mem_{\numobs/\Tcrit}(\theta^t)$, where we use a fresh subset at each
iteration.  In the following result, we use
$\varphi(\sigma;\norm{\theta^*}) = \sqrt{ \sigma^2 +
  \norm{\theta^*}^2}$, along with positive universal constants
($\plaincon_1$,$\plaincon_2$).
\bcors[Sample-splitting EM guarantees for MOR]
\label{CorMORSamp}
In addition to the conditions of Corollary~\ref{CorMORPop}, suppose
that the sample size is lower bounded as $\numobs \geq \sivacon_1
\usedim \log(T/\delta)$.  Then there is a contraction coefficient
$\kappa \leq 1/2$ such that, for any initial vector $\theta^0 \in
\BalltworONE{\frac{\|\thetastar\|_2}{32}}$, the \mbox{sample-splitting} EM
iterates $\{\theta^t\}_{t=1}^\Tcrit$ based on $\numobs/\Tcrit$ samples
per step satisfy the bound
\begin{align}
\label{EqnMORSamp}
\norm{\theta^t - \thetastar} \leq \kappa^t \norm{\theta^0 -
  \thetastar} + \plaincon_2 \varphi(\sigma;\norm{\fp})
\sqrt{ \frac{\usedim}{\numobs} \; \Tcrit \log (\Tcrit/\delta)}
\end{align}
with probability at least $1 - \delta$.
\ecors
\noindent We prove this corollary in Appendix~\ref{AppCorMORSamp}.
Note the bound~\eqref{EqnMORSamp} again provides guidance on the
number of iterations to perform. For a given sample size $\numobs$,
suppose we perform $\Tcrit = c \log (\numobs/ d\varphi^2(\sigma;
\norm{\theta^*}))$ iterations for a constant $c$.  The
bound~\eqref{EqnMORSamp} then implies that
\begin{align}
\norm{\theta^\Tcrit - \thetastar} & \leq \plaincon_3 \varphi(\sigma;
\norm{\theta^*})\sqrt{ \frac{\usedim}{\numobs} \;
  \log^2\Big(\frac{\numobs}{d \varphi^2(\sigma; \norm{\theta^*})}\Big)
  \; \log(1/\delta)}
\end{align}
with probability at least $1 - \delta$.  Apart from the logarithmic
penalty $\log^2\big(\frac{\numobs}{d \varphi^2(\sigma;
  \norm{\theta^*})}\big)$, this guarantee matches the minimax rate for
estimation of a $\usedim$-dimensional regression vector.  We note that
the logarithmic penalty can be removed by instead analyzing the
standard form of the EM updates, as we did for the Gaussian mixture
model. \\

We conclude our discussion of the MOR model by stating a result for
the stochastic form of gradient EM analyzed in
Theorem~\ref{ThmStochasticEM}.  In particular, given a data set of
size $\numobs$, we run the algorithm for $\numobs$ iterations, with a
step size $\stepsize^t \defn \frac{3}{t+2}$ for iterations $t = 1,
\ldots, \numobs$. Once again our result is terms of
$\varphi(\sigma;\norm{\theta^*}) = \sqrt{ \sigma^2 +
  \norm{\theta^*}^2}$ and positive universal constants
($\plaincon_1$,$\plaincon_2$).

\bcors[Stochastic gradient EM guarantees for MOR]
\label{CorMORSGD}
In addition to the conditions of Corollary~\ref{CorMORPop}, suppose
that the sample size is lower bounded as $\numobs \geq \sivacon_1
\usedim \log(1/\delta)$.  Then given any initialization $\theta^0 \in
\BalltworONE{\frac{\|\thetastar\|_2}{32}}$, performing $\numobs$
iterations of the stochastic gradient EM gradient
updates~\eqref{EqnStochasticEM} yields an estimate $\thetahat =
\theta^\numobs$ such that
\begin{align}
\Exs[\|\thetahat - \thetastar\|_2^2] & \leq c_2 \, \varphi^2(\sigma;\norm{\theta^*}) \, \frac{\usedim}{\numobs}.
\end{align}
\ecors
\noindent We prove this corollary in Appendix~\ref{AppCorMORSGD}.
Figure~\ref{FigMORSGD} illustrates this corollary showing the error
as a function of iteration number (sample size) for the stochastic 
gradient EM
algorithm.

\iftoggle{figs}{
\begin{figure}[h]
\begin{center}
\begin{tabular}{cc}
\widgraph{0.5\textwidth}{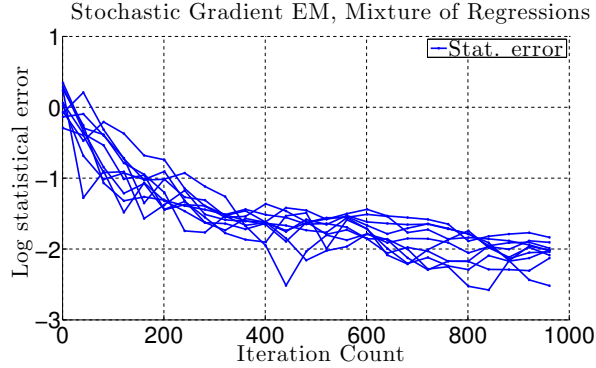}
\end{tabular}
\end{center}
\caption{A plot of the (log) statistical error for the stochastic gradient EM
algorithm as a function of iteration number (sample size) for the mixture of regressions
example. The plot shows 10 different problem instances 
with $d = 10$, $\frac{\norm{\theta^*}}{\sigma} = 2$ and 
$\frac{\norm{\theta^0 - \theta^*}}{\sigma} = 1.$
The statistical error decays at the sub-linear rate $\order(1/\sqrt{t})$ as a function
of the iteration number $t$.
An iteration of stochastic gradient EM is however
typically much faster and uses only a single sample.}
\label{FigMORSGD}
\end{figure}}

\iftoggle{figs}{
\begin{figure}[h]
\begin{center}
\begin{tabular}{cc}
\widgraph{0.45\textwidth}{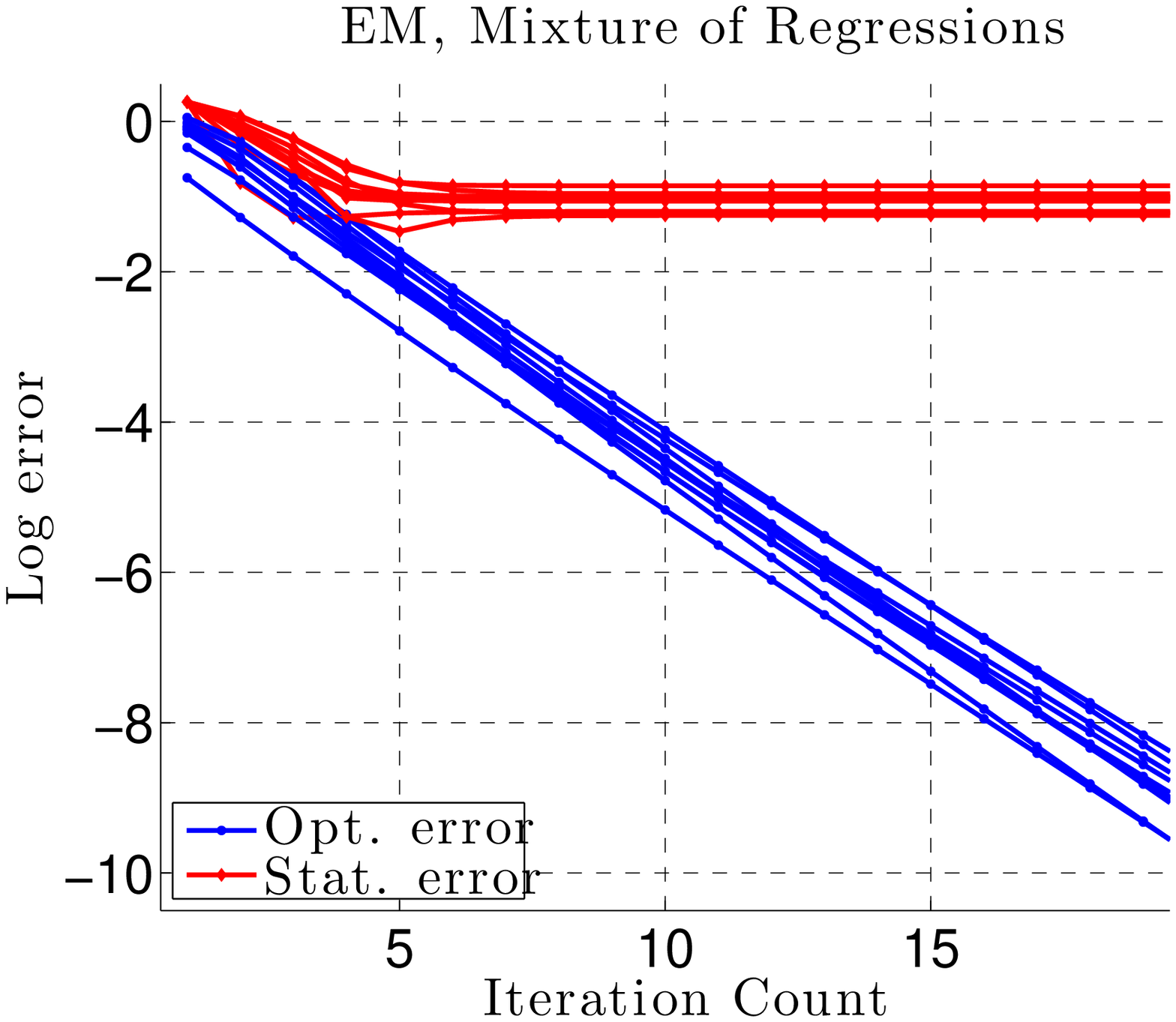} & 
\widgraph{0.45\textwidth}{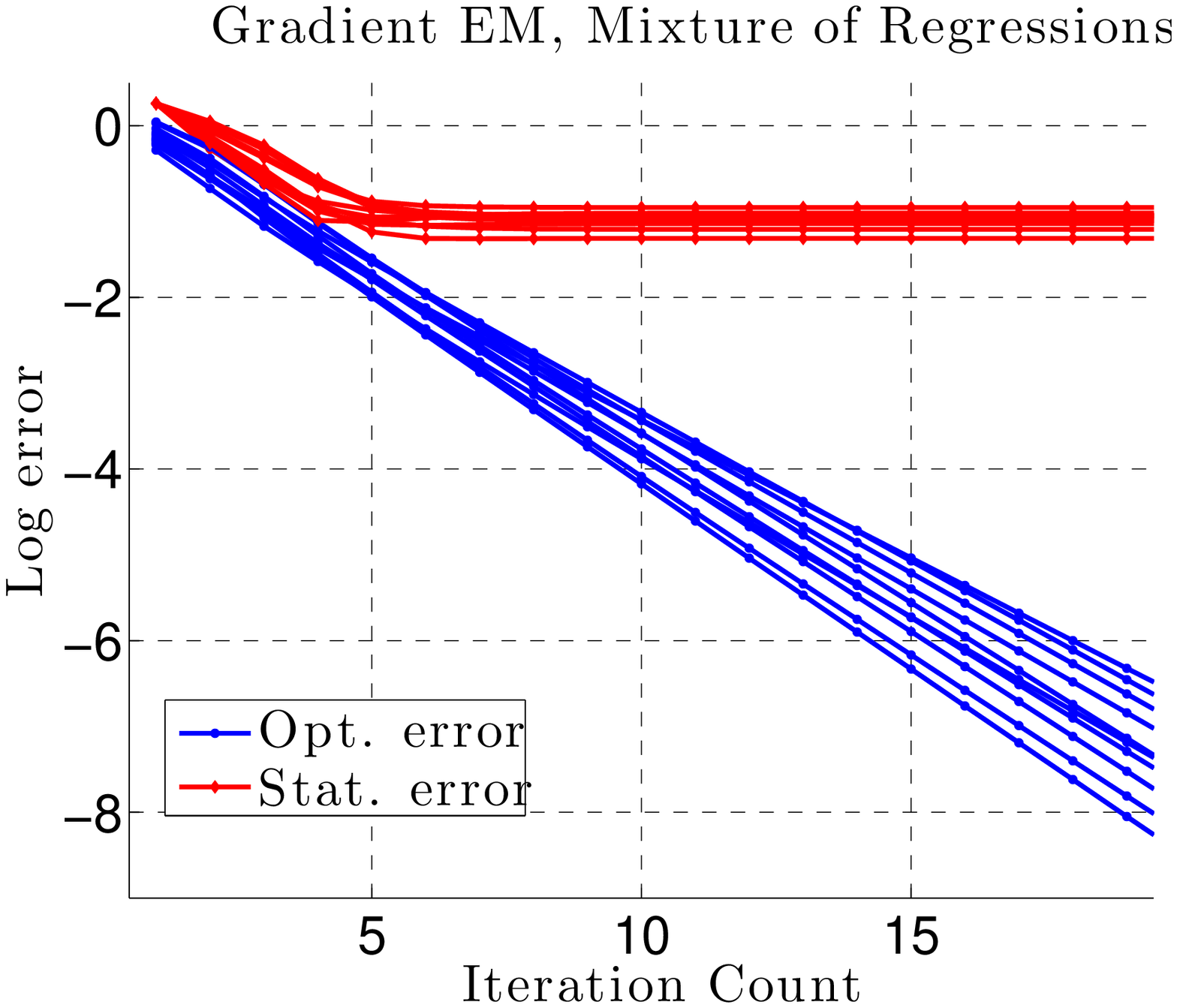}
\end{tabular}
\end{center}
\caption{Plots of the iteration count versus log optimization error
  $\log(\norm{\theta^t - \thetahat})$ and log statistical error
  $\log(\norm{\theta^t - \fp})$ for mixture of regressions. (a) Results
  for the EM algorithm.  (b) Results for the gradient EM algorithm.
  Each plot shows $10$ different problem instances of dimension
  $\usedim = 10$, sample size $\numobs = 1000$, and signal-to-noise
  ratio $\frac{\norm{\thetastar}}{\sigma} = 2$.  In both plots, the
  optimization error decays geometrically while the statistical error
  decays geometrically before leveling off.}
\end{figure}}


\subsection{Linear regression with missing covariates}
\label{SecMissingAnalysis}

This section is devoted to analysis of the gradient EM algorithm for
the problem of linear regression with missing covariates, as
previously introduced in Section~\ref{SecExaMissing}.  Here the
central parameter is the probability $\missing$ that any given
coordinate of the covariate vector is missing, and our analysis links
this quantity to the signal-to-noise ratio and the radius of
contractivity.  
Define $\ccon$ and $\ctwo$ to be such that the following bounds
hold,
\begin{align}
\label{EqnMissingConditions}
\frac{\|\thetastar\|_2}{\sigma} \leq \ccon, \quad \mbox{and} \quad
\|\theta - \thetastar\|_2 \leq r \defn \ctwo \sigma.
\end{align}
For any given choice of $(\ccon, \ctwo)$ define 
$\ncon \defn (\ccon + \ctwo)^2.$
Our guarantees apply
whenever the missing probability is bounded as
\begin{align}
\label{EqnMissingBound}
\missing < \frac{1}{1 + 2\ncon(1+\ncon)}.
\end{align}

\bcors[Population contractivity for missing covariates]
\label{CorMissingPop}
Given any missing covariate regression model with missing probability
$\missing$ satisfying the bound~\eqref{EqnMissingBound}, the gradient
EM operator~\eqref{EqnMemMissing} is $\kappa$-contractive over the
ball $\Balltwor$ with
\begin{align}
r = \ctwo \sigma, \quad \mbox{and} \quad \kappa = \frac{\ncon +
  \missing(1 + 2\ncon(1+\ncon))}{1+\ncon} < 1.
\end{align}
\ecors 
\noindent 
See Appendix~\ref{AppCorMissingPop} for the proof of
Corollary~\ref{CorMissingPop}.  Relative to our previous results, this
corollary is somewhat unusual, in that we require an \emph{upper
  bound} on the ratio $\frac{\|\thetastar\|_2}{\sigma}$.
Although this requirement might seem counter-intuitive at first sight,
known minimax lower bounds on regression with missing
covariates~\cite{loh_isit} show that it is unavoidable--- that is, it
is not an artifact of our analysis nor of the gradient EM algorithm.
Roughly these lower bounds formalize the intuition that 
as the norm
$\norm{\thetastar}$ increases, the amount of missing information
increases in proportion to the amount of observed information.
Figure~\ref{FigROC} provides the results of simulations that confirm
this behavior, in particular showing that for regression with missing
data, the radius of convergence eventually decreases as $\norm{\thetastar}$ grows.

\comment{\footnote{For the Gaussian mixture and
  mixture of regression models, we estimated the radius of convergence
  by sampling 100 points near the direction of $\theta^*$ at the
  specified radius. For the case of regression with missing data, we
  used 100 random initializations at the specified radius.  In each
  case, we checked convergence to a point sufficiently close to
  $\theta^*$. } }

\iftoggle{figs}{
\begin{figure}
\begin{center}
\widgraph{0.5\textwidth}{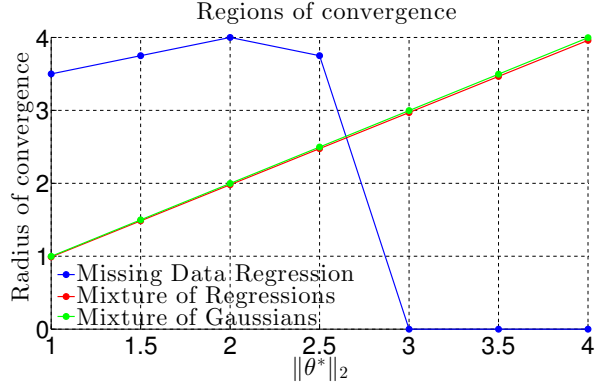}
\end{center}
\caption{Simulations of the radius of convergence for problems of
  dimension $\usedim = 10$, sample size $\numobs = 1000$, and variance
  $\sigma^2 = 1$.  Radius of convergence is defined as the maximum
  value of $\norm{\theta^0 - \fp}$ for which initialization at
  $\theta^0$ leads to convergence to an optimum near $\fp$.
  Consistent with the theory, for both the Gaussian mixture and
  mixture of regression models, the radius of convergence grows with
  $\norm{\fp}$.  In contrast, in the missing data case (here with
  $\missing = 0.2$), increasing $\norm{\fp}$ can cause the EM
  algorithm to converge to bad local optima, which is consistent with
  the prediction of Corollary~\ref{CorMissingPop}.}
\label{FigROC}
\end{figure}}


Let us now provide guarantees for a sample-splitting version of the EM
updates, based on the sample-based EM operator in
equation~\eqref{EqnMemSampMissing}.  As usual, for a given sample size
$\numobs$ and iteration number $\Tcrit$, suppose that we split our
full data set into $\Tcrit$ subsets, each of size $\numobs/\Tcrit$.
We then generate the sequence $\theta^{t+1} =
\Mem_{\numobs/\Tcrit}(\theta^t)$, where we use a fresh subset at each
iteration.
\bcors[Sample-splitting EM guarantees for missing covariates]
\label{CorMissingSamp}
In addition to the conditions of Corollary~\ref{CorMissingPop},
suppose that the sample size is lower bounded as $\numobs \geq
\sivacon_1 \usedim \log(1/\delta)$.  Then there is a contraction
coefficient $\kappa < 1$ such that, for any initial vector $\theta^0
\in \BalltworONE{\ctwo \sigma}$, the \mbox{sample-splitting} EM iterates
$\{\theta^t\}_{t=1}^\Tcrit$ based on $\numobs/\Tcrit$ samples per
iteration satisfy the bound
\begin{align}
\label{EqnMissingSamp}
\norm{\theta^t - \thetastar} \leq \kappa^t \norm{\theta^0 -
  \thetastar}+ \frac{c_2 \, \sqrt{1 + \sigma^2}}{1 - \kappa} \;
\sqrt{\frac{\usedim}{\numobs} \; \Tcrit \log(\Tcrit/\delta)}
\end{align}
with probability at least $1-\delta$.  
\ecors

\noindent We prove this corollary in Appendix~\ref{AppCorMissingSamp}.
We note that the constant $c_2$ is a monotonic function of the
parameters $(\ccon, \ctwo)$, but does not otherwise depend on
$\numobs$, $\usedim$, $\sigma^2$ or other problem-dependent
parameters.

As with Corollary~\ref{CorMORSamp}, this result provides guidance on
the appropriate number of iterations to perform: in particular, if we
set $\Tcrit = c \log \numobs$ for a sufficiently large constant $c$,
then the bound~\eqref{EqnMissingSamp} implies that
\begin{align*}
\norm{\theta^\Tcrit - \thetastar} & \leq c' \, \sqrt{1 + \sigma^2}
\sqrt{\frac{\usedim}{\numobs} \; \log^2(\numobs/\delta)}
\end{align*}
with probability at least $1-\delta$.  Modulo the logarithmic penalty
in $\numobs$, incurred due to the sample-splitting, this estimate
achieves the optimal $\sqrt{\frac{\usedim}{\numobs}}$ scaling of the
$\ell_2$-error. \\

\vspace*{.05in}

We conclude our discussion of the missing covariates model by stating
a result for the stochastic form of gradient EM analyzed in
Theorem~\ref{ThmStochasticEM}.  In particular, given a data set of
size $\numobs$, we run the algorithm for $\numobs$ iterations, with a
step size $\stepsize^t \defn \frac{3}{t+2}$ for iterations $t = 1,
\ldots, \numobs$.

\bcors[Stochastic gradient EM guarantees for missing covariates]
\label{CorMissingSGD}
In addition to the conditions of Corollary~\ref{CorMissingPop},
suppose that the sample size is lower bounded as $\numobs \geq
\sivacon_1 \usedim \log(1/\delta)$.  Then given any initialization
$\theta^0 \in \BalltworONE{\ctwo \sigma}$, performing $\numobs$ iterations
of the stochastic EM gradient updates~\eqref{EqnStochasticEM} with
step sizes $\stepsize^t = \frac{3}{2 \, (1-\kappa) (t+2)}$ yields an
estimate $\thetahat = \theta^\numobs$ such that
\begin{align}
\Exs[\|\thetahat - \thetastar\|_2^2] & \leq c_2 (1 + \sigma^2) \;
\frac{\usedim}{\numobs}.
\end{align}
\ecors
\noindent We prove this corollary in Appendix~\ref{AppCorMissingSGD}.
Figure~\ref{FigMISSGD} illustrates this.
 \iftoggle{figs}{
\begin{figure}[h]
\begin{center}
\begin{tabular}{c}
\widgraph{0.53\textwidth}{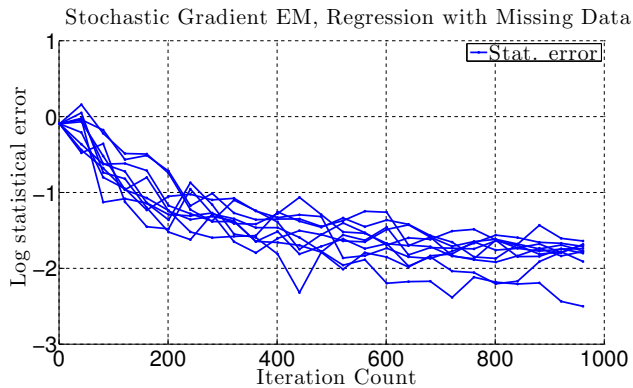}
\end{tabular}
\end{center}
\caption{
A plot of the (log) statistical error for the stochastic gradient EM
algorithm as a function of iteration number (sample size) for 
the problem of linear regression
with missing covariates. The plot shows 10 different problem instances 
with $d = 10$, $\frac{\norm{\theta^*}}{\sigma} = 2$ and 
$\frac{\norm{\theta^0 - \theta^*}}{\sigma} = 1.$
The statistical error decays at the sub-linear rate $\order(1/\sqrt{t})$ as a function
of the iteration number $t$.}
\label{FigMISSGD}
\end{figure}}

\iftoggle{figs}{
\begin{figure}[h]
\begin{center}
\begin{tabular}{ccc}
\widgraph{0.45 \textwidth}{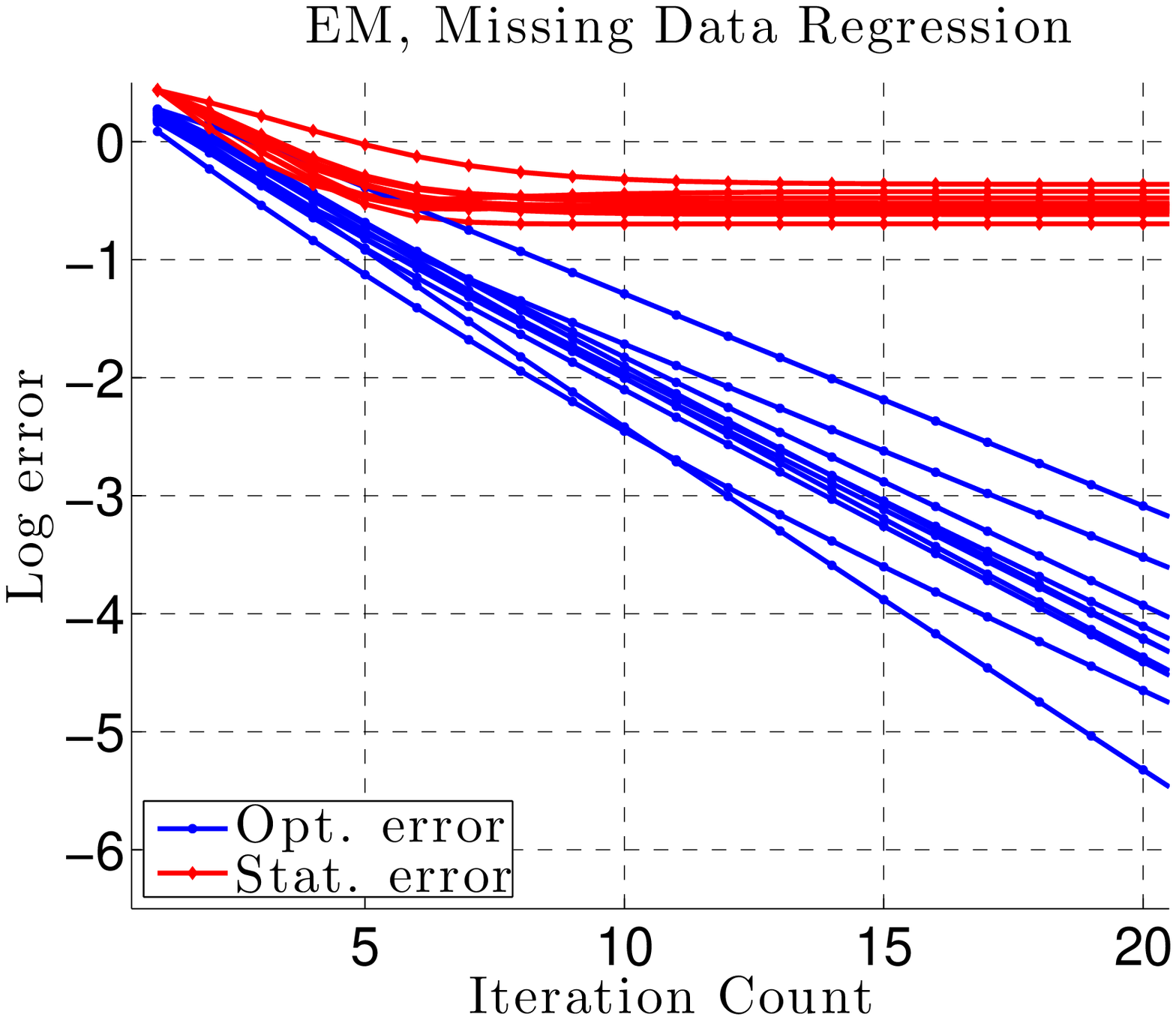} & &
\widgraph{0.45\textwidth}{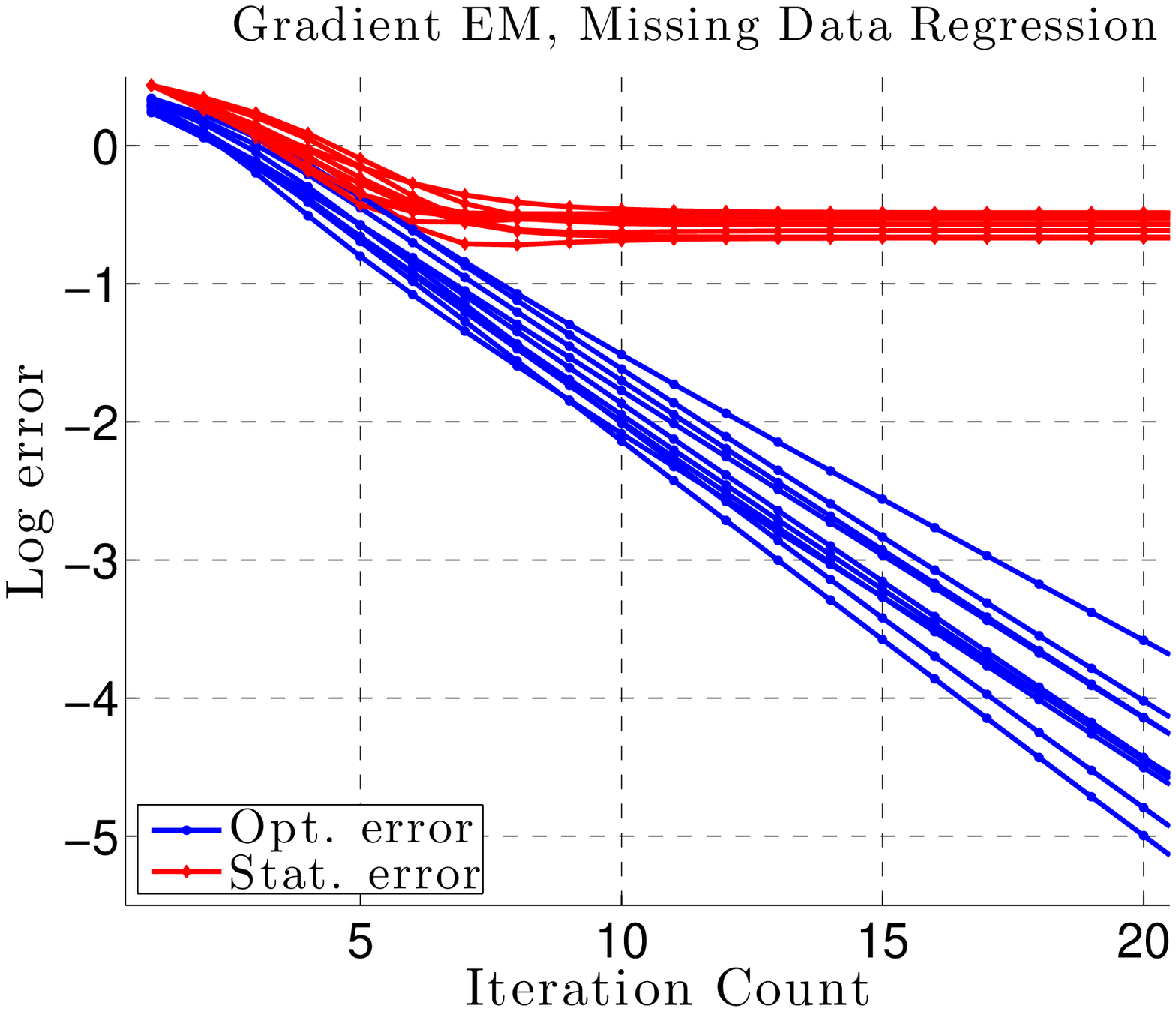} \\
(a) & & (b)
\end{tabular}
\end{center}
\caption{Plots of the iteration count versus log optimization error
  $\log(\norm{\theta^t - \thetahat})$ and log statistical error
  $\log(\norm{\theta^t - \fp})$ for regression with missing
  covariates. (a) Results for the EM algorithm.  (b) Results for the
  gradient EM algorithm.  Each plot shows $10$ different problem
  instances of dimension $\usedim = 10$, sample size $\numobs = 1000$,
  signal-to-noise ratio $\frac{\norm{\thetastar}}{\sigma} = 2$, and
  missing probability $\missing = 0.2$.  In both plots, the
  optimization error decays geometrically while the statistical error
  decays geometrically before leveling off.}
\label{FigOldMissing}
\end{figure}}


\section{Discussion}
\label{SecDiscussion}

In this paper, we have provided some general techniques for studying
the EM and gradient EM algorithms, at both the population and
finite-sample levels.
Although this paper focuses on these specific
algorithms, we expect that the techniques could be useful in
understanding the convergence behavior of other algorithms for
potentially non-convex problems.

The analysis of this paper can be extended in various directions.  For
instance, in the three concrete models that we treated, we assumed
that the model was correctly specified, and that the samples were
drawn in an i.i.d. manner, both conditions that may be violated in
statistical practice.  Maximum likelihood estimation is
known to have various robustness properties under model
mis-specification. Developing an understanding of the EM algorithm in
this setting is an important open problem.

Finally, we note that in concrete examples our 
analysis guarantees good behavior of the EM
and gradient EM algorithms when they are given suitable
initialization.  For the three model classes treated in this paper,
simple pilot estimators can be used to obtain such
initializations---in particular using PCA for Gaussian mixtures and
mixtures of regressions (e.g.,~\cite{hardem_mor}), and the plug-in
principle for regression with missing data
(e.g.,~\cite{ItuEtal99,XuYou07}).  These estimators can be seen as
particular instantiations of the method of moments
\cite{pearson94}. Although still an active area of research, a line of
recent work (e.g.,~\cite{tensor,hsu_mog,spectral_mor,Dictionary}) has
demonstrated the utility of moment-based estimators or initializations
for other types of latent variable models, and it would be interesting
to analyze the behavior of EM for such models.


\subsection*{Acknowledgments}  
This research was partially supported by ONR-MURI grant
N00014-11-1-0688 and NSF grant CIF-31712-23800 to MJW, and by US NSF
grants DMS-1107000, 
CDS\&E-MSS 1228246, ARO grant W911NF-11-1-0114, the
Center for Science of Information (CSoI), and US NSF Science and
Technology Center, under grant agreement CCF-0939370.  SB would also
like to thank John Duchi for helpful discussions.

\appendix


\section{Proofs for stochastic gradient EM} 
\label{AppStochEM}
In this section we provide proofs of results related to 
Theorem~\ref{ThmStochasticEM}
from Section~\ref{SecStochEMAnalysis}.
It only remains to prove Lemma~\ref{LKR}.

In order to establish Lemma~\ref{LKR} we require
an analogue of Theorem~\ref{ThmGradEM} that
allows for a wider range of step sizes.
Recall the classical gradient ascent operator on the function
$q(\theta) = \QFUN{\theta}{\thetastar}$.  For step size $\stepsize >
0$, it takes the form $T(\theta) = \theta + \stepsize \nabla
q(\theta)$.  Under the
stated $\strongparam$-concavity and $\smoothparam$-smoothness
conditions, for any step size $0 < \alpha \leq \frac{2}{\strongparam +
  \smoothparam}$, the classical gradient operator $T$ is contractive
with parameter
\begin{align*}
\phi(\alpha) & = 1 - \frac{2 \stepsize \smoothparam
  \strongparam}{\smoothparam+ \strongparam}.
\end{align*}
This follows from the classical analysis of gradient descent
(e.g.,~\cite{bubeck_co,Bertsekas_nonlin,Nesterov}).  Using this fact,
we can prove the following about the population gradient EM operator:
\blems
\label{LemBrazil}
For any step size $0 < \stepsize \leq \frac{2}{\strongparam +
  \smoothparam}$, the population gradient EM operator $G: \ParSpace
\rightarrow \ParSpace$ is contractive with parameter
$\kappa(\stepsize) = 1 - \stepsize \BRAZIL$, where
\begin{align}
\label{EqnDefnKappaStepsize}
\BRAZIL & \defn \frac{2 \smoothparam \strongparam}{\strongparam +
  \smoothparam} - \gsoneparam.
\end{align}
\elems
\noindent We omit the proof, since it follows from a similar argument
to that of Theorem~\ref{ThmGradEM}.  With this preliminary in
place we can now begin the proof of Lemma~\ref{LKR}. 

\subsection{Proof of Lemma~\ref{LKR}}
Let us write $\theta^{t+1} = \Pi(\thetatil^{t+1})$, where
$\tilde{\theta}^{t+1} \defn \theta^t + \alpha^t \nabla
\QFUNTHREE{\theta^t}{\theta^t}{1}$ is the update vector prior to
projecting onto the ball $\BalltwoTWO{\frac{r}{2}}{\theta^0}$.  Defining
the difference vectors $\Delta^{t+1} \defn \theta^{t+1} - \thetastar$
and $\DelTil^{t+1} \defn \thetatil^{t+1} - \thetastar$, we have
\begin{align*}
\|\Delta^{t+1}\|_2^2 - \|\Delta^t\|_2^2 & \leq \norm{\DelTil^{t+1}}^2
- \norm{\Delta^t}^2 = \inprod{\thetatil^{t+1} -
  \theta^t}{\thetatil^{t+1} + \theta^t - 2 \thetastar}.
\end{align*}
Introducing the shorthand $\What(\theta) \defn \nabla
\QFUNTHREE{\theta}{\theta}{1}$, we have $\thetatil^{t+1} - \theta^t =
\stepsize^t \What(\theta)$, and hence
\begin{align*}
\|\Delta^{t+1}\|_2^2 - \|\Delta^t\|_2^2 & \leq \stepsize^t
\inprod{\What(\theta^t)}{\stepsize^t \What(\theta^t) + 2 (\theta^t -
  \thetastar)} \\
& = (\stepsize^{t})^2 \|\What(\theta^t)\|_2^2 + 2 \stepsize^t
\inprod{\What(\theta^t)}{\Delta^t}.
\end{align*}
Letting $\Fil_t$ denote the $\sigma$-field of events up to the random
variable $\theta^t$, note that 
\begin{align*}
\Exs[\What(\theta^t) \mid
    \Fil_t] = \Wop(\theta^t) 
    \defn \nabla Q(\theta^t|\theta^t).
\end{align*} 
Consequently, by iterated
expectations, we have
\begin{align}
\label{EqnGermany}
\Exs[\|\Delta^{t+1}\|_2^2] & \leq \Exs[\|\Delta^t\|_2^2] +
(\stepsize^t)^2 \Exs \|\What(\theta^t) \|_2^2 + 2 \stepsize^t \Exs
\Big [ \inprod{\Wop(\theta^t)}{\Delta^t} \Big].
\end{align}
Now since $\thetastar$ maximizes the function $q$ and $\theta^t$
belongs to $\BalltwoTWO{\frac{r}{2}}{\theta^0}$, we have
\begin{align*}
\inprod{\Wop(\thetastar)}{\Delta^t} & = \inprod{\nabla q(\thetastar)}{
  \Delta^t} \leq 0.
\end{align*}
Combining with our earlier inequality~\eqref{EqnGermany} yields
\begin{align*}
\Exs[\|\Delta^{t+1}\|_2^2] & \leq \Exs[\|\Delta^t\|_2^2] +
(\stepsize^t)^2 \Exs \|\What(\theta^t) \|_2^2 + 2 \stepsize^t \Exs \Big
    [ \inprod{\Wop(\theta^t) - \Wop(\thetastar)}{\Delta^t} \Big].
\end{align*}
Defining $\Grem^t(\theta^t) \defn \theta^t + \stepsize^t
\Wop(\theta^t)$, we see that
\begin{align*}
\stepsize^t \inprod{\Wop(\theta^t) - W(\thetastar)}{\Delta^t } & =
\inprod{\Grem^t(\theta^t) - G^t(\thetastar) - (\theta^t -
  \thetastar)}{\theta^t - \thetastar} \\ & = \inprod{G^t(\theta^t) -
  G^t(\thetastar)}{\theta^t - \thetastar} - \|\theta^t -
\thetastar\|_2^2 \\
& \stackrel{\mathrm{(i)}}{\leq} (\kappa(\stepsize^t) - 1) \|\theta^t - \thetastar\|_2^2 \\
& \stackrel{\mathrm{(ii)}}{=} -\stepsize^t \BRAZIL \; \norm{\Delta^t}^2,
\end{align*}
where step (i) uses the contractivity of $G^t$ established in Lemma~\ref{LemBrazil} and step (ii) 
uses the definition of $\BRAZIL$ from
equation~\eqref{EqnDefnKappaStepsize}.  Putting together the pieces
yields the claim~\eqref{EqnBrazilChile}.
%


\section{Proofs for Gaussian mixture models}

In this section, we provide proofs of results related to the Gaussian
mixture model, as presented in Section~\ref{SecGMMAnalysis}.  More
specifically, we first prove Corollary~\ref{CorMOGPop} on the
population level behavior, followed by the proof of
Corollary~\ref{CorMOGUniform} on the behavior of the standard
sample-based EM updates.


\subsection{Proof of Corollary~\ref{CorMOGPop}}
\label{AppCorMOGPop}

In order to apply Theorem~\ref{ThmEM}, we need to verify the
$\strongparam$-concavity condition~\eqref{EqnStronglyConcave}, and the
FOS($\gstwoparam$) condition~\eqref{EqnGSTWO} over the ball
$\Balltwor$. The population EM operator for the Gaussian mixture
model was previously defined in equation~\eqref{EqnMemMOG}.  The
update $\theta \mapsto \mem{\theta}$ is based on maximizing the
function
\begin{align*}
\QFUN{\theta'}{\theta} = - \frac{1}{2}\Exs \big[ w_{\theta}(Y)
  \norm{Y - \theta'}^2 + (1 - w_{\theta}(Y)) \norm{Y + \theta'}^2 \big] \qquad \mbox{over $\theta' \in
  \real^\usedim$.}
\end{align*}
Here the weighting function takes the form
\begin{align*}
w_\theta(y) & \defn \frac{ \exp \big( - \frac{\norm{\theta - y}^2}{2
    \sigma^2} \big) }{ \exp \big( - \frac{\norm{\theta - y}^2}{2
    \sigma^2} \big) + \exp \big( - \frac{\norm{\theta + y}^2}{2
    \sigma^2} \big) }.
\end{align*}
By inspection, the function $\specq(\theta') = \QFUN{\theta'}{\fp}$ is
strongly concave on $\real^\usedim$ with $\strongparam = 1$.  

It remains to verify the FOS($\gstwoparam$)
condition~\eqref{EqnGSTWO}.  The following auxiliary lemma is central
to the proof:
\blems
\label{LemMogMainTech}
Under the conditions of Corollary~\ref{CorMOGPop}, there is a constant
$\gamma \in (0,1)$ with \mbox{$\gamma \leq \exp (-c_2 \eta^2)$} such
that
\begin{align}
\label{EqnMogMainTech}
\norm{\Exs \big[ 2\Delta_w(Y) Y \big]} & \leq \gamma \, \norm{ \theta -
  \thetastar},
\end{align}
where \mbox{$\Delta_w(y) \defn w_{\theta}(y) - w_{\fp}(y)$.}
\elems 
Taking this result as given for the moment, let us now verify the
\gstwo condition~\eqref{EqnGSTWO}.  By symmetry, we have $\Exs
\big[w_{\pms}(Y)\big] = 1 - \Exs \big[w_{\pms}(Y)\big] = \frac{1}{2}$
for any $\pms \in \pmss$.  Using this fact, it suffices to show that
\begin{align*}
\norm{\Exs \big[2\Delta_w(Y) Y \big]}
& < \norm{\theta - \thetastar}.
\end{align*}
This follows immediately from Lemma~\ref{LemMogMainTech}.
Thus, the FOS condition holds when $\gamma < 1.$
The bound on the contraction parameter follows from the fact that
$\gamma \leq \exp (-c_2 \eta^2)$
and applying Theorem~\ref{ThmEM} yields
Corollary~\ref{CorMOGPop}.


\paragraph{Proof of Lemma~\ref{LemMogMainTech}:}
We now prove Lemma~\ref{LemMogMainTech}.  Our proof makes
use of the following elementary facts:\\
\begin{carlist}
\item For the function $f(t) = \frac{t^2}{\exp(\mu t)}$, we have
\begin{subequations}
\begin{align}
\label{Eqn:x_over_exp}
\sup_{t \in [0, \infty]} f(t) & = \frac{4}{ (e \, \mu)^2}, \qquad
\mbox{achieved at $\tstar = \frac{2}{\mu}$ and} \\
\label{Eqn:x_over_exp2}
\sup_{t \in [t^*, \infty]} f(t) & = f(t^*), \qquad
\mbox{for $\tstar \geq \frac{2}{\mu}$.}
\end{align}
\end{subequations}
\item
For the function $g(t) = \frac{1}{(\exp(t) + \exp(-t))^2}$, we have
\begin{subequations}
\begin{align} 
\label{Eqnone_over_expA}
g(t) & \leq \frac{1}{4} \quad \mbox{for all $t \in \real$, and} \\
\label{Eqnone_over_expB}
\sup_{t \in [\mu, \infty]} g(t) & \leq \frac{1}{(\exp(\mu) +
  \exp(-\mu))^2} \leq \exp (- 2 \mu), \qquad \mbox{valid for any $\mu
  \geq 0$.}
\end{align}
\end{subequations}
\end{carlist}

\vspace*{.04in}

\noindent With these preliminaries in place, we can now begin the
proof.  For each $u \in [0,1]$, define $\theta_u = \fp + u \Delta$,
where $\Delta \defn \theta - \fp$.  Taylor's theorem applied to the
function $\theta \mapsto w_\theta(Y)$, followed by expectations,
yields
\begin{align*}
\Exs \Big[ Y \, \big(w_{\theta}(Y) - w_{\fp}(Y) \big) \Big] & = 2
\int_0^1 \Exs \Biggr[\underbrace{\frac{Y Y^T} {\sigma^2 \big( \exp
      \big( - \frac{\inprod{\theta_u}{Y}}{\sigma^2} \big) + \exp \big(
      \frac{\inprod{\theta_u}{Y}}{\sigma^2} \big) \big)^2}}_{\Gamma_u(Y)}
  \Biggr] \, \Delta \, du.
\end{align*}
For each choice of $u \in [0,1]$, the matrix-valued function $y
\mapsto \Gamma_u(y)$ is symmetric---that is, $\Gamma_u(y) =
\Gamma_u(-y)$.  Since the distribution of $Y$ is symmetric around
zero, we conclude that $\Exs[\Gamma_u(Y)] = \Exs [\Gamma_u(\Ytil)]$,
where $\Ytil \sim \NORMAL(\thetastar, \sigma^2 I)$, and hence that
\begin{align}
\label{eqn::main_mog}
\norm{\Exs \Big[ \big(w_{\theta}(Y) - w_{\fp}(Y) \big) \, Y \Big]} &
\leq 2 \sup_{u \in [0,1]} \opnorm{\Exs[\Gamma_u(\Ytil)} \, \|\Delta\|_2.
\end{align}
The remainder of the proof is devoted to bounding
$\opnorm{\Exs[\Gamma_u(\Ytil)}$ uniformly over $u \in [0,1]$.  For an
  arbitrary fixed $u \in [0,1]$ let $R$ be an orthonormal matrix such
  that $R \theta_u = \norm{\theta_u} e_1$, where $e_1 \in
  \real^\usedim$ denotes the first canonical basis vector.  Define the
  rotated random vector $\Ycar = R \Ytil$, and note that $\Ycar \sim
  \NORMAL(R \thetastar, \sigma^2 I)$.  Using this transformation, the
  operator norm of the matrix $\Exs[\Gamma_u(\Ytil)]$ is equal to that
  of
\begin{align*} 
\Dmat & = \Exs \Big [ \frac{ \Ycar \Ycar^T }{ \sigma^2 \big(\exp \big(
    \frac{ \inprod{\Ycar}{ \norm{ \theta_u}e_1}}{\sigma^2} \big) +
    \exp \big(- \frac{ \inprod{\Ycar}{ \norm{ \theta_u}
        e_1}}{\sigma^2} \big) \big)^2 } \Big].
\end{align*}
By construction, the matrix $\Dmat$ is diagonal, so that it suffices
to bound the diagonal terms.  Beginning with the first diagonal entry,
we have
\begin{align*} 
\Dmat_{11} = \Exs \Big[ \frac{V_1^2 }{ \sigma^2 \big(\exp \big(
    \frac{ \norm{ \theta_u} \Ycar_1}{\sigma^2} \big) + \exp \big(-
    \frac{ \norm{ \theta_u} \Ycar_1}{\sigma^2} \big) \big)^2 } \Big] &
\leq \Exs \Big[ \frac{ \Ycar_1^2/\sigma^2 } { \exp \big( \frac{ 2
      \norm{ \theta_u} \Ycar_1}{\sigma^2} \big)} \Big]. 
\end{align*}
Defining the event \mbox{$\Event = \{\Ycar_1 \leq
  \frac{\|\thetastar\|_2}{4} \}$,} we condition on it and its
complement to obtain
\begin{align*}
\Dmat_{11} & \leq \Exs \Big[ \frac{ \Ycar_1^2/\sigma^2 } { \exp \big(
    \frac{ 2 \norm{ \theta_u} \Ycar_1}{\sigma^2} \big)} \mid \Event
  \Big] \mprob[\Event] + \Exs \Big[ \frac{ \Ycar_1^2/\sigma^2 } { \exp
    \big( \frac{ 2 \norm{ \theta_u} \Ycar_1}{\sigma^2} \big)} \mid
  \Event^c \Big].
\end{align*}
Conditioned on $\Event$ and $\Event^c$, respectively, we then apply
the bounds~\eqref{Eqn:x_over_exp} and~\eqref{Eqn:x_over_exp2} to
obtain
\begin{align*}
\Dmat_{11} & \leq \frac{\sigma^2}{ e^2 \norm{\theta_u}^2}
\Prob[\Event] + \frac{ \norm{\thetastar}^2} { 16 \sigma^2\exp \big(
  \frac{ \norm{ \theta_u} \norm{\thetastar}}{2\sigma^2} \big)},
\end{align*}
provided $\norm{\thetastar}\norm{\theta_u} \geq 4 \sigma^2$. Noting
that
\begin{align}
\label{EqnHanaGiggle}
 \|\theta_u\|_2 = \|\thetastar + u(\theta - \thetastar)\|_2 & \geq
 \|\thetastar\|_2 - \frac{1}{4} \|\thetastar\|_2 \; = \frac{3}{4}
 \|\thetastar\|_2,
\end{align}
we obtain the bound $\Dmat_{11} \leq \frac{16\sigma^2}{ 9e^2
  \norm{\thetastar}^2} \mathbb{P}(\Event) + \frac{ \norm{\thetastar}^2
  \exp \big( - \frac{ 3 \norm{ \thetastar}^2}{8\sigma^2} \big)} { 16
  \sigma^2 }$, whenever $\norm{\thetastar}^2 \geq 16 \sigma^2/3$.

Note that the mean of $V_1$ is lower bounded as
\begin{align*}
 \Exs[V_1] = \inprod{R \thetastar}{e_1} = \inprod{R \theta_u}{e_1} +
 \inprod{R (\thetastar - \theta_u)}{e_1} \geq \norm{\theta_u} -
 \norm{\thetastar - \theta_u} \stackrel{\mathrm{(i)}}{\geq}
 \frac{\norm{\thetastar}}{2}, 
\end{align*}
where step (i) follows from the lower bound~\eqref{EqnHanaGiggle}.
Consequently, by standard Gaussian tail bounds, we have
\begin{align}
\label{EqnPP}
\Prob[\Event] & \leq \exp \Big(\frac{ - \norm{\thetastar}^2 }{ 32
  \sigma^2} \Big).
\end{align}
Combining the pieces yields
\begin{align*}
\Dmat_{11} & \leq \frac{16\sigma^2}{ 9e^2 \norm{\thetastar}^2} \,
\MYEXP{ - \frac{ \norm{\thetastar}^2 }{ 32 \sigma^2}} + \frac{
  \norm{\thetastar}^2 } { 16 \sigma^2 } \MYEXP{- \frac{ 3 \norm{
      \thetastar}^2}{8\sigma^2}} \qquad \mbox{whenever
  $\norm{\thetastar}^2 \geq 16 \sigma^2/3$.}
\end{align*}

On the other hand, for any index $j \neq 1$, we have
\begin{align*}
\Dmat_{jj} &= \Exs \Big[\frac{ 1 }{ \Big(\exp \Big( \frac{ \norm{
        \theta_u} Y_1}{\sigma^2} \Big) + \exp \Big(- \frac{ \norm{
        \theta_u} Y_1}{\sigma^2} \Big) \Big)^2 } \Big] \; = \; \Exs
\Big[ g \Big(\frac{\norm{\theta_u} \Ycar_1}{\sigma^2} \Big) \Big],
\end{align*}
where the reader should recall the function $g$ from
equation~\eqref{Eqnone_over_expA}.  Once again, conditioning on the
event \mbox{$\Event = \{\Ycar_1 \leq \frac{\|\thetastar\|_2}{4} \}$}
and its complement yields
  \begin{align*}
\Dmat_{jj} & \leq \Exs \Big[ g \Big(\frac{\norm{\theta_u}
    \Ycar_1}{\sigma^2} \Big) \mid \Event \Big] \mprob[\Event] + \Exs
\Big[ g \Big(\frac{\norm{\theta_u} \Ycar_1}{\sigma^2} \Big) \mid
  \Event^c \Big] \\
& \stackrel{(i)}{\leq} \frac{1}{4} \mprob[\Event] + \exp \Big(
-\frac{\|\thetastar\|_2 \|\theta_u\|_2}{4 \sigma^2} \Big) \\
& \stackrel{(ii)}{\leq} \frac{1}{4} \mprob[\Event] + \exp \Big(
-\frac{3 \|\thetastar\|^2_2 }{16\sigma^2} \Big),
\end{align*}
where step (i) follows by applying bound~\eqref{Eqnone_over_expA} to
the first term, and the bound~\eqref{Eqnone_over_expB} with \mbox{$\mu
  = \frac{\|\thetastar\|_2 \|\theta_u\|_2}{4 \sigma^2}$} to the second
term; and step (ii) follows from the bound~\eqref{EqnHanaGiggle}.
Applying the bound~\eqref{EqnPP} on $\Prob[\Event]$ yields
\begin{align*}
\Dmat_{jj} & \leq \frac{1}{4} \exp \big(-\frac{\|\thetastar\|_2^2}{32
  \sigma^2} \big) + \exp \big( -\frac{3 \|\thetastar\|^2_2
}{16\sigma^2} \big) \leq \; 2 \exp \big(-\frac{\|\thetastar\|_2^2}{32
  \sigma^2} \big).
\end{align*}

\noindent
Returning to equation~\eqref{eqn::main_mog}, we have shown that
\begin{align*}
\norm{2 \Exs \Big[ \big(w_{\theta}(Y) - w_{\fp}(Y) \big) \, Y \Big]}
\leq \sivacon_1 \big(1 + \frac{1}{\eta^2} + \eta^2\big)
\MYEXP{-\sivacon_2 \eta^2} \norm{\theta - \fp},
\end{align*}
whenever $\frac{\|\thetastar\|_2^2}{\sigma^2} \geq \eta^2 \geq 16/3$.
On this basis, the bound~\eqref{EqnMogMainTech} holds as long as the
signal-to-noise ratio is sufficiently large,


\subsection{Proof of Corollary~\ref{CorMOGUniform}}
\label{AppMOGUniform}

In order to prove this corollary, it suffices to bound the function
$\rateemunif(\numobs,\delta)$, as previously
defined~\eqref{EqnRateEMUniform}.  Defining the set $\SpecSet \defn
\big \{ \theta \in \real^\usedim \, \mid \norm{\theta - \thetastar}
\leq \norm{\thetastar}/4 \big\}$, our goal is to control the random
variable $Z \defn \sup_{\theta \in \SpecSet} \norm{ \Mem(\theta) -
  \MemSamp(\theta)}$.  For each unit-norm vector $u \in
\real^\usedim$, define the random variable
\begin{align*}
Z_{u} & \defn \sup_{\theta \in \SpecSet} \big \{ \frac{1}{\numobs}
\sum_{i=1}^n (2w_\theta(y_i) - 1) \inprod{y_i}{u} - \Exs
(2w_\theta(Y)-1) \inprod{Y}{u} \big\}.
\end{align*}
Noting that $Z = \sup_{u \in \Sphere{\usedim}} Z_u$, we begin by
reducing our problem to a finite maximum over the sphere
$\Sphere{\usedim}$.  Let $\{u^1, \ldots, u^M \}$ denote a
$1/2$-covering of the sphere $\Sphere{\usedim} = \{v \in \real^\usedim
\, \mid \, \|v\|_2 = 1\}$.  For any $v \in \Sphere{\usedim}$, there is
some index $j \in [M]$ such that $\|v - u^j\|_2 \leq 1/2$, and hence
we can write
\begin{align*}
Z_v & \leq Z_{u^j} + |Z_v - Z_{u^j}| \; \leq \; \max_{j \in [M]}
Z_{u^j} + Z \, \|v - u^j\|_2,
\end{align*}
where the final step uses the fact that $|Z_u - Z_v| \leq Z \, \|u -
v\|_2$ for any pair $(u,v)$. Putting together the pieces, we conclude
that
\begin{align}
\label{EqnFirstDiscrete}
Z & = \sup_{v \in \Sphere{\usedim}} Z_v \; \leq 2 \max_{j \in [M]}
  Z_{u^j}.
\end{align}
Consequently, it suffices to bound the random variable $Z_u$ for a
fixed $u \in \Sphere{\usedim}$.  Letting $\{\rade{i}\}_{i=1}^\numobs$
denote an i.i.d.  sequence of Rademacher variables, for any $\lambda >
0$, we have
\begin{align*}
\Exs \big[ \MYEXP{\lambda Z_u} \big] & \leq \Exs \Big [ \exp \Big(
  \frac{2}{\numobs} \sup_{ \theta \in \SpecSet} \sum_{i=1}^\numobs
  \rade{i} (2 w_\theta(y_i) - 1) \inprod{y_i}{u} \Big) \Big],
\end{align*}
using a standard symmetrization result for empirical processes
(e.g.,~\cite{LedTal91,Koltchinskii}).  Now observe that for any
triplet of $\usedim$-vectors $y$, $\theta$ and $\theta'$, we have the
Lipschitz property
\begin{align*}
\big| 2 w_\theta(y) - 2 w_{\theta'}(y) \big| & \leq \big|
\inprod{\theta}{y} - \inprod{\theta'}{y} \big|.
\end{align*}
Consequently, by the Ledoux-Talagrand contraction for Rademacher
processes~\cite{LedTal91,Koltchinskii}, we have
\begin{align*}
\Exs \Big [ \exp \Big( \frac{2}{\numobs} \sup_{ \theta \in \SpecSet}
  \sum_{i=1}^\numobs \rade{i} (2 w_\theta(y_i) - 1) \inprod{y_i}{u}
  \Big) \Big] & \leq \Exs \Big [ \exp \Big( \frac{4}{\numobs} \sup_{
    \theta \in \SpecSet} \sum_{i=1}^\numobs \rade{i}
  \inprod{\theta}{y_i} \inprod{y_i}{u} \Big) \Big]
\end{align*}
Since any $\theta \in \SpecSet$ satisfies $\|\theta\|_2 \leq 
\frac{5}{4} \|\thetastar\|_2$,  we have
\begin{align*}
\sup_{ \theta \in \SpecSet} \frac{1}{\numobs} \sum_{i=1}^\numobs
\rade{i} \inprod{\theta}{y_i} \inprod{y_i}{u}  & \leq \frac{5}{4}
\|\thetastar\|_2 \opnorm{\frac{1}{\numobs} \sum_{i=1}^\numobs \rade{i}
  y_i y_i^T},
\end{align*}
where $\opnorm{\cdot}$ denotes the $\ell_2$-operator norm of a matrix
(maximum singular value).  Repeating the same discretization argument
over $\{u^1, \ldots, u^M\}$, we find that
\begin{align*}
\opnorm{\frac{1}{\numobs} \sum_{i=1}^\numobs \rade{i} y_i y_i^T} &
\leq 2 \max_{j \in [M]} \frac{1}{\numobs} \sum_{i=1}^\numobs \rade{i}
\inprod{y_i}{u^j}^2.
\end{align*}
Putting together the pieces, we conclude that
\begin{align}
\label{EqnIndia}
\Exs \big[ \MYEXP{\lambda Z_u} \big] & \leq \Exs \Big[ \exp \Big(10
  \lambda \|\thetastar\|_2 \max_{j \in [M]} \frac{1}{\numobs}
  \sum_{i=1}^\numobs \rade{i} \inprod{y_i}{u^j}^2 \Big) \Big] \; \leq
\; \sum_{j=1}^M \Exs \Big[ \exp \Big(10 \lambda \|\thetastar\|_2
  \frac{1}{\numobs} \sum_{i=1}^\numobs \rade{i} \inprod{y_i}{u^j}^2
  \Big) \Big].
\end{align}
Now by assumption, the random vectors $\{y_i\}_{i=1}^\numobs$ are
generated i.i.d. according to the model $y = \eta \thetastar + w$,
where $\eta$ is a Rademacher sign variable, and $v \sim \NORMAL(0, \sigma^2
I$).  Consequently, for any $u \in \real^\usedim$, we have
\begin{align*}
\Exs[e^{\inprod{u}{y}}] & = \Exs[e^{\eta \inprod{u}{\thetastar}}] \;
\Exs[e^{\inprod{u}{v}}] \; \leq \; \MYEXP{ \frac{\|\thetastar\|_2^2 +
    \sigma^2}{2}},
\end{align*}
showing that the vectors $\inprod{y_i}{u}$ are sub-Gaussian with
parameter at most $\gamma = \sqrt{\|\thetastar\|_2^2 + \sigma^2}$.
Therefore, the vectors $\rade{i} \inprod{y_i}{u}^2$ are zero mean
sub-exponential, and have moment generating function bounded as
$\Exs[e^{t (\inprod{y_i}{u})^2}] \leq e^{ \frac{\gamma^2 t^2}{2}}$ for
all $t > 0$ sufficiently small.  Combined with our earlier
inequality~\eqref{EqnIndia}, we conclude that
\begin{align*}
\Exs \big[ \MYEXP{\lambda Z_u} \big] & \leq M \, e^{c \,
  \frac{\lambda^2 \|\thetastar\|_2^2 \gamma^2}{\numobs}} \; \leq \;
e^{c \, \frac{\lambda^2 \|\thetastar\|_2^2 \gamma^2}{\numobs} + 2
  \usedim}
\end{align*}
for all $\lambda$ sufficiently small.  Combined with our first
discretization~\eqref{EqnFirstDiscrete}, we have thus shown that
\begin{align*}
\Exs[e^{\frac{\lambda}{2} Z}] & \leq M e^{c \, \frac{\lambda^2
    \|\thetastar\|_2^2 \gamma^2}{\numobs} + 2 \usedim} \; \leq \; e^{c
  \, \frac{\lambda^2 \|\thetastar\|_2^2 \gamma^2}{\numobs} + 4
  \usedim}.
\end{align*}
Combined with the Chernoff approach, this bound on the MGF implies
that, as long as $\numobs \geq c_1 \usedim \log(1/\delta)$ for a
sufficiently large constant $c_1$, we have
\begin{align*}
Z & \leq c_2 \sigma \norm{\thetastar} \gamma \, \sqrt{ \frac{\usedim
    \log(1/\delta)}{\numobs}}
\end{align*}
with probability at least $1 - \delta$.


\subsection{Guarantees for EM with sample-splitting}
\label{AppCorMOGSam}

In this section, we state and prove a result for the EM algorithm with
sample-splitting for the mixture of Gaussians.

\bcors[Sample-splitting EM guarantees for Gaussian mixtures]
\label{CorMOGSam}
Consider a Gaussian mixture model satisfying the $\MYSNR(\MOGSNR)$
condition~\eqref{EqnMOGSNR}, and any initialization $\theta^0$ such
that $\|\theta^0 - \thetastar\|_2 \leq \frac{\|\thetastar\|_2}{4}$.
Given a sample size $\numobs \geq 16T \log (6T/\delta)$, then with
probability at least $1 - \delta$, the sample-splitting EM iterates
$\{\theta^t\}_{t=0}^\Tfinal$ satisfy the bound
\begin{align}
\norm{\theta^t - \thetastar} \leq \kappa^t \norm{\theta^0 -
  \thetastar} + \frac{c}{1 - \kappa} \Big( \sigma \sqrt{\frac{\usedim
    \Tfinal \log(\Tfinal/\delta)}{\numobs} } + \sqrt{ \frac{\Tfinal
    \log (\Tfinal/\delta)}{\numobs}}\norm{\gopt} \Big).
\end{align}
\ecors
\noindent It is worth comparing the result here to the result
established earlier in Corollary~\ref{CorMOGUniform}. The
sample-splitting EM algorithm is more sensitive to the number of
iterations which determines the batch size and needs to be chosen in
advance. Supposing that the number of iterations were chosen optimally
however the result has better dependence on $\norm{\theta^*}$ and
$\sigma$ at the cost of a logarithmic factor in $n$.

\begin{proof}
The proof follows by establishing a bound on the function
$\rateem(n,\delta)$.  Define $\mathcal{S} = \{\theta: \norm{\theta-
  \thetastar} \leq \frac{\norm{\thetastar}}{4} \}$. Recalling the
updates in~\eqref{EqnMemSampMOG} and~\eqref{EqnMemMOG}, note that
\begin{align*}
\norm{\mem{\theta} - \memsamp{\theta}} \leq \underbrace{\norm{
    \fsavgtwo{1}{} Y_i}}_{T_1} + \underbrace{\norm{ \fsavgtwo{1}{}
    w_{\theta} (Y_i) Y_i- \Exs w_{\theta} (Y) Y}}_{T_2}.
\end{align*}
We bound each of these terms in turn, in particular showing that
\begin{align}
\label{EqnHarding}
\max \{ T_1, T_2 \} & \leq \sqrt{ \frac{ \log
    (8/\delta)}{2n}}\norm{\gopt} + c \sigma \sqrt{\frac{d
    \log(1/\delta)}{n}},
\end{align}
with probability at least $1-\delta$.

\paragraph{Control of $T_1$:}
Observe that since $Y \sim (2Z - 1) \thetastar + \regnoise$ we have
\begin{align*}
T_1 = \norm{ \fsavgtwo{1}{} Y_i} \leq \norm{ \fsavgtwo{1}{}
  \regnoise_i} + \Big|\fsavgtwo{1}{} (2 Z_i - 1) \Big|
\norm{\thetastar}.
\end{align*}
Since $Z_i$ are i.i.d Bernoulli variables, Hoeffding's inequality
implies that
\begin{align*}
\Big|\fsavgtwo{1}{} (2 Z_i - 1) \Big| \leq \sqrt{ \frac{ \log
    (8/\delta)}{2n}}.
\end{align*}
with probability at least $1 - \frac{\delta}{4}$.  On the other hand,
the vector $U_1 \defn \fsavgtwo{1}{} {\regnoise_i}$ is zero-mean and
sub-Gaussian with parameter $\sigma/\sqrt{n}$, whence the squared norm
$\|U_1\|_2^2$ is sub-exponential.  Using standard bounds for
sub-exponential variates and the condition $\numobs > \sigma \usedim$,
we obtain
\begin{align*}
 \norm{U_1} \leq c_2 \sigma \sqrt{\frac{d \log(1/\delta)}{n} }.
\end{align*}
with probability at least $1-\delta/4$.  Combining the pieces yields
the claimed bound~\eqref{EqnHarding} on $T_1$.

\paragraph{Control of $T_2$:}
By triangle inequality, we have
\begin{align*}
T_2 \leq & \big| \fsavgtwo{1}{} w_{\theta} (Y_i) (2 Z_i - 1)- \Exs
w_{\theta} (Y) (2Z - 1)\big| \norm{\gopt} + \norm{
  \fsavgtwo{1}{}w_{\theta}(Y_i) \regnoise_i - \Exs w_{\gvec} (Y)
  \regnoise}.
\end{align*}
The random variable $w_{\theta}(Y) (2Z - 1)$ lies in the interval $[-1,1]$,
so that Hoeffding's inequality implies that
\begin{align*}
 \big| \fsavgtwo{1}{} w_{\theta} (Y_i) (2 Z_i - 1)- \Exs w_{\theta}
 (Y) (2Z - 1)\big| \norm{\gopt} \leq \sqrt{ \frac{ \log
     (6/\delta)}{2\numobs}} \norm{\gopt}.
\end{align*}
with probability at least $1 - \delta/4$. 

Next observe that the random vector $U_2 \defn
\fsavgtwo{1}{}w_{\theta^t}(X_i) \regnoise_i - \Exs w_{\gvec^t}(X)
\regnoise$ is zero mean and sub-Gaussian with parameter
$\sigma/\sqrt{n}$. Consequently, as in our analysis of $T_1$,we
conclude that
\begin{align*}
\norm{U_2} & \leq c \sigma \sqrt{\frac{d \log(1/\delta)}{n} }.
\end{align*}
with probability at least $1-\delta/4$.  Putting together the pieces
yields the claimed bound~\eqref{EqnHarding} on $T_2$, thereby
completing the proof of the corollary.
\end{proof}


\section{Proofs for mixtures of regressions}

In this appendix, we provide proofs of results related to the mixture
of regressions model, as presented in Section~\ref{SecMORAnalysis}.
More specifically, we first prove Corollary~\ref{CorMORPop} on the
population level behavior, followed by the proof of
Corollaries~\ref{CorMORSamp} and~\ref{CorMORSGD} on the behavior of
sample-splitting EM updates and stochastic gradient EM updates,
respectively.

\subsection{Proof of Corollary~\ref{CorMORPop}}
\label{AppCorMORPop}

We begin by proving part (a) of the corollary on the population EM
update, which is based on maximizing the function
\begin{align*}
\QFUN{\theta'}{\theta} \defn - \frac{1}{2} \Exs \big[ w_{\theta}(X,Y)
  (Y - \inprod{X}{\theta'})^2 + (1 - w_{\theta}(X,Y)) (Y +
  \inprod{X}{\theta'})^2 \big], 
\end{align*}
where $w_{\theta}(x,y) \defn \frac{ \exp \big( \frac{- (y -
    \inprod{x}{\rvec})^2 }{2\sigma^2} \big)}{ \exp \big( \frac{-(y -
    \inprod{x}{\rvec})^2}{2\sigma^2} \big) + \exp \big( \frac{-(y +
    \inprod{x}{\rvec})^2 }{2\sigma^2} \big)}$. Observe that function
$\qfunone{\rvec^*}$ is $\strongparam$-strongly concave, with
$\strongparam$ equal to the smallest eigenvalue of the matrix $\Exs
[XX^T]$.  Since $\Exs[ X X^T] = I$ by assumption, we see that strong
concavity holds with $\strongparam = 1$.

It remains to verify condition~\gstwonospace.  Define the difference
function $\Delta_w(X,Y) \defn w_{\theta}(X,Y) - w_{\thetastar}(X,Y)$,
and the difference vectors $\Delta = \theta - \thetastar$.  Using this notation, for
this model, we need to show that
\begin{align*}
\norm{2 \Exs \big[\Delta_w(X,Y) Y
    X\big]} & < \norm{\Delta}.
\end{align*}
Fix any $\DelPrime \in \mathbb{R}^d$. It suffices for us to show
that,
\begin{align*}
\inprod{2 \Exs \big[\Delta_w(X,Y) Y
    X\big]}{\DelPrime} & < \norm{\Delta}\norm{\DelPrime}.
\end{align*}
Note that we can write $Y \stackrel{d}{=} (2Z - 1)
\inprod{X}{\thetastar} + \regnoise$, where $Z \sim \mbox{Ber}(1/2)$ is
a Bernoulli variable.  Using this notation, it is equivalent to show
\begin{align}
\label{EqnAA}
\Exs \big [ \Delta_w(X,Y) (2Z - 1) \inprod{X}{\thetastar}
  \inprod{X}{\DelPrime} \big] + \Exs \big[\Delta_w(X,Y) \regnoise
  \inprod{X}{\DelPrime} \big] \leq \gamma \norm{\Delta} \;
\norm{\DelPrime}
\end{align}
for $\gamma \in [0,1/2)$ in order to establish contractivity.  In
  order to prove the theorem with the desired upper bound on $\kappa$
  we need to show~\eqref{EqnAA} with $\gamma \in [0,1/4).$ The
    following lemma provides control on the two terms:
\blems
\label{LemMORTech}
Under the conditions of Corollary~\ref{CorMORPop}, there is a constant
$\gamma < 1/4$ such that for any fixed vector $\DelPrime$ we have
\begin{subequations}
\begin{align} 
\label{EqnMORTechA}
\big| \Exs \big[ \Delta_w(X,Y) (2Z-1) \inprod{X}{\thetastar}
  \inprod{X}{\DelPrime} \big] \big| & \leq \frac{\gamma}{2}
\norm{\Delta} \norm{\DelPrime}, \quad \mbox{and} \\
\label{EqnMORTechB}
\big| \Exs \big[ \Delta_w(X,Y) \regnoise \inprod{X}{\DelPrime} \big]
\big| & \leq \frac{\gamma}{2}\, \norm{\Delta} \: \norm{\DelPrime}.
\end{align}
\end{subequations}
\elems 
In conjunction, these bounds imply that $\inprod{\Exs
  \big[\Delta_w(X,Y) Y X\big]}{\DelPrime} \leq \gamma \norm{\Delta} \,
\norm{\DelPrime}$ with $\gamma \in [0,1/4)$, as claimed. 

Part (b) of the corollary is nearly immediate given part (a). Our first task
is to verify smoothness of the objective
$\qfunone{\rvec^*}.$  
The smoothness parameter is given by the largest eigenvalue
of the Hessian of $\qfunone{\rvec^*}$ which is $\Exs[XX^T]$. 
Since $\Exs[ X X^T] = I$ by assumption, we see that 
smoothness holds with $\smoothparam = 1$. Finally, we need to verify 
the condition GS with the desired contraction coefficient. Some algebra shows
that it suffices to show that under the stated assumptions 
of the corollary we have
\begin{align*}
2 \norm{\Exs \big[\Delta_w(X,Y) Y
    X\big]} \leq \kappa \norm{\Delta},
\end{align*}
for $\kappa \in [0, \frac{1}{2}).$ This is an immediate consequence
of Lemma~\ref{LemMORTech}.

It remains to prove Lemma~\ref{LemMORTech}.  Since the standard
deviation $\sigma$ is known, a simple rescaling argument allows us to
take $\sigma = 1$, and replace the weight function
in~\eqref{EqnMORweight} with
\begin{align}
\label{EqnMORweightmod}
w_{\theta}(x,y) = \frac{ \exp \big( \frac{- (y -
    \inprod{x}{\theta})^2 }{2} \big)} { \exp \big( \frac{-(y
    - \inprod{x}{\theta})^2}{2} \big) + \exp \big( \frac{-(y
    + \inprod{x}{\theta})^2 }{2} \big)}.
\end{align}
\noindent Our proof makes use of the following elementary result on
Gaussian random vectors: \\
\blems
\label{LemSqNorm}
Given a Gaussian random vector $X \sim \NORMAL(0,I)$ and any fixed vectors
$u, v \in \real^\usedim$, we have
\begin{subequations}
\begin{align}
\label{EqnSqNormA}
\Exs [ \inprod{X}{u}^2 \inprod{X}{v}^2 ] & \leq 3 \norm{u}^2
\norm{v}^2 \quad \mbox{with equality when $u = v$, and} \\
\label{EqnSqNormB}
\Exs [ \inprod{X}{u}^4 \inprod{X}{v}^2 ] & \leq 15 \norm{u}^4
\norm{v}^2.
\end{align}
\end{subequations}
\elems

\begin{proof}
For any fixed orthonormal matrix $R \in \real^{d \times d}$, the
transformed variable $R^T X$ also has a $\NORMAL(0, I)$ distribution,
and hence $\Exs [ \inprod{X}{u}^2 \inprod{X}{v}^2 ] = \Exs [
  \inprod{X}{Ru}^2 \inprod{X}{Rv}^2 ]$.  Let us choose $R$ such that
$Ru = \norm{u} e_1$.  Introducing the shorthand $z = Rv$, we have
\begin{align*}
\Exs [ \inprod{X}{Ru}^2 \inprod{X}{Rv}^2 ] = \Exs [ \norm{u}^2 X_1^2
  \sum_{i=1}^\usedim \sum_{j=1}^\usedim X_i X_j z_i z_j ] & =
\norm{u}^2 (3 z_1^2 + (\norm{z}^2 - z_1^2)) \\
& \leq 3 \norm{u}^2 \norm{z}^2 \; = 3 \norm{u}^2 \norm{v}^2.
\end{align*}
A similar argument yields the second claim.
\end{proof}


With these preliminaries in place, we can now begin the proof of
Lemma~\ref{LemMORTech}.  Recall that $\Delta =
\theta - \thetastar$ and that $\DelPrime$ is any fixed vector in 
$\mathbb{R}^d$. Define $\theta_u = \thetastar + u \Delta$ for a scalar $u \in
[0,1]$. Recall that by our assumptions guarantee that
\begin{subequations}
\begin{align}
\label{EqnMORass}
\norm{\Delta} \leq \frac{ \norm{\thetastar}} {\deltaconstmor},
\qquad \mbox{and} \qquad \norm{\thetastar} \geq \snrconstmor.
\end{align}
For future reference, we observe that
\begin{align}
\label{EqnNormLB}
\norm{\theta_u} \geq \norm{\thetastar} - \norm{\Delta} \geq
\frac{\norm{\thetastar}}{2}.
\end{align}
\end{subequations}
Noting that Lemma~\ref{LemMORTech} consists of two separate
inequalities~\eqref{EqnMORTechA} and~\eqref{EqnMORTechB}, we treat
these cases separately.

\subsubsection{Proof of inequality~\eqref{EqnMORTechA}}
\label{SecMORTechA}

\newcommand{\morconstone}{14}
\newcommand{\morconsttwo}{32}

We split the proof of this bound into two separate cases: namely,
$\norm{\Delta} \leq 1$ and $\norm{\Delta} > 1$.
 
\paragraph{Case $\norm{\Delta} \leq 1$:} We then have.
\begin{align*}
\frac{d}{d \theta} w_{\theta}(X,Y) & = \frac{2 Y X}{ \big( \exp \big(
  Y \inprod{X}{\theta} \big) + \exp \big( -Y \inprod{X}{\theta}
  \big)\big)^2}.
\end{align*}
Thus, using a Taylor series with integral form remainder on the
function $\theta \mapsto w_\theta(X, Y)$ yields
\begin{align}
\label{eqn::mor_taylor}
\Delta_w(X, Y) = \int_0^1 \frac{2 Y \inprod{X}{\Delta}}{(\exp(Z_u) +
  \exp(-Z_u))^2} du,
\end{align}
where $Z_u \defn Y \inprod{X}{\thetastar + u \Delta}$.  Substituting
for $\Delta_w(X,Y)$ in inequality~\eqref{EqnMORTechA}, we see that it
suffices to show
\begin{align}
\label{EqnPortugalGhana}
\int_0^1 \underbrace{\Exs \big[ \frac{2 Y \inprod{X}{\thetastar}
    }{(\exp(Z_u) + \exp(-Z_u))^2} (2Z - 1) \inprod{X}{\Delta}
    \inprod{X}{\DelPrime}\big]}_{A_u} du & \leq \frac{\gamma}{2} \norm{\Delta}
\norm{ \DelPrime}.
\end{align}
for some $\gamma \in [0, 1/4)$.  The following auxiliary result is
  central to establishing this claim:
\blems
\label{LemMORTechA}
There is a $\gamma \in [0, 1/4)$ such that for each $u \in [0,1]$, we
  have
\begin{subequations}
\begin{align}
\label{EqnMORTechAone}
\sqrt{\Exs \big[ \frac{Y^2 \inprod{X}{\theta_u}^2 }{(\exp(Z_u) +
      \exp(-Z_u))^4} \big]} & \leq \frac{\gamma}{\morconstone}, \quad \mbox{and}
\\
\label{EqnMORTechAtwo}
\sqrt{\Exs \big[ \frac{Y^2 }{(\exp(Z_u) + \exp(-Z_u))^4} \big]} & \leq
\frac{\gamma}{\morconsttwo} \qquad \mbox{whenever $\norm{\Delta} \leq 1$.}
\end{align}
\end{subequations}
\elems
\noindent See Section~\ref{AppLemMORTechA} for the proof of this
lemma. \\

Using Lemma~\ref{LemMORTechA}, let us bound the quantity $A_u$ from
equation~\eqref{EqnPortugalGhana}.  Since $\thetastar = \theta_u -u
\Delta$, we have $A_u = B_1 + B_2$, where
\begin{align*}
B_1 & \defn \Exs \big[ \frac{2 Y \inprod{X}{\theta_u} }{(\exp(Z_u) +
    \exp(-Z_u))^2} (2Z - 1) \inprod{X}{\Delta}
  \inprod{X}{\DelPrime}\big], \quad \mbox{and} \\
B_2 & \defn - \Exs \big[ \frac{2 Y u \inprod{X}{\Delta} }{(\exp(Z_u) +
    \exp(-Z_u))^2} (2Z - 1) \inprod{X}{\Delta}
  \inprod{X}{\DelPrime}\big].
\end{align*}
In order to show that $A_u \leq \frac{\gamma}{2} \norm{\Delta} \,
\norm{\DelPrime}$, it suffices to show that $\max \{B_1, B_2 \} \leq
\frac{\gamma}{4} \norm{\Delta} \norm{\DelPrime}$. \\

\paragraph{Bounding $B_1$:} By the Cauchy-Schwarz inequality, we have
\begin{align*}
B_1 & \leq \sqrt{ \Exs \big[ \frac{y^2 \inprod{X}{\theta_u}^2
    }{(\exp(Z_u) + \exp(-Z_u))^4} \big]} \sqrt{ \Exs \big[4(2Z - 1)^2
    \inprod{X}{\Delta}^2 \inprod{X}{\DelPrime}^2 \big] } \\
& \leq \frac{\gamma}{\morconstone} \sqrt{ \Exs \big[4
    \inprod{X}{\Delta}^2 \inprod{X}{\DelPrime}^2 \big] },
\end{align*}
where the second step follows from the bound~\eqref{EqnMORTechAone},
and the fact that $(2Z - 1)^2 = 1$.  Next we observe that $\Exs \big[4
  \inprod{X}{\Delta}^2 \inprod{X}{\DelPrime}^2 \big] \leq 12
\norm{\Delta}^2 \norm{\DelPrime}^2$, where we have used the
bound~\eqref{EqnSqNormA} from Lemma~\ref{LemSqNorm}.  Combined with
our earlier bound, we conclude that $B_1 \leq \frac{\gamma}{4}
\norm{\Delta} \, \norm{\DelPrime}$, as claimed.

\paragraph{Bounding $B_2$:}
Similarly, another application of the Cauchy-Schwarz inequality yields
\begin{align*}
B_2 & \leq \sqrt{ \Exs \big[ \frac{y^2 }{(\exp(Z_u) + \exp(-Z_u))^4}
    \big]} \sqrt{ \Exs \big[4u^2(2Z - 1)^2 \inprod{X}{\Delta}^4
    \inprod{X}{\DelPrime}^2 \big] } \\
 & \leq \frac{\gamma}{\morconsttwo} \sqrt{ \Exs \big[4u^2
    \inprod{X}{\Delta}^4 \inprod{X}{\DelPrime}^2 \big] },
\end{align*}
where the second step follows from the bound~\eqref{EqnMORTechAtwo},
and the fact that $(2Z-1)^2 = 1$.
In this case, we have
\begin{align*}
\Exs \big[4u^2 \inprod{X}{\Delta}^4
  \inprod{X}{\DelPrime}^2 \big] \stackrel{\mathrm{(i)}}{\leq} 60
\norm{\Delta}^4 \norm{\DelPrime}^2 \stackrel{\mathrm{(ii)}}{\leq}
60 \norm{\Delta}^2 \norm{\DelPrime}^2,
\end{align*}
where \mbox{step (i)} uses the bound~\eqref{EqnSqNormB} from
Lemma~\ref{LemSqNorm}, and \mbox{step (ii)} that $\norm{\Delta} \leq
1$.  Combining the pieces, we conclude that $B_2 \leq \frac{\gamma}{4}
\norm{\Delta} \norm{\DelPrime}$, which completes the proof of
inequality~\eqref{EqnMORTechA} in the case $\norm{\Delta} \leq 1$.


\paragraph{Case $\norm{\Delta} > 1$:}

We now turn to the second case of the bound~\eqref{EqnMORTechA}.  Our
argument (here and in later sections) makes use of various probability
bounds on different events, which we state here for future reference.
These events involve the scalar $\tau \defn C_\tau \sqrt{\log
  \norm{\thetastar}}$ for a constant $C_\tau$, as well as the vectors
\begin{align*}
\Delta \defn \theta - \thetastar,~\mbox{and}~\theta_u \defn \thetastar
+ u \, \Delta \quad \mbox{for some fixed $u \in [0,1]$.}
\end{align*}
\blems[Event bounds]
\label{LemEventBounds}
\hfill
\begin{enumerate}
\item[(i)] For the event $\eone \defn \big \{ \sign
  (\inprod{X}{\thetastar}) = \sign (\inprod{X}{\theta_u}) \big\}$, we
  have
$\mprob[\eone^c]  \leq \frac{\|\Delta\|_2}{\|\thetastar\|_2}$.
\item[(ii)] For the event $\etwo \defn \big \{ |
  \inprod{X}{\thetastar} | > \tau \big \} \cap \big \{ |
  \inprod{X}{\theta_u} | > \tau \big \} \cap \big \{ | \regnoise |
  \leq \frac{\tau}{2} \big \}$, we have
\begin{align*}
\mprob[\etwo^c] & \leq \frac{\tau}{\norm{\thetastar}} +
\frac{\tau}{\norm{\theta_u}} + 2 \exp \Big( - \frac{\tau^2}{2} \Big).
\end{align*}
\item[(iii)] For the event $\ethree \defn \big \{
  |\inprod{X}{\thetastar}| \geq \tau \big \} \bigcup \big \{ |
  \inprod{X}{\theta_u} | \geq \tau \big \}$, we have $\mprob
  \big[\ethree^c \big] \leq \frac{\tau}{\norm{\thetastar}} +
  \frac{\tau}{\norm{\theta_u}}$.
\item[(iv)] For the event $\efour \defn \big \{ |\regnoise| \leq
  \tau/2 \big \}$, we have $\mprob[\efour^c] \leq 2 \MYEXP{-
    \frac{\tau^2}{2}}$.
\item[(v)] For the event $\efive \defn \big \{ |\inprod{X}{\theta_u}|
  > \tau \big \}$, we have $\mprob[\efive^c] \leq
  \frac{\tau}{\norm{\theta_u}}$.
\item[(vi)] For the event $\esix \defn \big \{
  |\inprod{X}{\thetastar}| > \tau \}$, we have $\mprob[\esix^c] \leq
  \frac{\tau}{\norm{\thetastar}}$.
\end{enumerate}
\elems

Various stages of our proof involve controlling the second moment
matrix $\Exs[XX^T]$ when conditioned on some of the events given
above:
\blems[Conditional covariance bounds]
\label{LemCC}
Conditioned on any event \mbox{$\mathcal{E} \in \{\eone \cap \etwo,
  \eone^c, \efive^c, \esix^c\}$}, we have $\opnorm{\Exs[XX^T \, \mid
    \, \mathcal{E}} \leq 2$.
\elems
\noindent See Section~\ref{AppLemCC} for the proof of this result. \\
 
\noindent With this set-up, our goal is to bound the quantity
\begin{align*}
T = \big| \Exs \big[ \Delta_w(X, Y) (2Z-1) \inprod{X}{\thetastar}
  \inprod{X}{\DelPrime}\big] \big| \leq \Exs \big[ \big|\Delta_w(X, Y)
  (2Z-1) \inprod{X}{\thetastar} \inprod{X}{\DelPrime}\big| \big].
\end{align*}
For any measurable event $\Event$, we define $\Psi(\Event) \defn \Exs
\big[\big| \Delta_w(X, Y) (2Z-1) \inprod{X}{\thetastar}
  \inprod{X}{\DelPrime} \big| \, \mid \, \Event \big] \;
\mprob[\Event]$, and note that by successive conditioning, we have
\begin{align}
\label{EqnSuccessive}
T & \leq \Psi(\eone \cap \etwo) + \Psi(\eone^c) + \Psi(\efour^c) +
\Psi(\efive^c) + \Psi(\esix^c).
\end{align}
We bound each of these five terms in turn.

\paragraph{Bounding $\Psi(\eone \cap \etwo)$: }
Applying the Cauchy-Schwarz inequality and using the fact
that $(2Z - 1)^2 = 1$ yields
\begin{align}
\label{EqnLassi}
\Psi(\eone \cap \etwo) & \leq \sqrt{ \Exs \big[ \Delta_w(X,Y)^2
    \inprod{X}{\DelPrime}^2 | \eone \cap \etwo \big] } \sqrt{ \Exs
  \big[ \inprod{X}{\thetastar}^2 | \eone \cap \etwo \big]}.
\end{align}
We now bound $\Delta_w(X,Y)$ conditioned on the event $\eone \cap
\etwo$.  Since $\sign(\inprod{X}{\thetastar}) =
\sign(\inprod{X}{\theta_u})$ on the event $\eone$, we have
\begin{subequations}
\begin{align}
\label{EqnSignMatch}
\sign(Y\inprod{X}{\thetastar}) & = \sign(Y \inprod{X}{\theta_u}).
\end{align}
Conditioned on the event $\etwo$, observe that $|Y| = |(2Z-1)
\inprod{X}{\thetastar} + \regnoise| \geq | \inprod{X}{\thetastar}| -
|\regnoise| \geq \frac{\tau}{2}$, which implies that
\begin{align}
\label{EqnMagnitude}
\min \big \{ | Y \inprod{X}{\thetastar} |, | Y \inprod{X}{\theta} |
\big \} & \geq \frac{\tau^2}{2}.
\end{align}
\end{subequations}
Recalling the weight function~\eqref{EqnMORweightmod}, we claim that
when conditions~\eqref{EqnSignMatch} and~\eqref{EqnMagnitude} hold,
then
\begin{align}
\label{EqnDeltaw}
|\Delta_w(X,Y)| = |w_{\theta_u}(X,Y) - w_{\thetastar}(X,Y)|
\stackrel{\mathrm{(i)}}{\leq} \frac{ \exp (- \tau^2/2) } {\exp (-
  \tau^2/2) + \exp(\tau^2/2)} \leq \exp( -\tau^2).
\end{align}
We need to verify inequality (i): suppose first that
$\sign(Y\inprod{X}{\thetastar}) = 1$.  In this case, both
$w_{\theta_u}(X,Y)$ and $w_{\thetastar}(X,Y)$ are at least $\frac{
  \exp ( \tau^2/2) } {\exp (- \tau^2/2) + \exp(\tau^2/2)}$. Since each
of these terms are upper bounded by $1$, we obtain the claimed bound
on $\Delta_w(X,Y)$. The case when
$\mathrm{sign}(Y\inprod{X}{\thetastar}) = -1$ follows analogously.

Combined with our earlier bound~\eqref{EqnLassi}, we have shown
\begin{align*}
\Psi(\eone \cap \etwo) & \leq \exp(-\tau^2) \sqrt{ \Exs \big[
    \inprod{X}{\DelPrime}^2 | \eone \cap \etwo \big] } \sqrt{ \Exs
  \big[ \inprod{X}{\thetastar}^2 | \eone \cap \etwo \big]}.
\end{align*}
Applying Lemma~\ref{LemCC} with $\Event = \eone \cap \etwo$ yields
$\Psi(\eone \cap \etwo) \leq 2 \norm{\DelPrime} \norm{\thetastar}
\MYEXP{-\tau^2}$.


\paragraph{Bounding $\Psi(\eone^c)$:} 
Combining the Cauchy-Schwarz inequality with
Lemma~\ref{LemEventBounds}(i), we have
\begin{align}
\label{EqnBrewedOne}
\Psi(\eone^c) & \leq \sqrt{ \Exs \big[ \inprod{X}{\DelPrime}^2 |
    \eone^c \big] } \sqrt{ \Exs \big[ \inprod{X}{\thetastar}^2 |
    \eone^c \big]} \frac{\norm{\Delta}}{\norm{\thetastar}}.
\end{align}
We first claim that $\Exs \big[ \inprod{X}{\thetastar}^2 \mid \eone^c
  \big] \leq \Exs \big[ \inprod{X}{\Delta}^2 \mid \eone^c \big]$.  To
establish this bound, it suffices to show that conditioned on $\eone$,
we have $\inprod{X}{\thetastar}^2 \leq \inprod{X}{\Delta}^2$.  Note
that event $\eone$ implies that $\inprod{X}{\thetastar} \,
\inprod{X}{\theta_u} \leq 0$.  Consequently, conditioned on event
$\eone$, we have
\begin{align*}
\inprod{X}{\thetastar}^2 = \frac{1}{4} \inprod{X}{(\thetastar
  -\theta_u) + (\theta_u + \thetastar)}^2 & \leq \frac{1}{2}
\inprod{X}{\thetastar - \theta_u}^2 + \frac{1}{2} \inprod{X}{\theta_u
  + \thetastar}^2 \\
& \stackrel{(i)}{\leq} \inprod{X}{\thetastar - \theta_u}^2 \\
& \stackrel{(ii)}{\leq} \inprod{X}{\Delta}^2
\end{align*}
where step (i) makes use of the bound $\inprod{X}{\thetastar} \,
\inprod{X}{\theta_u} \leq 0$; and step (ii) follows since
\mbox{$\theta_u = \thetastar + u \Delta$,} and $u \in [0,1]$.

Returning to equation~\eqref{EqnBrewedOne}, we have
\begin{align*}
\Psi(\eone^c) & \leq \sqrt{ \Exs \big[ \inprod{X}{\DelPrime}^2 |
    \eone^c \big] } \sqrt{ \Exs \big[ \inprod{X}{\Delta}^2 | \eone^c
    \big]} \frac{\norm{\Delta}}{\norm{\thetastar}}
\stackrel{\mathrm{(i)}}{\leq}
\frac{2\norm{\DelPrime}\norm{\Delta}^2}{\norm{\thetastar}}
\end{align*}
where step (i) follows from the conditional covariance bound of
Lemma~\ref{LemCC}.


\paragraph{Bounding $\Psi(\efour^c)$:}
Combining the Cauchy-Schwarz inequality with
Lemma~\ref{LemEventBounds}(iv) yields
\begin{align*}
\Psi(\efour^c) & \leq 2 \sqrt{ \Exs \big[ \inprod{X}{\DelPrime}^2 |
    \efour^c \big] } \sqrt{ \Exs \big[ \inprod{X}{\thetastar}^2 |
    \efour^c \big]} \; \MYEXP{- \frac{\tau^2}{2}}.
\end{align*} 
Observe that by the independence of $\regnoise$ and $X$, conditioning
on $\efour^c$ has no effect on the second moment of $X$.  Since
$\Exs[XX^T] = I$, we conclude that \mbox{$\Psi(\efour^c) \leq 2
  \norm{\DelPrime} \norm{\thetastar} \, \MYEXP{ - \frac{\tau^2}{2}}$.}

\paragraph{Bounding $\Psi(\efive^c)$: } 

Combining the Cauchy-Schwarz inequality with
Lemma~\ref{LemEventBounds}(v) yields $\Psi(\efive^c) \leq
\frac{\tau}{\norm{\theta_u}} \sqrt{ \Exs \big[ \inprod{X}{\DelPrime}^2
    | \efive^c \big] } \sqrt{ \Exs \big[ \inprod{X}{\thetastar}^2 |
    \efive^c \big]}$.  Conditioned on the event $\efive^c$, we have
\begin{align*}
\inprod{X}{\thetastar}^2 \leq 2 \inprod{X}{\theta_u}^2 + 2
\inprod{X}{\Delta}^2 \leq 2 \tau^2 + 2 \inprod{X}{\Delta}^2.
\end{align*}
Together with Lemma~\ref{LemCC}, we obtain the bound
\begin{align*}
\Psi(\efive^c) & \leq \frac{2 \tau \norm{\DelPrime} \sqrt{ \tau^2 + 2
    \norm{\Delta}^2 }}{\norm{\theta_u}} \stackrel{\mathrm{(i)}}{\leq}
\frac{2 \tau \norm{\DelPrime} \norm{\Delta} \sqrt{ \tau^2 + 2
}}{\norm{\theta_u}},
\end{align*}
where step (i) uses the fact that $\norm{\Delta} \geq 1$.


\paragraph{Bounding $\Psi(\esix^c)$: }
Combining the Cauchy-Schwarz inequality with
Lemma~\ref{LemEventBounds}(vi) yields $\Psi(\esix^c) \leq
\frac{\tau}{\norm{\thetastar}} \sqrt{ \Exs \big[
    \inprod{X}{\DelPrime}^2 | \esix^c \big] } \sqrt{ \Exs \big[
    \inprod{X}{\thetastar}^2 | \esix^c \big]}$.  Conditioned on the
event $\esix^c$, we have \mbox{$\inprod{X}{\thetastar}^2 \leq
  \tau^2$,} and so applying Lemma~\ref{LemCC} with $\Event = \esix^c$
yields $\Psi(\esix^c) \leq \frac{\sqrt{2} \tau^2
  \norm{\DelPrime}}{\norm{\thetastar}}$. \\

\vspace*{.06in}

We have thus obtained bounds on all five terms in the
decomposition~\eqref{EqnSuccessive}.  We combine these bounds with the
with lower bound $\|\theta_u\|_2 \geq \frac{\|\thetastar\|_2}{2}$ from
equation~\eqref{EqnNormLB}, and then perform some algebra to obtain
\begin{align*}
T & \leq c \, \norm{\Delta}\norm{\DelPrime} \: \Big \{ \frac{\tau^2}{
  \norm{\thetastar}} + \norm{\thetastar} \MYEXP{ -\tau^2/2} \Big\} + 
2\norm{\DelPrime}\frac{ \norm{\Delta}^2}{\norm{\thetastar}},
\end{align*}
where $c$ is a universal constant.  In particular, selecting $\tau =
\LASSI_\tau \sqrt{\log \norm{\thetastar}}$ for a sufficient large
constant $\LASSI_\tau$, selecting the constant $\snrconstmor$
in~\eqref{EqnMORass} sufficiently large yields the
claim~\eqref{EqnMORTechA}.


\subsubsection{Proof of inequality~\eqref{EqnMORTechB}}
\label{SecMORTechB}

As in Section~\ref{SecMORTechA}, we treat the cases $\norm{\Delta}\leq
1$ and $\norm{ \Delta} \geq 1$ separately.

\subsubsection{Case $\norm{\Delta} \leq 1$:} 
As before, by a Taylor expansion of the function $\theta \mapsto
\Delta_w(X, Y)$, it suffices to show that
\begin{align*} 
\int_0^1 \Exs \big[ \frac{2 Y \regnoise }{(\exp(Z_u) + \exp(-Z_u))^2}
  \inprod{X}{\Delta} \inprod{X}{\DelPrime}\big] du & \leq
\frac{\gamma}{2} \norm{\Delta} \norm{ \DelPrime}.
\end{align*}
For any fixed $u \in [0,1]$, the Cauchy-Schwarz inequality implies
that
\begin{align*} 
\Exs \big[ \frac{2 Y \regnoise \; \inprod{X}{\Delta}
    \inprod{X}{\DelPrime}}{(\exp(Z_u) + \exp(-Z_u))^2} \big] & \leq
\sqrt{ \Exs \big[ \frac{4 Y^2  }{(\exp(Z_u) + \exp(-Z_u))^4}
    \big]} \; \sqrt{ \Exs \big[  \regnoise^2 \inprod{X}{\Delta}^2
    \inprod{X}{\DelPrime}^2 \big]} \\
& \stackrel{\mathrm{(i)}}{\leq} \sqrt{ \Exs \big[ \frac{4 Y^2
    }{(\exp(Z_u) + \exp(-Z_u))^4} \big]} \; \sqrt{3 \, \norm{\Delta}^2 \:
\norm{\DelPrime}^2} \\
& \stackrel{\mathrm{(ii)}}{\leq} \frac{\sqrt{3}\gamma}{16}\norm{\Delta} \:
\norm{\DelPrime},
\end{align*}
where step (i) follows from inequality~\eqref{EqnSqNormA} in
Lemma~\ref{LemSqNorm}, the independence of $\regnoise$ and $X$, and
the fact that $\Exs[\regnoise^2] = 1$; and step (ii) follows from the
bound~\eqref{EqnMORTechAtwo} in Lemma~\ref{LemMORTechA}.


\subsubsection{Case $\norm{\Delta} > 1$:}
After applying the Cauchy-Schwarz inequality, it suffices show that
$\sqrt{ \Exs \big[ \Delta^2_w(X, Y) \big] } \leq \frac{\gamma}{2}$.
The remainder of this section is devoted to the proof of this claim.\\

Recall the scalar $\tau \defn C_\tau \sqrt{\log \norm{\thetastar}}$,
as well as the events $\eone$ and $\etwo$ from
Lemma~\ref{LemEventBounds}.  For any measurable event $\Event$, define
the function $\Psit(\Event) = \Exs \big[ \Delta^2_w(X, Y) \mid \Event
  \big] \mprob[\Event]$.  With this notation, by successive
conditioning, we have the upper bound
\begin{align}
\label{EqnAMain}
\Exs \big[ \Delta^2_w(X, Y) \big] & \leq \Psit(\eone^c) + \Psit(\eone
\cap \etwo^c) + \Psit(\eone \cap \etwo).
\end{align}
We control each of these terms in turn.

\paragraph{Controlling term $\Psit(\eone^c)$:}  
Noting that $\sup_{x,y} |\Delta_w(x,y)| \leq 2$ and applying
Lemma~\ref{LemEventBounds}(i), we have $\Psit(\eone^c) \leq 4
\mprob[\eone^c] \; \leq 4 \frac{\|\Delta\|_2}{\|\thetastar\|_2}$.


\paragraph{Controlling term $\Psit(\eone \cap \etwo^c)$:}  
Similarly, Lemma~\ref{LemEventBounds}(ii) implies that
\begin{align*}
\Psit(\eone \cap \etwo^c) & \leq 4 \mprob[ \etwo^c] \leq 4 \Big \{
\frac{\tau}{\norm{\thetastar}} + \frac{\tau}{\norm{\theta_u}} + 2
\MYEXP{- \frac{\tau^2}{2}} \Big \}.
\end{align*}


\paragraph{Controlling term $\Psit(\eone \cap \etwo)$:} 

Conditioned on the event $\eone \cap \etwo$, the
bound~\eqref{EqnDeltaw} implies that \mbox{$|\Delta_w(X,Y)| \leq \exp
  (-\tau^2)$,} and hence $\Psit(\eone \cap \etwo) \leq \MYEXP{ - 2
  \tau^2}$. \\

\vspace*{.1in}

Thus, we have derived bounds on each of the three terms in the
decomposition~\eqref{EqnAMain}: putting them together yields
\begin{align*}
\sqrt{ \Exs \big[ \Delta^2_w(X, Y) \big] } \leq \sqrt{ 4
  \frac{\|\Delta\|_2}{\|\thetastar\|_2} + 4 \big \{
  \frac{\tau}{\norm{\thetastar}} + \frac{\tau}{\norm{\theta_u}} + 2
  \MYEXP{- \frac{\tau^2}{2}} \big \} + \MYEXP{ - 2 \tau^2}}
\end{align*}
By choosing $C_\tau$ sufficiently large in the definition of $\tau$,
selecting the signal-to-noise constant $\snrconstmor$ in
condition~\eqref{EqnMORass} sufficiently large, the claim follows.


\subsubsection{Proof of Lemma~\ref{LemMORTechA}}
\label{AppLemMORTechA}
The lemma statement consists of two inequalities, and we divide our
proof accordingly.

\paragraph{Proof of inequality~\eqref{EqnMORTechAone}:}

For any measurable event $\Event$, let us introduce the function
$\Psi(\Event) \defn \Exs \Big[ \frac{Y^2 \inprod{X}{\theta_u}^2
  }{(\exp(Z_u) + \exp(-Z_u))^4 } \mid \Event \Big] \mprob[\Event]$.
With this notation, successive conditioning yields the decomposition
\begin{align}
\label{EqnTTMain}
\Exs \big[ \frac{Y^2 \inprod{X}{\theta_u}^2 }{(\exp(Z_u) +
    \exp(-Z_u))^4} \big] & = \Psi(\efour^c) + \Psi(\efour \cap
\ethree^c) + \Psi(\etwo),
\end{align}
and we bound each of these terms in turn.  The reader should recall
the constant \mbox{$\tau \defn C_\tau \sqrt{\log \norm{\thetastar}}$},
as well as the events $\ethree$ and $\efour$ from
Lemma~\ref{LemEventBounds}.

\paragraph{Bounding $\Psi(\efour^c)$:}
Observe that
\begin{align} 
\label{EqnWorldCup}
\frac{Y^2 \inprod{X}{\theta_u}^2 }{(\exp(Z_u) + \exp(-Z_u))^4 } \leq
\sup_{t \geq 0} \frac{t^2}{\exp(4t)} \leq \frac{1}{4 e^2},
\end{align} 
where the final step follows from inequality~\eqref{Eqn:x_over_exp}.
Combined with Lemma~\ref{LemEventBounds}(iv), we conclude that
$\Psi(\efour^c) \leq \frac{1}{2e^2} \MYEXP{- \frac{\tau^2}{2}}$.


\paragraph{Bounding $\Psi(\efour \cap \ethree^c)$:}
In this case, we have
\begin{align*}
\Psi(\efour \cap \ethree^c) & \stackrel{(i)}{\leq} \frac{1}{4e^2}
\Prob[\ethree^c] \; \stackrel{(ii)}{\leq} \frac{1}{4e^2} \Big \{
\frac{\tau}{\norm{\thetastar}} + \frac{\tau}{\norm{\theta_u}} \Big \},
\end{align*}
where step (i) follows from inequality~\eqref{EqnWorldCup}, and step
(ii) follows from Lemma~\ref{LemEventBounds}(iii).

 
\paragraph{Bounding $\Psi(\etwo)$: }  

Conditioned on the event $\etwo$, we have $Y^2 \inprod{X}{\theta_u}^2
\geq \frac{\tau^2}{2}$, where we have used the lower
bound~\eqref{EqnMagnitude}.  Introducing the shorthand $t^* =
\tau^2/2$, this lower bound implies that
\begin{align*}
\Psi(\etwo) & \leq \sup_{t \geq t^*} \frac{t^2}{ \MYEXP{4t}} \; \leq
\; \frac{(t^*)^2}{ \MYEXP{4t^*}} \; = \; \frac{\tau^4}{4
  \MYEXP{2\tau^2}},
\end{align*}
where inequality (i) is valid as long as $t^* = \frac{\tau^2}{2} \geq
\frac{1}{2}$, or equivalently $\tau^2 \geq 1$. \\

\vspace*{.06in}

Substituting our upper bounds on three components in the
decomposition~\eqref{EqnTTMain} yields
\begin{align*}
\Exs \big[ \frac{Y^2 \inprod{X}{\theta_u}^2 }{(\exp(Z_u) +
    \exp(-Z_u))^4} \big] & \leq \frac{1}{2e^2} \, \MYEXP{ -
  \frac{\tau^2}{2}} + \frac{1}{4e^2} \Big(
\frac{\tau}{\norm{\thetastar}} + \frac{\tau}{\norm{\theta_u}} \Big) +
\frac{\tau^4}{4} \MYEXP{-2\tau^2}.
\end{align*}
Setting $C_\tau$ sufficiently large in the definition of $\tau$ and
choosing sufficiently large values of the signal-to-noise constant
$\snrconstmor$ in the condition~\eqref{EqnMORass} yields the claim.


\paragraph{Proof of inequality~\eqref{EqnMORTechAtwo}:}

For any measurable event $\Event$, let us introduce the function
\mbox{$\newpsi(\Event) = \Exs \big[ \frac{Y^2 }{(\exp(Z_u) +
      \exp(-Z_u))^4} \mid \Event \big] \mprob[\Event]$.}  Recalling
the event $\efive$ from Lemma~\ref{LemEventBounds}, successive
conditioning yields the decomposition
\begin{align}
\label{EqnTTTMain}
\Exs \big[ \frac{Y^2 }{(\exp(Z_u) + \exp(-Z_u))^4} \big] & =
\newpsi(\efive^c) + \newpsi(\efive).
\end{align}
We bound each of these terms in turn.

\paragraph{Bounding $\newpsi(\efive^c)$: }

Simple algebra combined with Lemma~\ref{LemEventBounds}(v) yields the
upper bound $\newpsi(\efive^c) \leq \frac{ \tau}{16 \norm{\theta_u}}
\Exs [ Y^2 ]$.  Conditioned on $\efive$, we have the upper bound
$|\inprod{X}{\theta_u}| \leq \tau$, whence
\begin{align*} 
\inprod{X}{\thetastar}^2 \leq 2 \tau^2 + 2 \inprod{X}{\Delta}^2.
\end{align*} 
Combining Lemma~\ref{LemCC} with the bound $\norm{\Delta} \leq 1$, we
find that $\inprod{X}{\thetastar}^2 \leq 2 \tau^2 + 4$.  Since $Y
\stackrel{d}{=} (2Z - 1) \inprod{X}{\thetastar} + \regnoise$, we have
\begin{align*}
\Exs[Y^2 \mid \efive^c] \leq \Exs[ 2 \inprod{X}{\thetastar}^2 + 2
  \regnoise^2 \mid \efive^c] \stackrel{\mathrm{(i)}}{\leq} 4 \tau^2 +
10.
\end{align*}
Putting together the pieces, we conclude that $\newpsi(\efive^c) \leq
\frac{ 4 \tau^3 + 10 \tau}{16 \norm{\theta_u}}$.


\paragraph{Bounding $\newpsi(\efive)$:}
Recall that $Z_u = Y \inprod{X}{\theta_u}$, so we have that
\begin{align*}
\newpsi(\efive) \leq \Exs \Big[\frac{Y^2}{ ( \MYEXP{Y
      \inprod{X}{\theta_u} } + \MYEXP{-Y \inprod{X}{\theta_u} })^4}
  \mid \efive \Big] \stackrel{\mathrm{(i)}}{\leq} \frac{4}{ (e \,
  \tau)^2},
\end{align*}
where step (i) follows from the bound~\eqref{Eqn:x_over_exp} and the
observation that $|\inprod{X}{\theta_u}| \geq \tau$ conditioned on the
event $\efive$. \\

Substituting our bounds on the two terms into the
decomposition~\eqref{EqnTTMain} yields
\begin{align*}
\Exs \big[ \frac{Y^2 }{(\MYEXP{Z_u} + \MYEXP{-Z_u})^4} \big] & \leq
\frac{ 4 \tau^3 + 10 \tau}{16 \norm{\theta_u}} + \frac{4}{ (e \,
  \tau)^2} \; \leq \; \frac{ 8 \tau^3 + 20 \tau}{16 \norm{\thetastar}}
+ \frac{4}{ (e \, \tau)^2}.
\end{align*}
Once again, sufficiently large choices of the constant $\TAUCON$ and
the signal-to-noise constant $\snrconstmor$ in
equation~\eqref{EqnMORass} yields the claim.


\subsubsection{Proof of Lemma~\ref{LemEventBounds}}

In this section, we prove the probability bounds on events $\eone$
through $\esix$ stated in Lemma~\ref{LemEventBounds}.  In doing so, we
make use of the following auxiliary result, due to Yi et
al.~\cite{hardem_mor} (see Lemma 1 in their paper):
\blems
\label{LemYCS}
Given vectors $v,z \in \real^\usedim$ and a Gaussian random vector $X
\sim \NORMAL(0, I)$, the matrix $\Sigma = \Exs \big[ X X^T \mid
  \inprod{X}{v}^2 > \inprod{X}{z}^2 \big]$ has singular values
\begin{subequations}
\begin{align}
\label{EqnYCSA}
\Big(1 + \frac{\sin \alpha}{\alpha} , 1 - \frac{\sin \alpha}{\alpha},
1,\ldots,1\Big), \qquad \mbox{ where $\alpha = \cos^{-1} \frac{
    \inprod{z-v}{z+v}}{\norm{ z+v}\norm{z - v}}$.}
\end{align}
Moreover, whenever $\norm{v} \leq \norm{z}$, we have
\begin{align}
\label{EqnYCSB}
\mprob \big[ \inprod{X}{v}^2 > \inprod{X}{z}^2 \big] & \leq
\frac{\norm{v}}{\norm{z}}.
\end{align}
\end{subequations}
\elems

\paragraph{Proof of Lemma~\ref{LemEventBounds}(i):}
Note that the event $\eone^c$ holds if and only if
$\inprod{X}{\thetastar} \inprod{X}{\theta_u} < 0$, or equivalently, if
and only if
\begin{align*}
4 \inprod{X}{\thetastar}\inprod{X}{\theta_u} = \inprod{X}{\thetastar +
  \theta_u}^2 - \inprod{X}{\thetastar-\theta_u}^2 < 0.
\end{align*}
Now observe that
\begin{align*}
\norm{\thetastar - \theta_u} & \leq u \norm{ \Delta} \leq \norm{\Delta},
\quad \mbox{and} \quad \norm{ \thetastar + \theta_u} & \geq 2
\norm{\thetastar} - \norm{\Delta} \geq \norm{\thetastar} \geq \;
\norm{\Delta}.
\end{align*}
Consequently, we may apply the bound~\eqref{EqnYCSB} from
Lemma~\ref{LemYCS} with $v = \thetastar + \theta_u$ and $z =
\thetastar - \theta_u$ to obtain $\mprob[\eone^c] \leq
\frac{\|\thetastar - \theta_u\|_2}{\|\thetastar + \theta_u\|_2} \;
\leq \; \frac{\|\Delta\|_2}{\|\thetastar\|_2}$, as claimed.

\paragraph{Proof of Lemma~\ref{LemEventBounds}(iv):}
For $X \sim \NORMAL(0,\sigma^2)$, we have $\mprob \big[|X| \leq \tau
  \big] \leq 2 \exp \MYEXP{- \frac{\tau^2}{2 \sigma^2}}$ \mbox{for any
  $\tau \geq 0$}, from which the claim follows.

\paragraph{Proof of Lemma~\ref{LemEventBounds}(v):}
For $X \sim \NORMAL(0, \sigma^2)$, we have
\begin{align}
\label{EqnNormLower}
\mprob \big[|X| \leq \tau \big] & \leq \sqrt{ \frac{2}{\pi} }
\frac{\tau}{\sigma} \qquad \mbox{for any $\tau \geq 0$}
\end{align} 
from which the claim follows.

\paragraph{Proof of Lemma~\ref{LemEventBounds}(vi):}

Similarly, this inequality follows from the tail
bound~\eqref{EqnNormLower}.


\paragraph{Proof of Lemma~\ref{LemEventBounds}(iii):}
This claim follows from parts (v) and (vi) of
Lemma~\ref{LemEventBounds}, combined with the union bound.

\paragraph{Proof of Lemma~\ref{LemEventBounds}(ii):}
This bound follows from parts (iii) and (iv) of
Lemma~\ref{LemEventBounds}, combined with the union bound.


\subsubsection{Proof of Lemma~\ref{LemCC}}
\label{AppLemCC}

For an event $\Event$, define the matrix $\Gamma(\Event) = \Exs[ X X^T
  \mid \Event]$.  The lemma concerns the operator norm of this matrix
for different choices of the event $\Event$.

\paragraph{Conditioned on $\eone \cap \etwo$: } 
In this case, we write
\begin{align*}
\Exs \big[ X X^T ] & = \Gamma(\eone \cap \etwo) \mprob[\eone \cap
  \etwo] + \Gamma \big((\eone \cap \etwo)^c \big) \Prob[ (\eone \cap
  \etwo)^c] \; \; \succeq \; \Gamma(\eone \cap \etwo) \; \mprob[\eone
  \cap \etwo].
\end{align*}
Since $\Exs[X X^T] = I$, we conclude that $\opnorm{\Gamma(\eone \cap
  \etwo)} \leq \frac{1}{\mprob[\eone \cap \etwo]}$, and hence it
suffices show that $\mprob[\eone \cap \etwo] \geq \frac{1}{2}$.  Parts
(i) and (ii) of Lemma~\ref{LemEventBounds} imply that
\begin{align*}
\mprob[\eone \cap \etwo] \geq 1- \frac{\|\Delta\|_2}{\|\thetastar\|_2}
- \frac{\tau}{\norm{\thetastar}} - \frac{\tau}{\norm{\theta_u}} - 2
\MYEXP{ - \frac{\tau^2}{2}}.
\end{align*}
For appropriate choices of $\TAUCON$ and the constant $\snrconstmor$
in the signal-to-noise condition~\eqref{EqnMORass}, the claim follows.

\paragraph{Conditioned on $\eone^c$: } As before, 
note that the event $\eone^c$ holds if and only if the inequality
\mbox{$|\inprod{X}{\thetastar + \theta_u}| <
  |\inprod{X}{\thetastar-\theta_u}|$} holds.  Consequently,
Lemma~\ref{LemYCS} implies that $\opnorm{\Gamma(\eone^c)} \leq 2$.

\paragraph{Conditioned on $\efive^c$: } 
We make note of an elementary fact about Gaussians: for any scalar
$\alpha > 0$ and unit norm vector $\|v\|_2 = 1$, for $X \sim \NORMAL(0,
I_\usedim)$, we have
\begin{align}
\label{eqn::mor_cc2}
\opnorm{\Exs \big[ XX^T \mid \, | \inprod{X}{v} | \leq \alpha \big]} &
\leq \max \big(1, \alpha^2 \big).
\end{align}
In particular, when $\alpha \leq 1$, then the operator norm is at most
$1$.  This claim follows easily from the rotation invariance of the
Gaussian, which allows us to assume that $v = e_1$ without loss of
generality.  It is thus equivalent to bound the largest eigenvalue of
the matrix
\begin{align*}
\Dmat & \defn \Exs \big[ XX^T \mid | X_1 | \leq \alpha \big],
\end{align*}
which is a diagonal matrix by independence of the entries of $X$.
Noting that $D_{11} \leq \alpha^2$ and $\Dmat_{jj} = 1$ for $j \neq 1$
completes the proof of the bound~\eqref{eqn::mor_cc2}.

Applying the bound~\eqref{eqn::mor_cc2}, we find that
$\opnorm{\Gamma(\efive^c)} \leq \max \Big( 1, \frac{ \tau^2
}{\norm{\theta_u}^2} \Big)$.  Consequently, the claim follows by
making sufficiently large choices of $\TAUCON$ and the constant
$\snrconstmor$ in the signal-to-noise condition~\eqref{EqnMORass}.

\paragraph{Conditioned on $\esix^c$: } The bound~\eqref{eqn::mor_cc2} 
implies that $\opnorm{\Exs \big[ XX^T \mid \esix^c]} \leq \max \big \{
1, \frac{ \tau^2 }{\norm{\thetastar}^2} \big \}$.  As in the previous
case, choosing $\TAUCON$ and $\snrconstmor$ appropriately ensures that
$\frac{ \tau^2 }{\norm{\thetastar}^2} \leq 1$.


\subsection{Proof of Corollary~\ref{CorMORSamp}}
\label{AppCorMORSamp}

We need to compute an upper bound on the function $\rateem(n,\delta)$
previously defined in equation~\eqref{EqnRateEM}.  For this particular
model, we have
\begin{align*}
\norm{\mem{\rvec} - \memsamp{\rvec}} = \norm{ \Big( \sum_{i=1}^n x_i
  x_i^T \Big)^{-1} \Big( \sum_{i=1}^\numobs (2 w_{\rvec}(x_i,y_i) -
  1)y_i x_i\Big) - 2 \Exs [w_{\theta}(X,Y) YX]}.
\end{align*}
Define the matrices $\SigHat \defn \fsavgtwo{1}{} x_i x_i^T$ and $\Sig
= \Exs[X X^T] = I$, as well as the vector
\begin{align*}
\qvechat \defn \fsavgtwo{1}{} \big[\mm(x_i,y_i) y_i x_i \big], \quad
\mbox{and} \quad \qvec\defn \Exs \big[\mm(X,Y) Y X \big],
\end{align*}
where $\mm(x,y) \defn 2w_{\theta}(x,y) - 1$.  Noting that $\Exs[YX] =
0$, some straightforward algebra then yields the bound
\begin{align}
\label{eqn::mor_rate_bound}
\norm{\mem{\rvec} - \memsamp{\rvec}} \leq \underbrace{
  \opnorm{\SigHat^{-1}} \norm{\qvechat - \qvec}}_{T_1} +
\underbrace{\opnorm{ \SigHat^{-1} - \Sig^{-1}} \norm{\qvec}}_{T_2}.
\end{align}
We bound each of the terms $T_1$ and $T_2$ in turn.

\paragraph{Bounding $T_1$:} 
Recall the assumed lower bound on the sample size---namely $\numobs >
c \, \usedim \log(1/\delta)$ for a sufficiently large constant
$c$. Under this condition, standard bounds in random matrix
theory~\cite{vershynin_nonasymp}, guarantee that $\opnorm{\SigHat -
  \Sigma} \leq \frac{1}{2}$ with probability at least $1 - \delta$.
When this bound holds, we have $\opnorm{\SigHat^{-1}} \geq 1/2$.

As for the other part of $T_1$, let us write $\norm{\qvechat - \qvec}
= \sup_{u \in \Sphere{\usedim}} Z(u)$, where 
\begin{align*}
Z(u) \defn \frac{1}{\numobs} \sum_{i=1}^\numobs \mm(x_i,y_i) y_i
\inprod{x}{u} - \Exs[\mm(X,Y) Y\inprod{X}{u} ].
\end{align*}
By a discretization argument over a $1/2$-cover of the sphere
$\Sphere{\usedim}$---say $\{u^1, \ldots, u^M \}$---we have the upper
bound $\norm{\qvechat - \qvec} \leq 2 \max_{j \in [M]} Z(u^j)$.  Thus,
it suffices to control the random variable $Z(u)$ for a fixed $u \in
\Sphere{\usedim}$.  By a standard symmetrization
argument~\cite{vanderVaart96}, we have
\begin{align*}
\mprob \big[ Z(u) \geq t \big] & \leq 2 \mprob \Big[ \frac{1}{\numobs}
  \sum_{i=1}^\numobs \rade{i} \mm(x_i, y_i) y_i \inprod{x_i}{u} \geq
  t/2 \Big],
\end{align*}
where $\{ \rade{i} \}_{i=1}^\numobs$ are an i.i.d. sequence of
Rademacher variables.  Let us now define the event \mbox{$\Event \big
  \{ \frac{1}{\numobs} \sum_{i=1}^\numobs \inprod{x_i}{u}^2 \leq 2
  \}$.}  Since each variable $\inprod{x_i}{u}$ is sub-Gaussian with
parameter one, standard tail bounds imply that $\mprob[\Event^c] \leq
\MYEXP{-\numobs/32}$.  Therefore, we can write
\begin{align*}
\mprob \big[ Z(u) \geq t \big] & \leq 2 
\mprob \Big[ \frac{1}{\numobs}
  \sum_{i=1}^\numobs \rade{i} \mm(x_i, y_i) y_i \inprod{x_i}{u} \geq
  t/2 \mid \Event \Big] + 2 \MYEXP{-\numobs/32}.
\end{align*}
As for the remaining term, we have
\begin{align*}
\Exs \Big[ \exp \Big( \frac{\lambda}{\numobs} \sum_{i=1}^\numobs
  \rade{i} \mm(x_i, y_i) y_i \inprod{x_i}{u} \Big) \mid \Event \Big] &
\leq \Exs \Big[ \exp \Big( \frac{2 \lambda}{\numobs}
  \sum_{i=1}^\numobs \rade{i} y_i \inprod{x_i}{u} \Big) \mid \Event
  \Big],
\end{align*}
where we have applied the Ledoux-Talagrand contraction for Rademacher
processes~\cite{LedTal91,Koltchinskii}, using the fact that $|\mm(x,y)|
\leq 1$ for all pairs $(x,y)$.  Now conditioned on $x_i$, the random
variable $y_i$ is zero-mean and sub-Gaussian with parameter at most
$\sqrt{\|\thetastar\|_2^2 + \sigma^2}$.  Consequently, taking
expectations over the distribution $(y_i \mid x_i)$ for each index
$i$, we find that
\begin{align*}
\Exs \Big[ \exp \Big( \frac{2 \lambda}{\numobs} \sum_{i=1}^\numobs
  \rade{i} y_i \inprod{x_{i}}{u} \Big) \mid \Event \Big] & \leq \Big[\exp \Big(
\frac{4\lambda^2}{\numobs^2} \big( \|\thetastar\|_2^2 + \sigma^2 \big)
\sum_{i=1}^\numobs \inprod{x_{i}}{u}^2  \Big)\mid \Event\Big] \; \\ &\leq 
\exp \Big(\frac{8\lambda^2}{\numobs} \big( \|\thetastar\|_2^2 + \sigma^2 \big)
\Big),
\end{align*}
where the final inequality uses the definition of $\Event$.  Using
this bound on the moment-generating function, we find that
\begin{align*}
\mprob \Big[ \frac{1}{\numobs} \sum_{i=1}^\numobs \rade{i} \mm(x_i,
  y_i) y_i \inprod{x_i}{u} \geq t/2 \mid \Event \Big] & \leq \exp \Big ( -
\frac{\numobs t^2}{256 (\|\thetastar\|_2^2 + \sigma^2)} \Big).
\end{align*}
Since the $1/2$-cover of the unit sphere $\Sphere{\usedim}$ has at
most $2^\usedim$ elements, we conclude that there is a universal
constant $c$ such that $T_1 \leq c \, \sqrt{\|\thetastar\|_2^2 +
  \sigma^2} \sqrt{\frac{\usedim}{\numobs} \, \log(1/\delta)}$ with
probability at least $1-\delta$.

\paragraph{Bounding $T_2$:} 
Since $\numobs > \usedim$ by assumption, standard results in random
matrix theory~\cite{vershynin_nonasymp} imply that $\opnorm{ \SigHat^{-1} -
  \Sig^{-1}} \leq c \sqrt{\frac{\usedim}{\numobs} \, \log (1/\delta)}$
with probability at least $1 -\delta$.  On the other hand, observe
that
\begin{align*}
\norm{\qvec} & = \norm{\mem{\theta}} \leq 2 \norm{\thetastar},
\end{align*}
since the population operator $M$ is a contraction, and $\|\theta\|_2
\leq 2 \|\thetastar\|_2$.  Combining the pieces, we see that $T_2 \leq
c \norm{\thetastar} \sqrt{ \frac{\usedim}{\numobs} \; \log
  (1/\delta)}$ with probability at least $1 -\delta$. \\

\noindent Finally, substituting our bounds on $T_1$ and $T_2$ into the
decomposition~\eqref{eqn::mor_rate_bound} yields the claim.


\subsection{Proof of Corollary~\ref{CorMORSGD}}
\label{AppCorMORSGD}

We need to bound the uniform variance $\sigmaSGD^2 = \sup_{\theta \in
  \Balltwor} \Exs \norm{ \nabla Q_1 (\theta| \theta)}^2$, where $r =
\frac{\|\thetastar\|_2}{\deltaconstmor}$.  From the gradient
update~\eqref{EqnGradientMOR}, we have $\nabla Q_1 (\theta \mid
\theta) = (2 w_\theta(x_1,y_1) - 1) y_1 x_1 - \inprod{x_1}{\theta}
x_1$, and hence
\begin{align}
\label{EqnReturnHere}
\Exs \big[ \norm{\nabla Q_1 (\theta| \theta)}^2 \big] & \leq 2
\underbrace{\Exs [y_1^2 \norm{x_1}^2]}_{T_1} + 2
\underbrace{\opnorm{\Exs [x_1 x_1^T \norm{x_1}^2]}}_{T_2}
\norm{\theta}^2.
\end{align}
First considering $T_1$, recall that $y_1 = z_1
\inprod{x_1}{\thetastar} + \regnoise_1$, where $\regnoise \sim
\NORMAL(0,\sigma^2)$ and $z_1$ is a random sign, independent of
$(x_1, \regnoise_1)$.   Consequently, we have
\begin{align*}
T_1 & \leq 2 \Exs[\inprod{x_1}{\thetastar}^2 \norm{x_1}^2] + 2
\Exs[\regnoise_1^2 \norm{x_1}^2] \; \leq \; 2
\sqrt{\Exs[\inprod{x_1}{\thetastar}^4]} \; \sqrt{\Exs[\norm{x_1}^4]} +
2 \sigma^2 \usedim,
\end{align*}
where we have applied the Cauchy-Schwarz inequality, and observed that
$\Exs[\norm{x_1}^2] = \usedim$ and $\Exs[\regnoise_1^2] = \sigma^2$.
Since the random variable $\inprod{x_1}{\thetastar}$ is sub-Gaussian
with parameter at most $\|\thetastar\|_2$, we have
$\Exs[\inprod{x_1}{\thetastar}^4] \leq 3 \|\thetastar\|_2^4$.
Moreover, since the random vector $x_1$ has i.i.d.  components, we
have
\begin{align*}
\Exs[\|x_1\|_2^4] & = \sum_{j=1}^\usedim \Exs[x_{1j}^4] + 2 \sum_{i
  \neq j} \Exs[x_{1i}^2] \Exs[x_{1j}^2] = 3 \usedim + 2 {\usedim
  \choose 2} \; \leq 4 \usedim^2.
\end{align*}
Putting together the pieces, we conclude that $T_1 \leq 8
\|\thetastar\|_2^2 \usedim + 2 \sigma^2 \usedim$.

Turning to term $T_2$, by definition of the operator norm, there is a
unit-norm vector $u \in \real^\usedim$ such that
\begin{align*}
T_2 = \opnorm{\Exs [x_1 x_1^T \norm{x_1}^2]} \; = u^T \Big( \Exs [x_1
  x_1^T \norm{x_1}^2] \Big) u & = \; \Exs[ \inprod{x_1}{u}^2
  \norm{x_1}^2] \\
& \stackrel{\mathrm{(i)}}{\leq} \sqrt{\Exs[\inprod{x_1}{u}^4]} \;
\sqrt{\Exs[\|x_1\|_2^4]} \\
& \stackrel{\mathrm{(ii)}}{\leq} \sqrt{3} \sqrt{4 \usedim^2} \; \leq \; 4 \usedim.
\end{align*}
where step (i) applies the Cauchy-Schwarz inequality, and step (ii)
uses the fact that $\inprod{x_1}{u}$ is sub-Gaussian with parameter
$1$, and our previous bound on $\Exs[\|x_1\|_2^4]$. \\

\noindent Putting together the pieces yields $\sigmaSGD^2 \leq c \,
(\sigma^2 + \|\thetastar\|_2^2) \usedim$, so that
Corollary~\ref{CorMORSGD} follows as a consequence of
Theorem~\ref{ThmStochasticEM}.


\section{Proofs for missing covariates}

In this appendix, we provide proofs of results related to regression
with missing covariates, as presented in
Section~\ref{SecMissingAnalysis}.  More specifically, we first prove
Corollary~\ref{CorMissingPop} on the population level behavior,
followed by the proof of Corollaries~\ref{CorMissingSamp}
and~\ref{CorMissingSGD} on the behavior of sample-splitting EM updates
and stochastic gradient EM updates, respectively.


\subsection{Proof of Corollary~\ref{CorMissingPop}}
\label{AppCorMissingPop}

We need to verify the conditions of Theorem~\ref{ThmGradEM}, namely
that the function $\specq$ is $\smoothparam$-smooth,
$\strongparam$-strongly concave, and that the
GS condition is satisfied. In this case, $\specq$ is a
quadratic of the form
\begin{align*}
\specq(\theta) = \frac{1}{2} \inprod{\theta}{\Exs
  \big[\CovMat_{\thetastar}(\Xobs,Y)\big] \theta}- \inprod{\Exs \big[Y
    \mu_{\thetastar}(\Xobs,Y) \big]}{\theta},
\end{align*}
where the vector $\mu_{\fp} \in \real^\usedim$ and matrix
$\CovMat_{\thetastar}$ were previously defined (see
equations~\eqref{EqnMuMissing} and~\eqref{EqnCovMissing}
respectively). Here the expectation is over both the patterns of
missingness and the random $(\Xobs,Y)$.

\paragraph{Smoothness and strong concavity:}
Note that $\specq$ is a quadratic function with Hessian $\nabla^2
\specq(\theta) = \Exs \big[\CovMat_{\thetastar}(\Xobs,Y) \big]$.  Let
us fix a pattern of missingness, and then average over $(\Xobs, Y)$.
Recalling the matrix $\mism_{\mopt}$ from
equation~\eqref{EqnMuMissingAux}, we find that yields
\begin{align*}
\Exs \big[ \CovMat_{\thetastar}(\Xobs,Y) \big] & = \begin{bmatrix} I &
  \mism_{\thetastar} \begin{bmatrix} I \\ \bobs^{*T} \end{bmatrix}
  \\ \begin{bmatrix} I & \bobs^* \end{bmatrix} \mism_{\thetastar}^T &
  I \end{bmatrix} \; = \; \begin{bmatrix} I & 0 \\ 0 & I 
\end{bmatrix},
\end{align*}
showing that the expectation does not depend on the pattern of
missingness.  Consequently, the quadratic function $\specq$ has an
identity Hessian, showing that smoothness and strong concavity hold
with $\smoothparam = \strongparam = 1$.

\paragraph{Condition \gsone:}
We need to prove the existence of a scalar $\gamma \in [0,1)$ such
  that \mbox{$\norm{\Exs [\mismain]} \leq \gamma \norm{ \theta -
      \thetastar}$,} where the vector $\mismain = \mismain(\theta,
  \thetastar)$ is given by
\begin{align}
\mismain & \defn \CovMat_{\thetastar}(\Xobs,Y) \theta -
Y\mu_{\thetastar}(\Xobs,Y) - \CovMat_{\theta}(\Xobs,Y) \theta +
Y\mu_{\theta}(\Xobs,Y).
\end{align}
For a fixed pattern of missingness, we can compute the expectation
over $(\Xobs, Y)$ in closed form.  Supposing that the first block is
missing, we have
\begin{align}
\label{EqnMisExp}
\Exs_{\Xobs, Y}[\mismain] & =
\begin{bmatrix} (\bmis - \bmis^*) + \mytemp_1 \bmis \\
\mytemp_2 (\bobs - \bobs^*)
\end{bmatrix}.
\end{align}
where $\mytemp_1 \defn \frac{\norm{ \bmis^*}^2 - \norm{\bmis}^2 +
  \norm{\bobs - \bobs^*}^2}{\norm{\bmis}^2 + \sigma^2}$ and $\mytemp_2
\defn \frac{ \norm{\bmis}^2}{\norm{\bmis}^2 + \sigma^2}$.  We claim
that these scalars can be bounded, independently of the missingness
pattern, as
\begin{align}
\label{EqnAlphaMissing}
\mytemp_1 & \leq 2 (\ccon + \ctwo) \frac{\norm{ \theta -
    \thetastar}}{\sigma}, \qquad \mbox{and} \quad \mytemp_2 \leq \SIVA
\defn \frac{1}{1 + \big(\frac{1}{\ccon + \ctwo} \big)^2} < 1.
\end{align}
Taking these bounds~\eqref{EqnAlphaMissing} as given for the moment,
we can then average over the missing pattern. Since each coordinate is
missing independently with probability $\missing$, the expectation of
the $i^{\mathrm{th}}$ coordinate is at most $\big|\Exs[\mismain]|_i \leq \big|
\missing |\theta_i - \thetastar_i| + \missing \mytemp_1 |\theta_i| +
(1-\missing) \mytemp_2 |\theta_i - \theta_i^*|\big|$.  Thus, defining $\eta
\defn (1-\missing) \SIVA + \missing < 1$, we have
\begin{align*}
\|\Exs[\mismain]\|_2^2 & \leq \eta^2 \|\theta - \thetastar\|_2^2 +
\missing^2 \mytemp_1^2 \|\theta\|_2^2 + 2 \mytemp_1 \eta \missing
|\inprod{\theta}{\theta - \thetastar}| \\
& \leq \underbrace{\Big \{ \eta^2 + \missing^2 \|\theta\|_2^2 \frac{4
    \, (\ccon + \ctwo)^2}{\sigma^2} + \frac{4 \eta \missing
    \|\theta\|_2 (\ccon + \ctwo)}{\sigma} \Big \}}_{\gamma^2} \|\theta -
\thetastar\|_2^2,
\end{align*}
where we have used our upper bound~\eqref{EqnAlphaMissing} on
$\mytemp_1$.  We need to ensure that $\gamma < 1$. By assumption, we
have $\norm{\thetastar} \leq \ccon \sigma$ and $\norm{\theta -
  \thetastar} \leq \ctwo \sigma$, and hence $\norm{\theta} \leq (\ccon
+ \ctwo) \sigma$.  Thus, the coefficient $\gamma^2$ is upper bounded as
\begin{align*}
\gamma^2 & \leq \eta^2 + 4 \missing^2 \, (\ccon + \ctwo)^4 + 4 \eta
\missing (\ccon + \ctwo)^2.
\end{align*}
 Under the stated conditions of the corollary, we have $\gamma < 1$,
 thereby completing the proof. \\
 
It remains to prove the bounds~\eqref{EqnAlphaMissing}.  By our
assumptions, we have $\norm{\bmis}- \norm{ \bmis^*} \leq \norm{ \bmis
  - \bmis^*}$, and moreover
\begin{align}
\label{EqnBmisBound}
\norm{ \bmis} \leq \norm{\bmis^*} + \ctwo \sigma \leq (\ccon + \ctwo)
\sigma.
\end{align}
As consequence, we have
\begin{align*}
\norm{\bmis^*}^2 - \norm{\bmis}^2 = ( \norm{\bmis} - \norm{ \bmis^*}
)( \norm{ \bmis} + \norm{ \bmis^*} ) \leq (2\ccon + \ctwo) \sigma
\norm{\bmis - \bmis^*}
\end{align*}
Since $\norm{ \bobs - \bobs^*}^2 \leq \ctwo \sigma \norm{ \bobs -
  \bobs^*}$, the stated bound on $\mytemp_1$ follows.

On the other hand, we have
\begin{align*} 
\mytemp_2 = \frac{\norm{\bmis}^2}{\norm{\bmis}^2 + \sigma^2} & =
\frac{1}{1 + \frac{\sigma^2}{\norm{\bmis}^2}}
\stackrel{\mathrm{(i)}}{\leq} \underbrace{\frac{1}{1 +
    \big(\frac{1}{\ccon + \ctwo} \big)^2 }}_{\SIVA} < 1,
\end{align*}
where step (i) follows from~\eqref{EqnBmisBound}.


\subsection{Proof of Corollary~\ref{CorMissingSamp}}
\label{AppCorMissingSamp}

We need to upper bound the deviation function
$\rategrad(\numobs,\delta)$ previously defined~\eqref{EqnRateGrad}.
For any fixed $\theta \in \Balltwor = \{ \theta \in \real^\usedim
\mid \|\theta - \thetastar\|_2 \leq \ctwo \sigma \}$, we have the
bound $\norm{\mgrad{\mvec} - \mgradsamp{\mvec}} \leq T_1 + T_2$, where
\begin{align*}
T_1 & \defn \norm{ \big[
\Exs \CovMat_{\mvec} (\xobs, y) \mvec -
    \fsavgtwo{1}{} \CovMat_{\mvec}(\xobsi,y_i) \mvec \big]}, \quad
\mbox{and} \\
T_2 & \defn \norm{\big[\Exs (y \mu_{\mvec}(\xobs, y)) -
    \fsavgtwo{1}{} y_i \mu_{\mvec}(\xobsi,y_i) \big]}.
\end{align*} 
For convenience, we let $z_i \in \real^\usedim$ be a
$\{0,1\}$-valued indicator vector, with ones in the positions of
observed covariates.  For ease of notation, we frequently use the
abbreviations $\CovMat_\theta$ and $\mu_\theta$ when the arguments are
understood. We use the notation $\odot$ to denote the element-wise product.

\paragraph{Controlling $T_1$:}
Define the matrices $\bar{\Sigma} = \Exs[\CovMat_\theta(\xobs, y)]$
and $\SigHat = \frac{1}{\numobs} \sum_{i=1}^\numobs
\CovMat_{\mvec}(\xobsi, y_i)$.  With this notation, we have
\begin{align*}
T_1 & \leq \opnorm{\bar{\Sigma} - \SigHat} \; \|\theta\|_2 \; \leq \;
\opnorm{\bar{\Sigma} - \SigHat} \; (\ccon + \ctwo) \, \sigma,
\end{align*}
where the second step follows since any vector $\theta \in \Balltwor$
has $\ell_2$-norm bounded as $\|\theta\|_2 \leq (\ccon +
\ctwo)\sigma$.  We claim that for any fixed vector $u \in
\Sphere{\usedim}$, the random variable $\inprod{u}{(\bar{\Sigma} -
  \SigHat)u}$ is zero-mean and sub-exponential.  When this tail
condition holds and $\numobs > \usedim$, standard arguments in random
matrix theory~\cite{vershynin_nonasymp} ensure that
$\opnorm{\bar{\Sigma} - \SigHat} \leq c \sqrt{\frac{\usedim}{\numobs}
  \, \log(1/\delta)}$ with probability at least $1-\delta$.

It is clear that $\inprod{u}{(\bar{\Sigma} - \SigHat) u}$ has zero mean.
It remains to prove that $\inprod{u}{(\bar{\Sigma} - \SigHat) u}$ is
sub-exponential.  Note that $\SigHat$ is a rescaled sum of rank one
matrices, each of the form
\begin{align*}
\CovMat_{\mvec}(\xobs,y) = I_{\mathrm{mis}} +  \mu_{\mvec} \mu_{\mvec}^T - ((\vecone - z) \odot \mu_{\mvec})((\vecone - z) 
\odot \mu_{\mvec})^T,
\end{align*}
where $I_{\mathrm{mis}}$ denotes the identity matrix on the diagonal sub-block
corresponding to the missing entries.
The square of any sub-Gaussian random variable has sub-exponential
tails.  Thus, it suffices to show that each of the random variables
$\inprod{\mu_\mvec}{u}$,  and
$\inprod{(\vecone - z) \odot \mu_\mvec}{u}$ are sub-Gaussian.  The
random vector $z \odot x$ has i.i.d. sub-Gaussian components 
with parameter at most $1$ and
$\|u\|_2 = 1$, so that $\inprod{z \odot x}{u}$ is sub-Gaussian
with parameter at most $1$.  It remains to verify that $\mu_\mvec$ is
sub-Gaussian, a fact that we state for future reference as a lemma:
\blems
\label{LemMuSubgauss}
Under the conditions of Corollary~\ref{CorMissingPop}, the random
vector $\mu_\mvec(\xobs, y)$ is sub-Gaussian with a constant
parameter.
\elems
\begin{proof}
Introducing the shorthand $\myinter = (\vecone - z) \odot \mvec$, we
have
\begin{align*}
\mu_{\mvec}(\xobs,y) = z \odot x + \frac{1}{ \sigma^2 + \norm{
    \myinter }^2 } \big[ y - \inprod{z \odot \mvec }{z \odot x }\big]
\myinter.
\end{align*}
Moreover, since $y = \inprod{x}{\thetastar} + \regnoise$, we have
\begin{align*}
\inprod{\mu_{\mvec}(\xobs,y)}{u} & = \underbrace{\inprod{z \odot
    x}{u}}_{B_1} + \underbrace{ \frac{\inprod{x}{\myinter}
    \inprod{\myinter}{u}}{ \sigma^2 + \norm{ \myinter}^2 } }_{B_2} +
\underbrace{\frac{\inprod{x}{\thetastar - \theta}
    \inprod{\myinter}{u}}{ \sigma^2 + \norm{ \myinter }^2 }}_{B_3} +
\underbrace{\frac{\regnoise \inprod{\myinter}{u}}{ \sigma^2 + \norm{
      \myinter }^2 }}_{B_4}.
\end{align*}
It suffices to show that each of the variables $\{B_j\}_{j=1}^4$ is
sub-Gaussian with a constant parameter.  As discussed previously, the
variable $B_1$ is sub-Gaussian with parameter at most one. On the
other hand, note that $x$ and $\myinter$ are independent.  Moreover,
with $\myinter$ fixed, the variable $\inprod{x}{\myinter}$ is
sub-Gaussian with parameter $\|\myinter\|_2^2$, whence
\begin{align*}
\Exs[ \MYEXP{\lambda B_2}] & \leq \exp \Big(\lambda^2
\frac{\|\myinter\|_2^2 \inprod{\myinter}{u}^2}{2 (\sigma^2 +
  \|\myinter\|_2^2)^2} \Big) \; \leq e^{\frac{\lambda^2}{2}},
\end{align*}
where the final inequality uses the fact that $\inprod{\myinter}{u}^2
\leq \|\myinter\|_2^2$.  We have thus shown that $B_2$ is sub-Gaussian
with parameter one. Since $\|\theta - \thetastar\|_2 \leq \ctwo
\sigma$, the same argument shows that $B_3$ is sub-Gaussian with
parameter at most $\ctwo$.  Since $\regnoise$ is sub-Gaussian with
parameter $\sigma$ and independent of $\myinter$, the same argument
shows that $B_4$ is sub-Gaussian with parameter at most one, thereby
completing the proof of the lemma.
\end{proof}

\paragraph{Controlling $T_2$:}  We now turn to the second term.
Note the variational representation
\begin{align*}
T_2 & = \sup_{\|u\|_2 = 1} \Big| \Exs\big[y \inprod{\mu_\mvec(\xobs,
    y)}{u} \big] - \fsavgtwo{1}{} y_i
\inprod{\mu_{\mvec}(\xobsi,y_i)}{u} \Big |.
\end{align*}
By a discretization argument--say with a $1/2$ cover $\{u^1, \ldots,
u^M\}$ of the sphere with $M \leq 2^\usedim$ elements---we obtain
\begin{align*}
T_2 & \leq 2 \max_{j \in [M]} \Big| \Exs\big[y
  \inprod{\mu_\mvec(\xobs, y)}{u^j} \big] - \fsavgtwo{1}{} y_i
\inprod{\mu_{\mvec}(\xobsi,y_i)}{u^j} \Big |.
\end{align*}
Each term in this maximum is the product of two zero-mean variables,
namely $y$ and $\inprod{\mu_\mvec}{u}$.  On one hand, the variable $y$
is sub-Gaussian with parameter at most $\sqrt{\|\thetastar\|_2^2 +
  \sigma^2} \leq c \sigma$; on the other hand,
Lemma~\ref{LemMuSubgauss} guarantees that $\inprod{\mu_\mvec}{u}$ is
sub-Gaussian with constant parameter.  The product of any two
sub-Gaussian variables is sub-exponential, and thus, by standard
sub-exponential tail bounds~\cite{BulKoz}, we have
\begin{align*}
\mprob[T_2 \geq t] & \leq 2 M \exp \biggr ( - c \, \min \Big \{
\frac{\numobs t}{\sqrt{1+ \sigma^2}}, \frac{\numobs t^2}{1 + \sigma^2}
  \Big \} \biggr).
\end{align*}
Since $M \leq 2^\usedim$ and $\numobs > c_1 \usedim$, we conclude that
$T_2 \leq c \sqrt{1 + \sigma^2} \, \sqrt{\frac{\usedim}{\numobs} \;
  \log(1/\delta)}$ with probability at least $1 - \delta$. \\

\vspace*{.04in}

Combining our bounds on $T_1$ and $T_2$, we conclude that
$\rategrad(\numobs,\delta) \leq c \sqrt{1 + \sigma^2} \;
\sqrt{\frac{\usedim}{\numobs} \; \log(1/\delta)}$ with probability at
least $1-\delta$.  Thus, we see that Corollary~\ref{CorMissingSamp}
follows from Theorem~\ref{ThmSampleEM}.



\subsection{Proof of Corollary~\ref{CorMissingSGD}}
\label{AppCorMissingSGD}

Once again we focus on bounding the uniform variance $\sigmaSGD^2$.
From the form of $Q$ given in equation~\eqref{EqnMissingQsamp} (with
$\numobs = 1$), we have
\begin{align}
\label{EqnNetArg}
\Exs \Big[ \norm{\nabla Q_1(\theta| \theta)}^2 \Big] & \leq 2 \Big \{
\underbrace{\Exs \Big[ \| \CovMat_{\theta}(\xobs,y) \theta
    \|_2^2}_{T_1} + \underbrace{\Exs[ y^2 \; \| \mu_\theta(\xobs,
      y)\|_2^2] \Big]}_{T_2} \Big \}.
\end{align}
We bound each of these terms in turn.  To simplify notation, we omit
the dependence of $\mu_\theta$ and $\CovMat_\theta$ on $(\xobs, y)$,
but it should be implicitly understood.

\paragraph{Bounding $T_1$:}
Letting $\vecone \in \real^\usedim$ be the vector of 
all ones, and $z \in
\real^\usedim$ be an indicator of observed indices, we have
$\CovMat_{\theta}  =  I_{\mathrm{mis}} +  \mu_{\mvec} \mu_{\mvec}^T - ((\vecone - z) \odot \mu_{\mvec}) ((\vecone - z) 
\odot \mu_{\mvec})^T.$
Consequently,
\begin{align*}
\frac{1}{3} \Exs[\|\CovMat_\theta \theta\|_2^2] & \leq \|\theta\|_2^2
+ \Exs \big[ \|\mu_\mvec\|_2^2 \; \inprod{\mu_\mvec}{\theta}^2 \big] +
\Exs \big[\| (\vecone - z) \odot \mu_\mvec\|_2^2 \; \inprod{(\vecone -
    z) \odot \mu_\mvec}{\theta}^2\big].
\end{align*}
By the Cauchy-Schwarz inequality, we have
\begin{align*}
\Exs \big[ \|\mu_\mvec\|_2^2 \; \inprod{\mu_\mvec}{\theta}^2 \big] &
\leq \sqrt{\Exs[\|\mu_\mvec\|_2^4]} \;
\sqrt{\Exs[\inprod{\mu_\mvec}{\theta}^4]}.
\end{align*}
From Lemma~\ref{LemMuSubgauss}, the random vector $\mu_\theta$ is
sub-Gaussian with constant parameter, so that $\Exs[\|\mu_\mvec\|_2^4]
\leq c \, \usedim^2$.  Since $\|\theta\|_2 \leq c \|\thetastar\|_2$,
the random variable $\inprod{\mu_\mvec}{\theta}$ is sub-Gaussian with
parameter $c \|\thetastar\|_2$, and hence
$\Exs[\inprod{\mu_\mvec}{\theta}^4] \leq c \, \|\thetastar\|_2^4$.
Putting together the pieces, we see that $\Exs \big[ \|\mu_\mvec\|_2^2
  \; \inprod{\mu_\mvec}{\theta}^2 \big] \leq c \, \usedim
\|\thetastar\|_2^2$.  A similar argument applies to other expectation,
so that we conclude that $T_1 = \Exs \big[ \|\CovMat_\theta
  \theta\|_2^2] \leq c \, \|\thetastar\|_2^2 \usedim$, a bound that
holds uniformly for all $\theta \in \Balltwor$.

\paragraph{Bounding $T_2$:}
By the Cauchy-Schwarz inequality, we have
\begin{align*}
T_2 \; = \; \Exs[y^2 \|\mu_\theta(\xobs, y)\|_2^2] & \leq \sqrt{
  \Exs[y^4] }\sqrt{\Exs[ \|\mu_\theta(\xobs, y)\|_2^4]}.
\end{align*}
Note that $y$ is sub-Gaussian with parameter at most
$\sqrt{\|\thetastar\|_2^2 + \sigma^2}$, whence 
\begin{align*}
\sqrt{\Exs[y^4]} \leq c \; (\|\thetastar\|_2^2 + \sigma^2).
\end{align*}
Similarly, Lemma~\ref{LemMuSubgauss} implies that $\sqrt{\Exs[
    \|\mu_\theta(\xobs, y)\|_2^4]} \leq c \usedim$, and hence $T_2
\leq c' \big (\|\thetastar\|_2^2 + \sigma^2 \big) \usedim$. \\
 
\vspace*{.05in}

Substituting our upper bounds on $T_1$ and $T_2$ into the
decomposition~\eqref{EqnNetArg}, we find that $\sigmaSGD^2 \leq c \,
\big( \|\thetastar\|_2^2 + \sigma^2 \big) \usedim$.  Thus,
Corollary~\ref{CorMissingSGD} follows from
Theorem~\ref{ThmStochasticEM}.


\bibliographystyle{siva}
\bibliography{EM_Journal} 

\begin{thebibliography}{52}
\providecommand{\natexlab}[1]{#1}
\providecommand{\url}[1]{\texttt{#1}}
\providecommand{\urlprefix}{URL }
\providecommand{\eprint}[2][]{\url{#2}}

\bibitem[{Anandkumar et~al.(2012)Anandkumar, Ge, Hsu, Kakade, and
  Telgarsky}]{tensor}
A.~Anandkumar, R.~Ge, D.~Hsu, S.~M. Kakade, and M.~Telgarsky.
\newblock Tensor decompositions for learning latent variable models.
\newblock \emph{CoRR}, abs/1210.7559, 2012.

\bibitem[{Anandkumar et~al.(2013)Anandkumar, Jain, Netrapalli, and
  Tandon}]{Dictionary}
A.~Anandkumar, P.~Jain, P.~Netrapalli, and R.~Tandon.
\newblock Learning sparsely used overcomplete dictionaries via alternating
  minimization.
\newblock Technical report, Microsoft Research, 2013.

\bibitem[{Balan et~al.(2006)Balan, Casazza, and Edidin}]{balan_pr}
R.~Balan, P.~Casazza, and D.~Edidin.
\newblock On signal reconstruction without phase.
\newblock \emph{Applied and Computational Harmonic Analysis}, 20(3):345 -- 356,
  2006.

\bibitem[{Baum et~al.(1970)Baum, Petrie, Soules, and Weiss}]{baum1970}
L.~E. Baum, T.~Petrie, G.~Soules, and N.~Weiss.
\newblock A maximization technique occurring in the statistical analysis of
  probabilistic functions of markov chains.
\newblock \emph{The Annals of Mathematical Statistics}, 41(1):164--171, 1970.

\bibitem[{Beale and Little(1975)}]{beale75}
E.~M.~L. Beale and R.~J.~A. Little.
\newblock Missing values in multivariate analysis.
\newblock \emph{Journal of the Royal Statistical Society. Series B
  (Methodological)}, 37(1):pp. 129--145, 1975.

\bibitem[{Bertsekas(1995)}]{Bertsekas_nonlin}
D.~Bertsekas.
\newblock \emph{Nonlinear Programming}.
\newblock Athena Scientific, 1995.

\bibitem[{Bubeck(2014)}]{bubeck_co}
S.~Bubeck.
\newblock Theory of convex optimization for machine learning.
\newblock 2014.
\newblock \eprint{arXiv:1405.4980}.

\bibitem[{Buldygin and Kozachenko(2000)}]{BulKoz}
V.~V. Buldygin and Y.~V. Kozachenko.
\newblock \emph{Metric characterization of random variables and random
  processes}.
\newblock American Mathematical Society, Providence, RI, 2000.

\bibitem[{Cand{\`e}s et~al.(2013)Cand{\`e}s, Strohmer, and
  Voroninski}]{candes_pr}
E.~J. Cand{\`e}s, T.~Strohmer, and V.~Voroninski.
\newblock Phaselift: Exact and stable signal recovery from magnitude
  measurements via convex programming.
\newblock \emph{Communications on Pure and Applied Mathematics},
  66(8):1241--1274, 2013.

\bibitem[{Celeux et~al.(1995)Celeux, Chauveau, and Diebolt}]{sem1}
G.~Celeux, D.~Chauveau, and J.~Diebolt.
\newblock On stochastic versions of the {EM} algorithm.
\newblock Technical Report 2514, INRIA, 1995.

\bibitem[{Celeux and Govaert(1992)}]{hardem}
G.~Celeux and G.~Govaert.
\newblock A classification {EM} algorithm for clustering and two stochastic
  versions.
\newblock \emph{Comput. Stat. Data Anal.}, 14(3):315--332, 1992.

\bibitem[{Chaganty and Liang(2013)}]{spectral_mor}
A.~T. Chaganty and P.~Liang.
\newblock Spectral experts for estimating mixtures of linear regressions.
\newblock 2013.
\newblock \eprint{arXiv:1306.3729}.

\bibitem[{Chen et~al.(2013)Chen, Yi, and Caramanis}]{convex_mor}
Y.~Chen, X.~Yi, and C.~Caramanis.
\newblock A convex formulation for mixed regression: Near optimal rates in the
  face of noise.
\newblock 2013.
\newblock \eprint{arXiv:1312.7006}.

\bibitem[{Chr{\'e}tien and Hero(2008)}]{hero}
S.~Chr{\'e}tien and A.~O. Hero.
\newblock On {EM} algorithms and their proximal generalizations.
\newblock \emph{ESAIM: Probability and Statistics}, 12:308--326, 2008.

\bibitem[{Dasgupta and Schulman(2007)}]{dasgupta}
S.~Dasgupta and L.~J. Schulman.
\newblock A probabilistic analysis of em for mixtures of separated, spherical
  gaussians.
\newblock \emph{Journal of Machine Learning Research}, 8:203--226, 2007.

\bibitem[{Dempster et~al.(1977)Dempster, Laird, and Rubin}]{dlr}
A.~P. Dempster, N.~M. Laird, and D.~B. Rubin.
\newblock Maximum likelihood from incomplete data via the {EM} algorithm.
\newblock \emph{Journal of the Royal Statistical Society, Series B},
  39(1):1--38, 1977.

\bibitem[{Hartley(1958)}]{hartley}
H.~O. Hartley.
\newblock Maximum likelihood estimation from incomplete data.
\newblock \emph{Biometrics}, 14(2):pp. 174--194, 1958.

\bibitem[{Healy and Westmacott(1956)}]{healy}
M.~Healy and M.~Westmacott.
\newblock Missing values in experiments analysed on automatic computers.
\newblock \emph{Journal of the Royal Statistical Society. Series C (Applied
  Statistics)}, 5(3):pp. 203--206, 1956.

\bibitem[{Hsu and Kakade(2012)}]{hsu_mog}
D.~Hsu and S.~M. Kakade.
\newblock Learning gaussian mixture models: Moment methods and spectral
  decompositions.
\newblock \emph{CoRR}, abs/1206.5766, 2012.

\bibitem[{Iturria et~al.(1999)Iturria, Carroll, and Firth}]{ItuEtal99}
S.~J. Iturria, R.~J. Carroll, and D.~Firth.
\newblock Polynomial regression and estimating functions in the presence of
  multiplicative measurement error.
\newblock \emph{Journal of the Royal Statistical Society Series B - Statistical
  Methodology}, 61:547--561, 1999.

\bibitem[{Jain et~al.(2013)Jain, Netrapalli, and Sanghavi}]{jns13}
P.~Jain, P.~Netrapalli, and S.~Sanghavi.
\newblock Low-rank matrix completion using alternating minimization.
\newblock In \emph{STOC}, pages 665--674. 2013.

\bibitem[{Keshavan et~al.(2010)Keshavan, Montanari, and Oh}]{kmo10}
R.~H. Keshavan, A.~Montanari, and S.~Oh.
\newblock Matrix completion from a few entries.
\newblock \emph{Information Theory, IEEE Transactions on}, 56(6):2980--2998,
  2010.

\bibitem[{Koltchinskii(2011)}]{Koltchinskii}
V.~Koltchinskii.
\newblock \emph{Oracle inequalities in empirical risk minimization and sparse
  recovery problems {\'E}cole d'{\'e}t{\'e} de probabilit{\'e}s de Saint-Flour
  XXXVIII-2008}.
\newblock Springer Verlag, Berlin Heidelberg New York, 2011.

\bibitem[{Ledoux and Talagrand(1991)}]{LedTal91}
M.~Ledoux and M.~Talagrand.
\newblock \emph{Probability in Banach Spaces: Isoperimetry and Processes}.
\newblock Springer-Verlag, New York, NY, 1991.

\bibitem[{Liu and Rubin(1994)}]{ecme}
C.~Liu and D.~B. Rubin.
\newblock The {ECME} algorithm: a simple extension of {EM} and {ECM} with
  faster monotone convergence.
\newblock \emph{Biometrika}, 81:633--648, 1994.

\bibitem[{Loh and Wainwright(2012)}]{loh_isit}
P.-L. Loh and M.~J. Wainwright.
\newblock Corrupted and missing predictors: Minimax bounds for high-dimensional
  linear regression.
\newblock In \emph{ISIT}, pages 2601--2605. 2012.

\bibitem[{Louis(1982)}]{aem2}
T.~A. Louis.
\newblock Finding the observed information matrix when using the {EM}
  algorithm.
\newblock \emph{Journal of the Royal Statistical Society: Series B},
  44:226--233, 1982.

\bibitem[{Ma and Xu(2005)}]{xu_gmm2}
J.~Ma and L.~Xu.
\newblock Asymptotic convergence properties of the {EM} algorithm with respect
  to the overlap in the mixture.
\newblock \emph{Neurocomputing}, 68, 2005.

\bibitem[{McLachlan and Krishnan(2007)}]{em_textbook}
G.~McLachlan and T.~Krishnan.
\newblock \emph{The {EM} Algorithm and Extensions}.
\newblock Wiley Series in Probability and Statistics. Wiley, 2007.

\bibitem[{Meilijson(1989)}]{aem1}
I.~Meilijson.
\newblock A fast improvement of the {EM} algorithm on its own terms.
\newblock \emph{Journal of the Royal Statistical Society: Series B},
  51:127--138, 1989.

\bibitem[{Meng and Rubin(1993)}]{mengecm}
X.~L. Meng and D.~B. Rubin.
\newblock Maximum likelihood via the {ECM} algorithm: a general framework.
\newblock \emph{Biometrika}, 80:267--278, 1993.

\bibitem[{Neal and Hinton(1999)}]{neal99}
R.~M. Neal and G.~E. Hinton.
\newblock A view of the {EM} algorithm that justifies incremental, sparse, and
  other variants.
\newblock In M.~I. Jordan, editor, \emph{Learning in Graphical Models}, pages
  355--368. MIT Press, Cambridge, MA, USA, 1999.

\bibitem[{Nemirovski et~al.(2009)Nemirovski, Juditsky, Lan, and
  Shapiro}]{nemirovski}
A.~Nemirovski, A.~Juditsky, G.~Lan, and A.~Shapiro.
\newblock Robust stochastic approximation approach to stochastic programming.
\newblock \emph{SIAM Journal on Optimization}, 19(4):1574--1609, 2009.

\bibitem[{Nesterov(2004)}]{Nesterov}
Y.~Nesterov.
\newblock \emph{Introductory Lectures on Convex Optimization: A Basic Course}.
\newblock Applied Optimization. Springer, 2004.

\bibitem[{Netrapalli et~al.(2013)Netrapalli, Jain, and Sanghavi}]{altmin_pr}
P.~Netrapalli, P.~Jain, and S.~Sanghavi.
\newblock Phase retrieval using alternating minimization.
\newblock In \emph{NIPS}, pages 2796--2804. 2013.

\bibitem[{Noorshams and Wainwright(2013)}]{NooWai13a}
N.~Noorshams and M.~J. Wainwright.
\newblock Stochastic belief propagation: A low-complexity alternative to the
  sum-product algorithm.
\newblock \emph{IEEE Transactions on Information Theory}, 59(4):1981--2000,
  2013.

\bibitem[{Orchard and Woodbury(1972)}]{orchard1972}
T.~Orchard and M.~A. Woodbury.
\newblock A missing information principle: theory and applications.
\newblock \emph{Proceedings of the Sixth Berkeley Symposium on Mathematical
  Statistics and Probability, Volume 1: Theory of Statistics}, pages 697--715,
  1972.

\bibitem[{Pearson(1894)}]{pearson94}
K.~Pearson.
\newblock \emph{Contributions to the Mathematical Theory of Evolution}.
\newblock Harrison and Sons, 1894.

\bibitem[{Redner and Walker(1984)}]{redner84}
R.~A. Redner and H.~F. Walker.
\newblock Mixture densities, maximum likelihood and the {EM} algorithm.
\newblock \emph{SIAM Review}, 26(2):195--239, 1984.

\bibitem[{Rubin(1974)}]{rubin74}
D.~B. Rubin.
\newblock Characterizing the estimation of parameters in incomplete-data
  problems.
\newblock \emph{Journal of the American Statistical Association}, 69(346):pp.
  467--474, 1974.

\bibitem[{Sundberg(1974)}]{sundberg}
R.~Sundberg.
\newblock Maximum likelihood theory for incomplete data from an exponential
  family.
\newblock \emph{Scand. J. Statist}, 1:49--58, 1974.

\bibitem[{Tanner and Wong(1987)}]{sem2}
M.~A. Tanner and W.~H. Wong.
\newblock The calculation of posterior distributions by data augmentation.
\newblock \emph{Journal of the American Statistical Association}, 82:528--550,
  1987.

\bibitem[{Tseng(2004)}]{tseng}
P.~Tseng.
\newblock An analysis of the em algorithm and entropy-like proximal point
  methods.
\newblock \emph{Mathematics of Operations Research}, 29(1):pp. 27--44, 2004.

\bibitem[{van~de Geer(2000)}]{vandeGeer}
S.~van~de Geer.
\newblock \emph{Empirical Processes in M-Estimation}.
\newblock Cambridge University Press, 2000.

\bibitem[{van~der Vaart and Wellner(1996)}]{vanderVaart96}
A.~W. van~der Vaart and J.~Wellner.
\newblock \emph{Weak Convergence and Empirical Processes}.
\newblock Springer-Verlag, New York, NY, 1996.

\bibitem[{van Dyk and Meng(2000)}]{sem4}
D.~A. van Dyk and X.~L. Meng.
\newblock Algorithms based on data augmentation: A graphical representation and
  comparison.
\newblock In \emph{Computing Science and Statistics: Proceedings of the 31st
  Symposium on the Interface}, pages 230--239. Berk and M. Pourahmadi, 2000.

\bibitem[{Vershynin(2010)}]{vershynin_nonasymp}
R.~Vershynin.
\newblock Introduction to the non-asymptotic analysis of random matrices.
\newblock Chapter 5 of: Compressed Sensing, Theory and Applications. Edited by
  Y. Eldar and G. Kutyniok. Cambridge University Press, 2012, 2010.
\newblock \eprint{arXiv:1011.3027}.

\bibitem[{Wei and Tanner(1990)}]{sem3}
G.~Wei and M.~A. Tanner.
\newblock A {M}onte {C}arlo implementation of the {EM} algorithm and the poor
  man's data augmentation algorithm.
\newblock \emph{Journal of the American Statistical Association}, 85:699--704,
  1990.

\bibitem[{Wu(1983)}]{wu1983}
C.~F.~J. Wu.
\newblock On the convergence properties of the {EM} algorithm.
\newblock \emph{The Annals of Statistics}, 11(1):95--103, 1983.

\bibitem[{Xu and Jordan(1996)}]{xu_gmm1}
L.~Xu and M.~I. Jordan.
\newblock On convergence properties of the {EM} algorithm for gaussian
  mixtures.
\newblock \emph{Neural Comput.}, 8(1):129--151, 1996.

\bibitem[{Xu and You(2007)}]{XuYou07}
Q.~Xu and J.~You.
\newblock Covariate selection for linear errors-in-variables regression models.
\newblock \emph{Communications in Statistics - Theory and Methods},
  36(2):375--386, 2007.

\bibitem[{Yi et~al.(2013)Yi, Caramanis, and Sanghavi}]{hardem_mor}
X.~Yi, C.~Caramanis, and S.~Sanghavi.
\newblock Alternating minimization for mixed linear regression.
\newblock 2013.
\newblock \eprint{arXiv:1310.3745}.

\end{thebibliography}


\end{document}